\documentclass[11pt]{article}

\usepackage[ngerman,english]{babel}
\usepackage[utf8]{inputenc}
\usepackage[T1]{fontenc} 
\usepackage{ae,aecompl} 

\usepackage{url}

\usepackage{amsmath} 
\usepackage{amsfonts}
\usepackage{amsbsy}
\usepackage{amsthm}
\usepackage{amssymb}

\usepackage[a4paper]{geometry}

\usepackage{cancel} 

\usepackage[mathcal]{eucal}
\usepackage{mathrsfs} 

\usepackage{multicol} 
\usepackage{pstricks}
\usepackage{graphicx}
\usepackage{color}

\usepackage{float}
\usepackage{subfig}

\usepackage{enumerate}

\newtheorem{thm}{Theorem}[section]
\newtheorem{lem}[thm]{Lemma}

\theoremstyle{definition}
\newtheorem{rem}[thm]{Remark}

\newenvironment{prf}{\noindent\textbf{Proof}.}
{\hfill $\Box$}

\DeclareMathOperator{\ddd}{\ensuremath \normalfont{D}}

\newcommand{\oox}{\tilde{o}^{(1)}}
\newcommand{\ooy}{\tilde{o}^{(2)}}
\newcommand{\ooz}{\tilde{o}^{(3)}}

\newcommand{\oy}{{o}^{(2)}}
\newcommand{\oz}{{o}^{(3)}}

\newsavebox{\psymbol}

\savebox{\psymbol}(0,0){

  \put(0,0){\circle{1}}}

\newcommand{\Theorem}{\begin{theorem}}
\newcommand{\CTheorem}{\begin{theorem}}
\newcommand{\EndTheorem}{\end{theorem}}
\newcommand{\Proof}{\begin{proof} {}}
\newcommand{\EndProof}{\end{proof}}
\newcommand{\Lemma}{\begin{lemma}}
\newcommand{\CLemma}{\begin{lemma}}
\newcommand{\EndLemma}{\end{lemma}}
\newcommand{\Proposition}{\begin{proposition}}
\newcommand{\EndProposition}{\end{proposition}}
\newcommand{\Corollary}{\begin{corollary}}
\newcommand{\CCorollary}{\begin{corollary}}
\newcommand{\EndCorollary}{\end{corollary}}
\newcommand{\Remark}{\begin{remark}}
\newcommand{\EndRemark}{\end{remark}}
\newcommand{\Definition}{\begin{definition}}
\newcommand{\CDefinition}{\begin{definition}}
\newcommand{\EndDefinition}{\end{definition}}
\newcommand{\Example}{\begin{example}}
\newcommand{\EndExample}{\end{example}}
\newcommand{\Algorithm}{\begin{algorithm}}
\newcommand{\EndAlgorithm}{\end{algorithm}}
\newcommand{\Exercise}{\begin{exercise}}
\newcommand{\EndExercise}{\end{exercise}}

\newcommand{\be}{\begin{displaymath}}
\newcommand{\ee}{\end{displaymath}}
\newcommand{\ben}{\begin{equation}}
\newcommand{\een}{\end{equation}}
\newcommand{\bean}{\begin{eqnarray}}
\newcommand{\eean}{\end{eqnarray}}
\newcommand{\bea}{\begin{eqnarray*}}
\newcommand{\eea}{\end{eqnarray*}}

\newcommand{\R}{{\Bbb{R}}}

\newcommand{\NN}{{\Bbb N}}
\newcommand{\No}{{\NN_0}}

\def\={ \ =\  }

\def\be{\beta}

\def\0xe{(\partial_0;\xi ,\eta)}

\def\F{{\cal F}}

\def\U{{\cal U}}

\def\L{{\rm L}}

\newcommand{\ltwo}[1]{\mathcal{L}^2(\mathbb{R})}

\newcommand{\rtwo}[1]{\mathbb{R}^2}

\newfont{\fontfuerpfeil}{cmsy10 scaled 3583}
\newfont{\fontfuerpfeilklein}{cmsy10 scaled 2500}
\definecolor{dGray}{gray}{0.89}
\definecolor{lGray}{gray}{0.91}
\definecolor{llGray}{gray}{0.95}

\newfont{\astfont}{cmr10 scaled 2074}

\def\F{{\cal F}}

\begin{document}

\begin{center}
\Large{\textbf{Spherical Potential Theory: Tools and Applications}}\normalsize
\\[3ex]
{Christian Gerhards}\footnote{Computational Science Center, University of Vienna, 1090 Vienna, Autria, e-mail: christian.gerhards@univie.ac.at}
\\[5ex]
\end{center}

Classical potential theoretic concepts in the Euclidean space $\mathbb{R}^3$ have been described in an earlier chapter of this handbook. They appear frequently in geodesy when treating the harmonic gravitational potential in the exterior of the spherical Earth. The sphere $\Omega_R=\{x\in\mathbb{R}^3:|x|=R\}$ occurs as the boundary surface of a subdomain in $\mathbb{R}^3$. Opposed to this, in the present chapter, the sphere is not representing a boundary surface, it is rather regarded as the underlying domain on which a problem is formulated.  Examples for this are the spherical Navier-Stokes equations and shallow water equations in meteorology and ocean modeling (see, e.g.,  \cite{comblen09, fenglerfreeden05, ganesh11, ilin91, pedlosky79, stewart}). But also simpler spherical differential equations occur in geodesy and geomagnetism (see, e.g., \cite{backus96, fehlinger07, fehlinger09, freeden09, freeden05, gerhards14c, heiskanenmoritz67}) and vortex dynamics (see, e.g., \cite{kidambinewton00, kimuraokamoto87, pedlosky79, polvanidritschel93, 
stewart}), more precisely, 
those based on the Beltrami operator $\Delta^*$ (the spherical 
counterpart to the Laplace operator $\Delta$). Latter is going to be the focus of this chapter. In particular, we are interested in the Beltrami equation on subdomains $\Gamma_R\subset\Omega_R$ of the sphere, which eventually leads to potential theoretic concepts analogous to those of the Euclidean case. Subdomains appear naturally, e.g., due to only regionally available data or coastal/continental boundaries. The problems we take a closer look at are the following (note that $\Omega$ and $\Gamma$ are simply abbreviations for the unit sphere $\Omega_1$ and a corresponding subdomain $\Gamma_1$):
\begin{description}
\item[\textit{Poisson Problem} {\rm (PP):}]  Let $H$ be of class ${\rm C}^{(1)}(\overline{\Gamma})$. We are looking for a function $U$ of class ${\rm C}^{(2)}(\Gamma)$ such that
\begin{align}\label{eqn:beltramiequation}
\Delta^*U(\xi)&=H(\xi),\quad\xi\in\Gamma.
\end{align}
 \item[\textit{Dirichlet Problem} {\rm (DP):}]  Let $F$ be of class ${\rm C}^{(0)}(\partial\Gamma)$. We are looking for a function $U$ of class ${\rm C}^{(2)}(\Gamma)\cap {\rm C}^{(0)}(\overline{\Gamma})$ such that
\begin{align}
\Delta^*U(\xi)&=0,\qquad\,\,\,\xi\in\Gamma,
\\U^-(\xi)&=F(\xi), \quad\xi\in\partial\Gamma.\label{eqn:dv1}
\end{align}
\item[\textit{Neumann Problem} {\rm (NP):}]  Let $F$ be of class ${\rm C}^{(0)}(\partial\Gamma)$. We are looking for a function $U$ of class ${\rm C}^{(2)}(\Gamma)\cap {\rm C}^{(0)}(\overline{\Gamma})$, with a well-defined normal derivative $\frac{\partial}{\partial\nu}U^-$ on $\partial\Gamma$, such that
\begin{align}
\Delta^*U(\xi)&=0,\qquad\,\,\,\xi\in\Gamma,
\\\frac{\partial}{\partial\nu}U^-(\xi)&=F(\xi), \quad \xi\in\partial\Gamma.\label{eqn:nv1}
\end{align}
\end{description}
In the setting above, $\partial\Gamma$ denotes the boundary curve of $\Gamma$, $\nu(\xi)$ the outward directed unit normal vector at $\xi\in\partial\Gamma$, and $\frac{\partial}{\partial\nu}$ the corresponding normal derivative. The minus of $U^-$ simply indicates that we are approaching the boundary $\partial\Gamma$ from within $\Gamma$.

Certainly, the problems above and its potential theoretic consequences can be and have been treated on more general manifolds than the sphere (e.g., in \cite{duduchava06, mitreataylor99,mitreataylor00}). However, we focus on the geophysically relevant case of the sphere where explicit representations of the fundamental solution and some Green's functions are known. In large parts, we follow the course of \cite{freedengerhards12} and emphasize similarities and differences to the Euclidean case.

The first section supplies the reader with necessary notations and several mathematical tools related to spherical potential theory. In Section \ref{sec:spot2}, we treat the problems (PP), (DP), and (NP). In particular, we are interested in integral representations of their solutions. In Section \ref{sec:sdecomp}, we turn towards spherical differential operators of order one, namely the surface gradient $\nabla^*$ (the spherical counterpart to the gradient $\nabla$) and the surface curl gradient $\L^*$. We investigate the corresponding differential equations on $\Gamma$ as well as the so-called spherical Helmholtz decomposition and the spherical Hardy-Hodge decomposition. Section \ref{sec:cfs} comments briefly on complete function systems and approximation methods on the sphere. Finally, in Section \ref{sec:app}, applications of the previous concepts to some geophysical problems are discussed, namely, vertical deflections, (geostrophic) ocean flow, and a toy problem for point vortex motion.

\section{Fundamental Tools}\label{sec:spot0}

Of fundamental importance to us is the Beltrami operator $\Delta^*$ which denotes the tangential contribution to the Euclidean Laplace operator $\Delta$. More precisely,
\begin{align}
 \Delta_x=\frac{1}{r^2}\frac{\partial}{\partial r}r^2\frac{\partial}{\partial r}+\frac{1}{r^2}\Delta^*_\xi,\label{eqn:deltabeltramirel}
\end{align}
where $\Delta_x$ acts on $x\in\mathbb{R}^3$ while $\Delta^*_\xi$ acts on $\xi=\frac{x}{|x|}\in\Omega$. The length $|x|$ is usually denoted by $r$. Furthermore, $\nabla^*$ stands for the (spherical) surface gradient, which denotes the tangential contribution to the gradient $\nabla$:
\begin{align}
 \nabla_x=\xi\frac{\partial}{\partial r}+\frac{1}{r}\nabla^*_\xi.
\end{align}
The occasionally occuring (spherical) surface curl gradient $\L^*$ acts via $\xi\wedge \nabla^*_\xi$ at a point $\xi\in\Omega$ (\lq\lq$\wedge$'' denotes the vector product). It should be noted that $\Delta^*=\nabla^*\cdot\nabla^*=\L^*\cdot \L^*$ (\lq\lq$\,\cdot\,$'' denotes the Euclidean inner product). If it is clear on which variables the operators act, we usually omit the subindices $\xi$ and $x$.  For convenience, we typically use greek letters $\xi,\eta$ to indicate unit vectors in $\Omega$ while we use latin letters $x,y$ for general vectors in $\mathbb{R}^3$. Upper case letters $F,G$ denote scalar-valued functions mapping $\Gamma\subset \Omega$ into $\mathbb{R}$ while lower case letters $f,g$ denote vector-valued functions mapping $\Gamma\subset \Omega$ into $\mathbb{R}^3$. Correspondingly, the set of $k$-times continuously differentiable scalar-valued functions on $\Gamma$ is designated by ${\rm C}^{(k)}(\Gamma)$ and the set of $k$-times continuously differentiable vector-valued functions on $\Gamma$  by 
$c^{(k)}(\Gamma)$. The closure of $\Gamma$ is denoted by $\overline{\Gamma}$ and 
the open complement by $\Gamma^c=\Omega\setminus\overline{\Gamma}$.

\begin{figure}
\begin{center}
\includegraphics[scale=0.6]{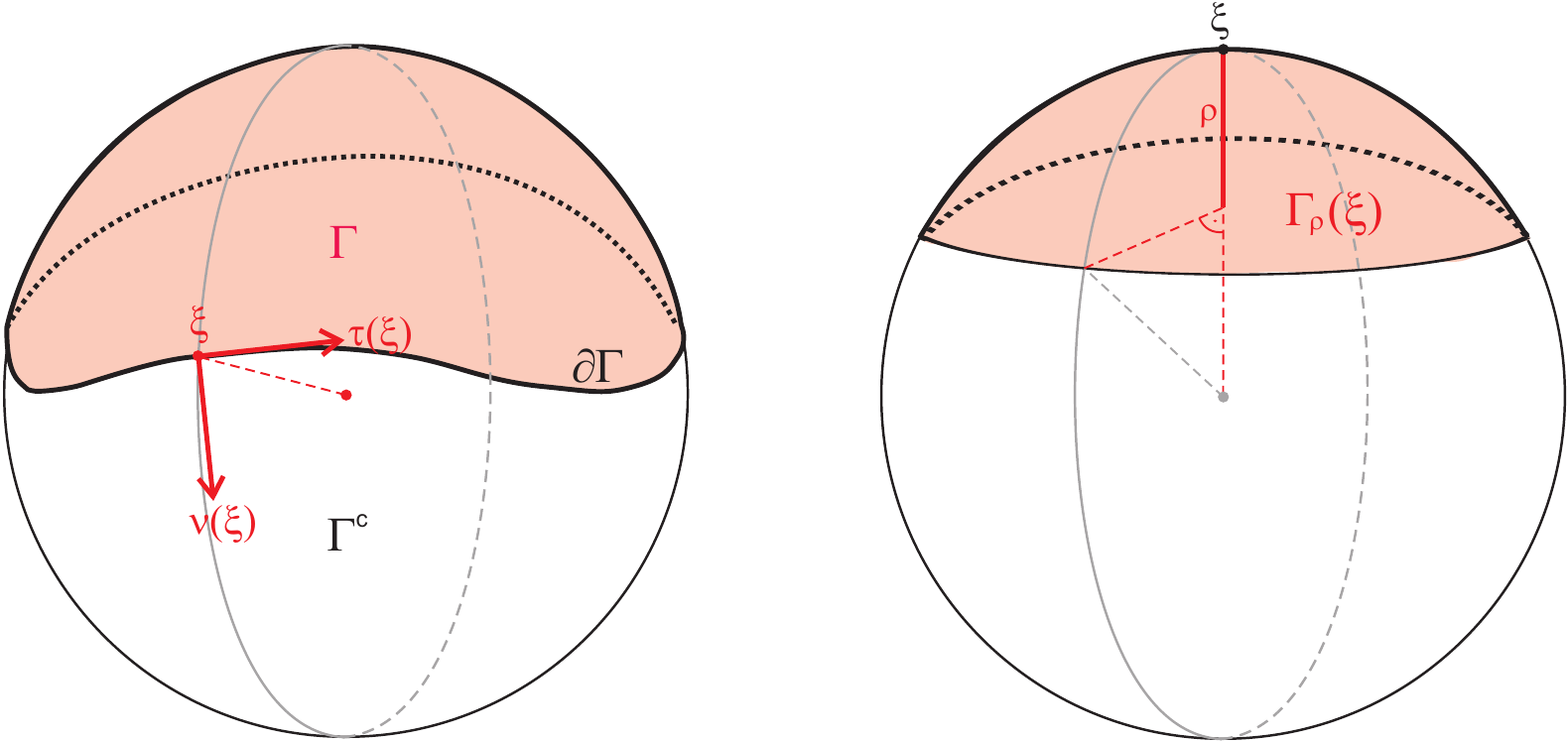}
\end{center}
\caption{Examples for a general regular region $\Gamma$ (left) and a spherical cap $\Gamma_\rho(\xi)$ with center $\xi$ and radius $\rho$ (right).}
\label{fig:loctancoords}
\end{figure}

Whenever we talk about subdomains $\Gamma\subset\Omega$ in this chapter, we mean, without further mention, regular regions, i.e., subdomains with a sufficiently smooth boundary curve $\partial\Gamma$ (for details, the reader is referred to \cite{freedengerhards12}; an exemplary illustration is supplied in Figure \ref{fig:loctancoords}). For such regular regions, the positively oriented unit tangential vector $\tau(\xi)$ at a point $\xi\in\partial\Gamma$ is well-defined. The unit normal vector $\nu(\xi)$  at $\xi\in\partial\Gamma$ points into the exterior of $\Gamma$ and is perpendicular to $\tau(\xi)$ and $\xi$ (i.e., $\nu(\xi)$ is perpendicular to the  boundary curve $\partial\Gamma$ but tangential to the unit sphere $\Omega$). The normal derivative  of a scalar-valued function $F$ at $\xi\in\partial\Gamma$ is defined as
\begin{align}
\frac{\partial}{\partial\nu}F(\xi)=\nu(\xi)\cdot\nabla^*_\xi F(\xi).
\end{align}

\subsection{Green's Formulas}\label{sec:gf}

We frequently need integral expressions that describe the shifting of differential operators from one integrand to another, so-called Green's formulas. Some spherical versions are stated in the next theorem.

\begin{thm}[Spherical Green Formulas I]\textbf{}\\[-4ex]\label{thm:greenformulas}
\begin{enumerate}[(a)]
\item If $f$ is of class ${\rm c}^{(1)}(\overline{\Gamma})$ and tangential, i.e., $\xi\cdot f(\xi)=0$ for $\xi\in\overline{\Gamma}$, then
\begin{align}
\int_\Gamma\nabla^*\cdot f(\eta)d\omega(\eta)&=\int_{\partial\Gamma}\nu(\eta)\cdot f(\eta)d\sigma(\eta),\label{eqn:green0a}
\\\int_\Gamma \L^*\cdot f(\eta)d\omega(\eta)&=\int_{\partial\Gamma}\tau(\eta)\cdot f(\eta)d\sigma(\eta).\label{eqn:green0b}
\end{align}
\item If $F$ is  of class ${\rm C}^{(1)}(\overline{\Gamma})$ and $f$ of class ${\rm c}^{(1)}(\overline{\Gamma})$, then
\begin{align}
&\int_\Gamma f(\eta)\cdot \nabla^*F(\eta)d\omega(\eta)+\int_\Gamma F(\eta)\nabla^*\cdot f(\eta) d\omega(\eta)\label{eqn:green1}
\\&=\int_{\partial\Gamma}\nu(\eta)\cdot \left(F(\eta)f(\eta)\right)d\sigma(\eta)+2\int_{\Gamma}\eta\cdot \left(F(\eta)f(\eta)\right)d\omega(\eta),\nonumber
\\[1.25ex]&\int_\Gamma f(\eta)\cdot \L^*F(\eta)d\omega(\eta)+\int_\Gamma F(\eta)\L^*\cdot f(\eta) d\omega(\eta)
\\&=\int_{\partial\Gamma}\tau(\eta)\cdot \left(F(\eta)f(\eta)\right)d\sigma(\eta).\nonumber
\end{align}
\item If $F,H$ are functions of class ${\rm C}^{(2)}(\overline{\Gamma})$, then
\begin{align}\label{eqn:green3}
&\int_\Gamma F(\eta)\Delta^*H(\eta)d\omega(\eta)-\int_\Gamma H(\eta)\Delta^*F(\eta) d\omega(\eta)
\\&=\int_{\partial\Gamma}F(\eta)\frac{\partial}{\partial\nu}H(\eta)d\sigma(\eta)-\int_{\partial\Gamma}H(\eta)\frac{\partial}{\partial\nu}F(\eta)d\sigma(\eta).\nonumber
\end{align}
\end{enumerate}
Generally, '$d\omega$' denotes the surface element in $\Gamma\subset\Omega$ and '$d\sigma$' the line element on $\partial\Gamma$.
\end{thm}

\begin{rem}
The formulas \eqref{eqn:green1}--\eqref{eqn:green3} are direct consequences of \eqref{eqn:green0a} and \eqref{eqn:green0b}. Dropping the boundary terms $\int_{\partial\Gamma}\ldots d\sigma$, all of these formulas also hold true for the choice $\Gamma=\Omega$.
\end{rem}

A crucial step for later considerations is the combination of Green's formulas with the \emph{fundamental solution for the Beltrami operator} $G(\Delta^*;\cdot):[-1,1)\to\R$, which is uniquely determined by the following properties:
\begin{enumerate}
\item For any fixed $\xi \in \Omega$, the function $\eta\mapsto G(\Delta^*;\xi\cdot\eta)$ is twice continuously differentiable on $\Omega\setminus\{\xi\}$ and
\begin{align}\label{eqn:fsolprop}
\Delta^*_\eta G(\Delta^*;\xi\cdot\eta)= -\frac{1}{4 \pi},\quad\eta\in\Omega\setminus\{\xi\}.
\end{align}

\item For any fixed $\xi \in \Omega$, the function $\eta\mapsto G(\Delta^*;\xi\cdot\eta)- \frac{1}{4 \pi}\ln(1-\xi \cdot \eta)$ is continuously differentiable on $\Omega$.

\item For any fixed $\xi \in \Omega$, it holds $\frac{1}{4\pi}\int_\Omega G(\Delta^*;\xi\cdot\eta) d\omega(\eta) = 0$.
\end{enumerate}
Some basic calculations show that the function given by
\begin{align}
G(\Delta^*;t)=\frac{1}{4\pi}\ln(1-t)+\frac{1}{4\pi}(1-\ln(2)),\quad t\in[-1,1),\label{eqn:repfundsol}
\end{align}
satisfies the properties (i)--(iii). The property (i) denotes the major difference between the fundamental solution for the Laplace operator $G(\Delta;\cdot)$ and its spherical counterpart. While $G(\Delta;\cdot)$ generates a 'true' Dirac distribution in the sense that $\Delta_yG(\Delta;|x-y|)=0$, $y\in\mathbb{R}^3\setminus\{x\}$, the fundamental solution $G(\Delta^*;\cdot)$ only generates a Dirac distribution up to an additive constant (reflecting the nullspace of the Beltrami operator $\Delta^*$). Eventually, applying Green's formulas from Theorem \ref{thm:greenformulas}, the properties of $G(\Delta^*;\cdot)$ lead to the following integral representations.

\begin{thm}[Spherical Green Formulas II]\textbf{}\\[-4ex]\label{thm:locfundthmd}
\begin{enumerate}[(a)]
\item If $F$ is of class ${\rm C}^{(2)}(\overline{\Gamma})$, then we have for $\xi\in\Omega$, 
\begin{align}
\frac{\alpha(\xi)}{2\pi}F(\xi)=&\frac{1}{4\pi}\int_\Gamma F(\eta)d\omega(\eta)+\int_\Gamma G(\Delta^*;\xi\cdot\eta)\Delta^*_\eta F(\eta)d\omega(\eta)\label{eqn:3greensphere}
\\&+\int_{\partial\Gamma} F(\eta)\frac{\partial}{\partial\nu(\eta)}G(\Delta^*;\xi\cdot\eta) d\sigma(\eta)\nonumber
-\int_{\partial\Gamma}G(\Delta^*;\xi\cdot\eta)\frac{\partial}{\partial\nu(\eta)}F(\eta)d\sigma(\eta).\nonumber
\end{align}
\item If $F$ is of class ${\rm C}^{(1)}(\overline{\Gamma})$, then we have for $\xi\in\Omega$,
\begin{align}
\frac{\alpha(\xi)}{2\pi}F(\xi)=&\frac{1}{4\pi}\int_\Gamma F(\eta)d\omega(\eta)-\int_\Gamma \nabla^*_\eta G(\Delta^*;\xi\cdot\eta)\cdot \nabla^*_\eta F(\eta)d\omega(\eta)\qquad
\\&+\int_{\partial\Gamma} F(\eta)\frac{\partial}{\partial\nu(\eta)}G(\Delta^*;\xi\cdot\eta)d\sigma(\eta),\nonumber
\\=&\frac{1}{4\pi}\int_\Gamma F(\eta)d\omega(\eta)-\int_\Gamma \L^*_\eta G(\Delta^*;\xi\cdot\eta)\cdot \L^*_\eta F(\eta)d\omega(\eta)\qquad\nonumber
\\&+\int_{\partial\Gamma} F(\eta)\frac{\partial}{\partial\nu(\eta)}G(\Delta^*;\xi\cdot\eta)d\sigma(\eta).\nonumber
\end{align}
\end{enumerate}
The solid angle $\alpha$ of a regular region $\Gamma$ is defined such that $\alpha(\xi)=2\pi$ for $\xi\in\Gamma$, $\alpha(\xi)=\pi$ for $\xi\in\partial\Gamma$, and $\alpha(\xi)=0$ for $\xi\in\Gamma^c$.
\end{thm}

\begin{rem}
Again, dropping the boundary terms $\int_{\partial\Gamma}\ldots d\sigma$ in the expressions of Theorem \ref{thm:locfundthmd} leads to results that hold true for the global choice $\Gamma=\Omega$. 
\end{rem}

\begin{rem}
Choosing $F\equiv 1$ in any of the formulas in Theorem \ref{thm:locfundthmd} implies
\begin{align}\label{eqn:intnormderivgf1}
\int_{\partial\Gamma}\frac{\partial}{\partial\nu(\eta)}G(\Delta^*;\xi\cdot\eta)d\sigma(\eta)=\left\{\begin{array}{ll}1-\frac{\|\Gamma\|}{4\pi},&\xi\in\Gamma,
\\\frac{1}{2}-\frac{\|\Gamma\|}{4\pi},&\xi\in\partial\Gamma,
\\-\frac{\|\Gamma\|}{4\pi},&\xi\in\Gamma^c,                                                               
\end{array}\right.
\end{align}
where $\|\Gamma\|$ denotes the surface area of $\Gamma$. The behaviour of \eqref{eqn:intnormderivgf1} across the boundary $\partial\Gamma$ states a first hint at the limit and jump relations of the layer potentials in Section \ref{sec:cp}. Apart from the additive constant $\frac{\|\Gamma\|}{4\pi}$, they are identical to the Euclidean setting.
\end{rem}

\subsection{Harmonic Functions}\label{sec:hf}

In this subsection, we turn towards functions that are  \emph{harmonic (with respect to Beltrami operator)} in $\Gamma$, i.e., functions $U$ of class ${\rm C}^{(2)}(\Gamma)$ that satisfy
\begin{align}\label{eqn:harmeq}
 \Delta^*U=0\quad\textnormal{in }\Gamma.
\end{align}
If no confusion with the Euclidean case is likely to arise, we just say that $U$ is harmonic. Plugging such functions into \eqref{eqn:3greensphere}, together with the choice of $\Gamma$ being a spherical cap $\Gamma_\rho(\xi)=\{\eta\in\Omega:1-\xi\cdot\eta<\rho\}$ with center $\xi\in\Omega$ and radius $\rho\in(0,2)$, we end up with the mean value property for harmonic functions.

\begin{thm}[Mean Value Property I]\label{thm:meanvalpropsphere}
A function $U$ of class ${\rm C}^{(0)}(\Gamma)$ is harmonic if and only if
\begin{align}
U(\xi)=\frac{1}{4\pi}\int_{\Gamma_\rho(\xi)}U(\eta)d\omega(\eta)+\frac{\sqrt{2-\rho}}{4\pi\sqrt{\rho}}\int_{\partial\Gamma_\rho(\xi)}U(\eta)d\sigma(\eta),\quad\xi\in\Gamma,\label{eqn:mvtsphere}
\end{align}
for any spherical cap $\overline{\Gamma_\rho(\xi)}\subset\Gamma$.
\end{thm}

The mean value property above contains the typical additive term for spherical problems. However, we can get rid of this additive constant when using Green's functions for spherical caps as described later on in Section \ref{sec:gf2}. We are led to the following representation which resembles a Mean Value Property that is more closely related to the Euclidean case of functions that are harmonic with respect to the Laplace operator.

\begin{thm}[Mean Value Property II]\label{thm:meanvalpropsphere2}
A function $U$ of class ${\rm C}^{(0)}(\Gamma)$ is harmonic if and only if
\begin{align}
U(\xi)=\frac{1}{2\pi\sqrt{\rho(2-\rho)}}\int_{\partial\Gamma_\rho(\xi)}U(\eta)d\sigma(\eta),\quad\xi\in\Gamma,\label{eqn:mvtsphere3}
\end{align}
for any spherical cap $\overline{\Gamma_\rho(\xi)}\subset\Gamma$.
\end{thm}

Once a Mean Value Property is established, it can be used to derive a Maximum Principle. For details, we refer to \cite{freedengerhards12} or, in the Euclidean case, any book on classical potential theory such as \cite{freedenmichel04a, helms69,kellogg67,wermer74}. 

\begin{thm}[Maximum Principle]\label{thm:maxprinciple}
 If $U$ of class ${\rm C}^{(2)}({\Gamma})\cap {\rm C}^{(0)}(\overline{\Gamma})$ is harmonic, then
\begin{align}
\sup_{\xi\in\Gamma}|U(\xi)|\leq \sup_{\xi\in\partial\Gamma}|U(\xi)|.
\end{align}
\end{thm}

\subsection{Surface Potentials}\label{sec:sp}

Analogous to the Euclidean setting, we can define a Newton potential and layer potentials for the spherical setting which take over the corresponding roles. The obvious difference is that now the Newton potential is a surface potential and the layer potentials represent curve potentials. Throughout this section, we take a closer look at the surface potential
\begin{align}
U(\xi)=\int_\Gamma G(\Delta^*;\xi\cdot\eta) H(\eta)d\omega(\eta),\quad \xi\in\Omega, \label{eqn:newtonpotsphere2}
\end{align}
From the properties of the fundamental solution for the Beltrami operator it becomes directly clear that $U$ is of class ${\rm C}^{(2)}(\Gamma^c)$ and that
\begin{align}\label{eqn:udeltaappgc}
\Delta^*_\xi\int_\Gamma G(\Delta^*;\xi\cdot\eta) H(\eta)d\omega(\eta)=-\frac{1}{4\pi}\int_\Gamma H(\eta)d\omega(\eta),\quad \xi\in\Gamma^c.
\end{align}
Yet, the interesting question is what happens if $\xi\in{\Gamma}$, i.e., when the integration region contains the singularity of $G(\Delta^*;\cdot)$.

\begin{thm}\label{thm:nablaintexchange}
If $H$ is of class ${\rm C}^{(0)}(\overline{\Gamma})$ and $U$ is given by \eqref{eqn:newtonpotsphere2}, then $U$ is of class ${\rm C}^{(1)}(\Omega)$ and
\begin{align}
\nabla^*_\xi\int_\Gamma G(\Delta^*;\xi\cdot\eta) H(\eta)d\omega(\eta)=\int_\Gamma  \nabla^*_\xi G(\Delta^*;\xi\cdot\eta)H(\eta) d\omega(\eta),\quad \xi\in\Omega.\label{eqn:interchangenablaandint}
\end{align}
\end{thm}

The proof of the theorem above can be based on a regularization of the fundamental solution $G(\Delta^*;\cdot)$. This approach also works for the application of the Beltrami operator to $U$. However, in connection with Theorem \ref{thm:locfundthmd}, we find that $\Delta^*U$ is not continuous across $\partial\Gamma$ anymore (compare equation \eqref{eqn:udeltaappgc} and Theorem \ref{thm:poissoneqsphere}). For brevity, we do not supply the proofs at this point but refer the reader, e.g., to \cite{freedengerhards12}. A related regularized Green function plays an important role in the applications in Section \ref{sec:app} and is explained in more detail later on.
 
\begin{thm}\label{thm:poissoneqsphere}
If $H$ is  of class ${\rm C}^{(1)}(\overline{\Gamma})$ and $U$ is given by \eqref{eqn:newtonpotsphere2}, then $U$ is of class ${\rm C}^{(2)}(\Gamma)$ and satisfies
\begin{align}
\Delta^*_\xi\int_\Gamma G(\Delta^*;\xi\cdot\eta) H(\eta)d\omega(\eta)=H(\xi)-\frac{1}{4\pi}\int_\Gamma H(\eta)d\omega(\eta),\quad \xi\in\Gamma.
\end{align}
\end{thm}

\subsection{Curve Potentials}\label{sec:cp}

While surface potentials are useful to deal with the Poisson problem (PP), curve potentials are particularly useful when dealing with functions that are harmonic (with respect to the Beltrami operator). More precisely, we take a closer look at the two layer potentials
\begin{align}\label{eqn:singlay}
 U_1[\tilde{Q}](\xi)=\int_{\partial\Gamma}G(\Delta^*;\xi\cdot \eta)\tilde{Q}(\eta)\,d\sigma(\eta),\quad\xi\in\Gamma,
\end{align}
and 
\begin{align}\label{eqn:doublay}
 U_2[Q](\xi)=\int_{\partial\Gamma}\left(\frac{\partial}{\partial\nu(\eta)}G(\Delta^*;\xi\cdot \eta)\right)Q(\eta)\,d\sigma(\eta),\quad\xi\in\Gamma.
\end{align}
From the properties of the fundamental solution $G(\Delta^*;\cdot)$ it can be seen that the so-called \emph{double-layer potential} $U_2[Q]$ is harmonic in $\Gamma$ for any $Q$ of class ${\rm C}^{(0)}(\partial\Gamma)$. The \emph{single-layer potential} $U_1[\tilde{Q}]$ is harmonic in $\Gamma$ if $\tilde{Q}$ is of class ${\rm C}^{(0)}(\partial\Gamma)$ and if the integral over $\partial\Gamma$ vanishes, i.e., if $\int_{\partial\Gamma}\tilde{Q}(\eta)d\sigma(\eta)=0$ (we say that $\tilde{Q}$ is of class ${\rm C}_0^{(0)}(\partial\Gamma)$). Therefore, these two potentials represent good candidates for solutions to the  boundary value problems (DP) and (NP). The aim of the present section is to investigate the behaviour of the single- and double-layer potentials $U_1[Q]$ and $U_2[\tilde{Q}]$, respectively, when they approach the boundary $\partial\Gamma$. 
The essential behaviour of the double-layer potential $U_2[Q]$ is already reflected by the relation \eqref{eqn:intnormderivgf1}. Based on this relation and a set of several more technical estimates, one can prove the following set of limit- and jump-relations at the boundary $\partial\Gamma$.

\begin{thm}[Limit- and Jump-Relations]\label{thm:limitrelc}
Let $Q$, $\tilde{Q}$ be of class ${\rm C}^{(0)}(\partial\Gamma)$ and $U_1$, $U_2$ be given as in \eqref{eqn:singlay} and \eqref{eqn:doublay}, respectively. Furthermore, let $\xi\in\partial\Gamma$.
\begin{enumerate}[(a)]
\item For the single-layer potential, we have the {limit-relations}
\begin{align}
\lim_{\tau\to0+}U_1[\tilde{Q}]\left(\frac{\xi\pm\tau\nu(\xi)}{\sqrt{1+\tau^2}}\right)-U_1[\tilde{Q}](\xi)&=0,
\\\lim_{\tau\to0+}\left(\frac{\partial}{\partial\nu}U_1[\tilde{Q}]\right)\left(\frac{\xi\pm\tau\nu(\xi)}{\sqrt{1+\tau^2}}\right)-\left(\frac{\partial}{\partial\nu}U_1[\tilde{Q}]\right)(\xi)&=\pm \frac{1}{2}\tilde{Q}(\xi).
\end{align}
For the double-layer potential, we have
\begin{align}
\lim_{\tau\to0+}U_2[Q]\left(\frac{\xi\pm\tau\nu(\xi)}{\sqrt{1+\tau^2}}\right)-U_2[Q](\xi)&=\mp \frac{1}{2}Q(\xi).
\end{align}
\item For the single-layer potential, we have the {jump-relations}
\begin{align}
\lim_{\tau\to0+}\left(U_1[\tilde{Q}]\left(\frac{\xi+\tau\nu(\xi)}{\sqrt{1+\tau^2}}\right)-U_1[\tilde{Q}]\left(\frac{\xi-\tau\nu(\xi)}{\sqrt{1+\tau^2}}\right)\right)&=0,
\\\lim_{\tau\to0+}\left(\left(\frac{\partial}{\partial\nu}U_1[\tilde{Q}]\right)\left(\frac{\xi+\tau\nu(\xi)}{\sqrt{1+\tau^2}}\right)-\left(\frac{\partial}{\partial\nu}U_1[\tilde{Q}]\right)\left(\frac{\xi-\tau\nu(\xi)}{\sqrt{1+\tau^2}}\right)\right)&=\tilde{Q}(\xi).
\end{align}
For the double-layer potential, we have
\begin{align}
\lim_{\tau\to0+}\left(U_2[Q]\left(\frac{\xi+\tau\nu(\xi)}{\sqrt{1+\tau^2}}\right)-U_2[Q]\left(\frac{\xi-\tau\nu(\xi)}{\sqrt{1+\tau^2}}\right)\right)&=-Q(\xi),
\\\lim_{\tau\to0+}\left(\left(\frac{\partial}{\partial\nu}U_2[Q]\right)\left(\frac{\xi+\tau\nu(\xi)}{\sqrt{1+\tau^2}}\right)-\left(\frac{\partial}{\partial\nu}U_2[Q]\right)\left(\frac{\xi-\tau\nu(\xi)}{\sqrt{1+\tau^2}}\right)\right)&=0.
\end{align}
\end{enumerate}
All of the relations above hold uniformly with respect to $\xi\in\partial\Gamma$.
\end{thm}
\begin{rem}
Theorem \ref{thm:limitrelc} essentially tells us that the single-layer potential $U_1[\tilde{Q}]$ and the normal derivative of the double-layer potential $\frac{\partial}{\partial\nu}U_2[Q]$ are continuous across the boundary $\partial\Gamma$ while the double-layer potential $U_2[Q]$ and the normal derivative of the single-layer potential $\frac{\partial}{\partial\nu}U_1[\tilde{Q}]$ are not. However, one has to be careful about $\frac{\partial}{\partial\nu}U_2[Q]$: it is only well-defined on $\partial\Gamma$ under higher smoothness assumptions on $Q$ than just ${\rm C}^{(0)}(\partial\Gamma)$. Therefore, we only supplied the jump relation for this particular case but not the limit relation, which is sufficient for most theoretical considerations.
\end{rem}

\begin{rem}\label{lem:adjont}
The relations in Theorem \ref{thm:limitrelc} were formulated with respect to the uniform topology for $Q,\tilde{Q}\in {\rm C}^{(0)}(\Omega)$. However, they can also be formulated with respect to the $\L^2(\Omega)$-topology for $Q,\tilde{Q}\in\L^2(\Omega)$. For details, the reader is again referred to \cite{freedengerhards12, freedenmichel04a} and earlier references therein.
\end{rem}

\section{Boundary Value Problems for the Beltrami Operator}\label{sec:spot2}

In this section, we investigate the problems (PP), (DP), and (NP) and try to obtain integral representations of their solutions.

\subsection{Poisson Problem}\label{sec:pp}
We remember the Poisson problem (PP) from the beginning of this chapter: Let $H$ be of class ${\rm C}^{(1)}(\overline{\Gamma})$, then we are looking for a function $U$ of class ${\rm C}^{(2)}(\Gamma)$ such that
\begin{align}\label{eqn:beltramiequation2}
\Delta^*U(\xi)&=H(\xi),\quad\xi\in\Gamma.
\end{align}
If we choose $\bar{H}=H-\frac{1}{\|\Gamma\|}\int_\Gamma H(\eta)d(\eta)$, we find that $\int_\Gamma \bar{H}(\eta)d\omega(\eta)=0$ and, by Theorem \ref{thm:poissoneqsphere}, that
\begin{align}
 \bar{U}( \xi)=\int_\Gamma G(\Delta^*;\xi\cdot\eta)\bar{H}(\eta)d\omega(\eta),\quad\xi\in\Gamma,
\end{align}
satisfies $\Delta^*\bar{U}(\xi)=\bar{H}(\xi)$, for $\xi\in\Gamma$. Setting $U(\xi)=\bar{U}(\xi)-\frac{1}{\|\Gamma\|}\ln(1-\xi\cdot\bar{\xi})\int_\Gamma H(\eta)d\omega(\eta)$, for some fixed $\bar{\xi}\in\Gamma^c$, 
we eventually obtain the desired solution satisfying
\begin{align}\label{eqn:deltau}
\Delta^*U(\xi)=\bar{H}(\xi)+\frac{1}{\|\Gamma\|}\int_\Gamma H(\eta)d\omega(\eta)=H(\xi),\quad \xi\in\Gamma.
\end{align}
The solution of \eqref{eqn:beltramiequation2}, however, is not unique. Subscribing further boundary values on $U$, e.g., Dirichlet boundary values $U^-(\xi)=F(\xi)$, for $\xi\in\partial\Gamma$, it is possible to obtain uniqueness. Letting $\tilde{U}$ denote the function $U$ from \eqref{eqn:deltau} that we constructed before, we can formulate the boundary value problem of finding a function $\tilde{\tilde{U}}$ that solves
\begin{align}
 \Delta^*\tilde{\tilde{U}}(\xi)&=0,\quad\xi\in\Gamma,\label{eqn:bv1}
 \\\tilde{\tilde{U}}^-(\xi)&=F(\xi)-\tilde{U}^-(\xi),\quad\xi\in\partial\Gamma.\label{eqn:bv2}
\end{align}
The newly obtained function $U=\tilde{U}+\tilde{\tilde{U}}$ would then satisfy the desired differential equation \eqref{eqn:beltramiequation2} and the desired Dirichlet boundary values. Boundary value problems such as \eqref{eqn:bv1}, \eqref{eqn:bv2} are studied in more detail in the upcoming section.

\subsection{Dirichlet and Neumann Problem}\label{sec:dnp}
We take a closer look at the following boundary value problems that have already been mentioned in the introduction: 
\begin{description}
\item[\textit{Dirichlet Problem} {\rm (DP):}]  Let $F$ be of class ${\rm C}^{(0)}(\partial\Gamma)$. We are looking for a function $U$ of class ${\rm C}^{(2)}(\Gamma)\cap {\rm C}^{(0)}(\overline{\Gamma})$ such that
\begin{align}
\Delta^*U(\xi)&=0,\qquad\,\,\,\xi\in\Gamma,
\\U^-(\xi)&=F(\xi), \quad\xi\in\partial\Gamma.\label{eqn:dv11}
\end{align}
\item[\textit{Neumann Problem} {\rm (NP):}]  Let $F$ be of class ${\rm C}^{(0)}(\partial\Gamma)$. We are looking for a function $U$ of class ${\rm C}^{(2)}(\Gamma)\cap {\rm C}^{(0)}(\overline{\Gamma})$, with a well-defined normal derivative $\frac{\partial}{\partial\nu}U^-$  on $\partial\Gamma$, such that
\begin{align}
\Delta^*U(\xi)&=0,\qquad\,\,\,\xi\in\Gamma,
\\\frac{\partial}{\partial\nu}U^-(\xi)&=F(\xi), \quad \xi\in\partial\Gamma.\label{eqn:nv11}
\end{align}
\end{description}
First, we formalize the term $U^-(\xi)$. For $\xi\in\partial\Gamma$, it is meant in the sense
\begin{align}
 U^-(\xi)=\lim_{\tau\to0+}U\left(\frac{\xi-\tau\nu(\xi)}{\sqrt{1+\tau^2}}\right),
\end{align}
i.e., we approach the boundary $\partial\Gamma$ in normal direction from within $\Gamma$. The term $U^+(\xi)$ is meant in the sense
\begin{align}
 U^+(\xi)=\lim_{\tau\to0+}U\left(\frac{\xi+\tau\nu(\xi)}{\sqrt{1+\tau^2}}\right),
\end{align}
i.e., we approach the boundary $\partial\Gamma$ in normal direction from the outside of $\Gamma$ (or, in other words, from within $\Gamma^c$). The expressions $\frac{\partial}{\partial\nu}U^\pm(\xi)$ are meant analogously. We can already see the connection to the limit- and jump-relations from Theorem \ref{thm:limitrelc}. More precisely, making the ansatz $U=U_2[Q]$ for the Dirichlet problem (DP) and $U=U_1[\tilde{Q}]$ for the Neumann problem (NP), Theorem \ref{thm:limitrelc} yields the following closely related problems:
\begin{description}
 \item[\textit{Integral Dirichlet Problem} {\rm (IDP):}]  Let $F$ be of class ${\rm C}^{(0)}(\partial\Gamma)$. We are looking for some $Q$ of class ${\rm C}^{(0)}(\partial\Gamma)$ that satisfies
\begin{align}
F(\xi)=U_2[Q](\xi)+\frac{1}{2}Q(\xi),\quad\xi\in\partial\Gamma.\label{eqn:inteqdirichlet}
\end{align}
 \item[\textit{Integral Neumann Problem} {\rm (INP):}] Let $F$ be of class ${\rm C}^{(0)}(\partial\Gamma)$. We are looking for some $\tilde{Q}$ of class ${\rm C}_0^{(0)}(\partial\Gamma)$ that satisfies
\begin{align}
F(\xi)=U_1[\tilde{Q}](\xi)-\frac{1}{2}\tilde{Q}(\xi),\quad\xi\in\partial\Gamma.\label{eqn:inteqneumann}
\end{align}
\end{description}
In other words, the Dirichlet problem (DP) and the Neumann problem (NP) have been reduced to the Fredholm equations (IDP) and (INP). These boundary integral formulations have been used, e.g., in \cite{gemmrich08, kropinski14} to numerically solve the original boundary value problems for the Beltrami operator. In this section, however, we are mainly interested in (IDP) and (INP) as tools to guarantee the existence of solutions to (DP) and (NP) via the Fredholm alternative. Uniqueness of the solutions can be obtained via the application of the maximum principle from Theorem \ref{thm:maxprinciple} and the Green formulas from Theorem \ref{thm:greenformulas}. 

\begin{rem}
There are two noteworthy differences in comparison to the Euclidean case. First, considerations on the sphere do not require a clear distinction between interior and exterior problems since the open complement $\Gamma^c$ of a bounded regular region $\Gamma\subset\Omega$ is again a bounded regular region. Second, the single-layer potential $U_1[\tilde{Q}]$ is only harmonic if $\tilde{Q}\in \textrm{C}^{(0)}_0(\partial\Gamma)$. A solution of \eqref{eqn:inteqneumann} in $\textrm{C}^{(0)}_0(\partial\Gamma)$ exists if and only if $F$ is of class $\textrm{C}_0^{(0)}(\partial\Gamma)$, which suits the general necessary condition for the existence of a solution to (NP) that can be obtained from Green's formulas. However, it should be mentioned that the integral equation \eqref{eqn:inteqneumann} additionally has a unique solution $\tilde{Q}\in\textrm{C}^{(0)}(\partial\Gamma)$ if $F$ is of class $\textrm{C}^{(0)}(\partial\Gamma)$. This is not true for the Euclidean counterpart.
\end{rem}

Summarizing, and including the considerations from Section \ref{sec:pp}, we obtain the following results. For details, the reader is again referred to \cite{freedengerhards12} and, for the Euclidean counterparts, to \cite{freedenmichel04a, helms69, kellogg67, wermer74}.

\begin{thm}[Uniqueness]\label{thm:uniquedirichlet}\textbf{}\\[-4ex]
\begin{enumerate}[(a)]
\item A solution of {\rm (DP)} is uniquely determined. 
\item A solution of {\rm (NP)} is uniquely determined up to an additive constant. 
\end{enumerate}
\end{thm}

\begin{thm}[Existence for Generalized (DP)]\label{thm:existuniquedirichlet}
Let $F$ be of class ${\rm C}^{(0)}(\partial\Gamma)$ and $H$ of class ${\rm C}^{(1)}(\overline{\Gamma})$. Then there exists a unique solution $U$ of class ${\rm C}^{(2)}(\Gamma)\cap {\rm C}^{(0)}(\overline{\Gamma})$ of the Dirichlet problem
\begin{align}
\Delta^*U(\xi)&=H(\xi),\quad\xi\in\Gamma,
\\U^-(\xi)&=F(\xi),\quad \xi\in\partial\Gamma.
\end{align}
\end{thm}

\begin{thm}[Existence for Generalized (NP)]\label{thm:existuniqueneumann}
Let $F$ be of class ${\rm C}^{(0)}(\partial\Gamma)$ and $H$ of class ${\rm C}^{(1)}(\overline{\Gamma})$. Then there exists an up to an additive constant uniquely determined solution $U$ of class ${\rm C}^{(2)}(\Gamma)\cap {\rm C}^{(0)}(\overline{\Gamma})$, with a well-defined normal derivative $\frac{\partial}{\partial\nu}U^-$ on $\partial\Gamma$, to the Neumann problem
\begin{align}
\Delta^*U(\xi)&=H(\xi),\quad\xi\in\Gamma,
\\\frac{\partial}{\partial\nu}U^-(\xi)&=F(\xi),\quad \xi\in\partial\Gamma,
\end{align}
if and only if
\begin{align}
\int_{\partial\Gamma} F(\eta)d\sigma(\eta)-\int_\Gamma H(\eta)d\omega(\eta)=0.\label{eqn:solvcondneumann2}
\end{align}
\end{thm}

\begin{prf}
The condition \eqref{eqn:solvcondneumann2} is a simple consequence from 
\begin{align}
 \int_\Gamma H(\eta)d\omega(\eta)=\int_\Gamma \Delta^*U(\eta)d\omega(\eta)=\int_{\partial\Gamma} \frac{\partial}{\partial\nu}U(\eta)d\sigma(\eta)=\int_{\partial\Gamma} F(\eta)d\sigma(\eta),
\end{align}
where Green's formulas have been used for the second equation. The general existence follows from the application of the Fredholm alternative to (INP).
\end{prf}

\subsection{Green's Functions}\label{sec:gf2}

Next, we are interested in the representation of a solution to (DP) and (NP). A possibility is indicated in Theorem \ref{thm:locfundthmd}(a). However, this representation requires the simultaneous knowledge of $U$ and $\frac{\partial}{\partial\nu}U$ on the boundary $\partial\Gamma$, which is not necessary and can be problematic since the two quantities are not independent from each other. As a remedy, Green's functions for Dirichlet and Neumann boundary values can be used. 

\begin{figure}
\begin{center}
\scalebox{0.75}{\includegraphics{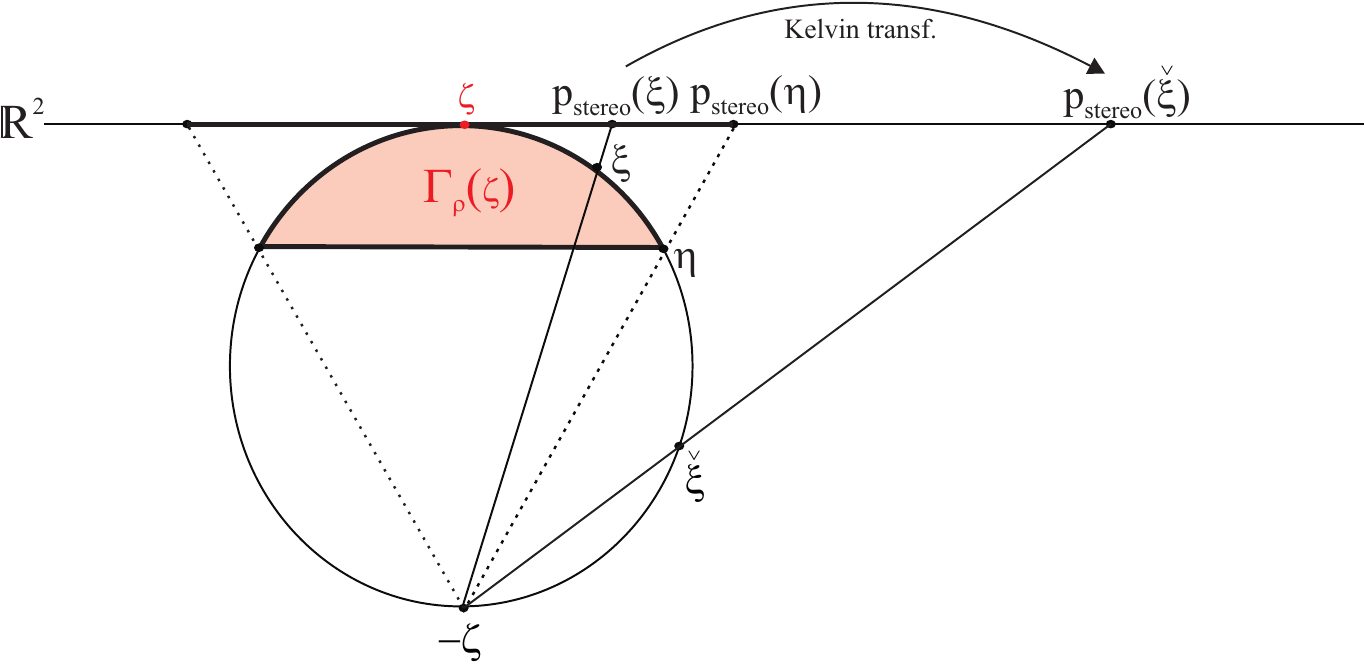}}
\end{center}
\caption{Schematic description of the construction of the reflection point $\check{\xi}$.}\label{fig:imagetrafo}
\end{figure}

More precisely, a function  $G_D(\Delta^*;\cdot,\cdot)$ is called a \emph{Dirichlet Green function (with respect to the Beltrami operator)} if it can be decomposed in the form
\begin{align}
G_D(\Delta^*;\xi,\eta)=G(\Delta^*;\xi\cdot\eta)-\Phi_D(\xi,\eta),\quad \eta\in\overline{\Gamma},\xi\in\Gamma,\xi\not=\eta,\label{eqn:dirgreenfuncdef}
\end{align}
where $\Phi_D(\xi,\cdot)$ is of class ${\rm C}^{(2)}(\Gamma)\cap {\rm C}^{(1)}(\overline{\Gamma})$ and satisfies
\begin{align}
\Delta_\eta^*\Phi_D(\xi,\eta)&=-\frac{1}{4\pi},\quad\eta\in\Gamma,\label{eqn:dg1}
\\\Phi_D^-(\xi,\eta)&=G(\Delta^*;\xi\cdot\eta),\quad \eta\in\partial\Gamma,\label{eqn:dg2}
\end{align}
for every $\xi\in\Gamma$. Analogously, a function $G_N(\Delta^*;\cdot,\cdot)$ is called a \emph{Neumann Green function (with respect to the Beltrami operator)} if it can be decomposed in the form
\begin{align}
G_N(\Delta^*;\xi,\eta)=G(\Delta^*;\xi\cdot\eta)-\Phi_N(\xi,\eta),\quad \eta\in\overline{\Gamma},\xi\in\Gamma,\xi\not=\eta,
\end{align}
where $\Phi_N(\xi,\cdot)$ is of class ${\rm C}^{(2)}(\Gamma)\cap {\rm C}^{(1)}(\overline{\Gamma})$ and satisfies the conditions
\begin{align}
\Delta^*_\eta \Phi_N(\xi,\eta)&=\frac{1}{\|\Gamma\|}-\frac{1}{4\pi},\quad \eta\in\Gamma,
\\\frac{\partial}{\partial\nu(\eta)}\Phi_N^-(\xi,\eta)&=\frac{\partial}{\partial\nu(\eta)}G(\Delta^*;\xi\cdot\eta),\quad\eta\in\partial\Gamma,
\end{align}
for every $\xi\in\Gamma$. Using Theorem \ref{thm:locfundthmd}(a) and Theorem \ref{thm:greenformulas}(c) for $\Phi_D$ and $\Phi_N$, we eventually achieve the representations
\begin{align}
U(\xi)=\int_\Gamma G_D(\Delta^*;\xi,\eta)\Delta^*_\eta U(\eta)d\omega(\eta)
+\int_{\partial\Gamma}U(\eta)\frac{\partial}{\partial\nu(\eta)}G_D(\Delta^*;\xi,\eta)d\sigma(\eta) \label{eqn:repdirgreenfuncsol}
\end{align}
and
\begin{align}
U(\xi)=&\frac{1}{\|\Gamma\|}\int_\Gamma U(\eta)d\omega(\eta)+\int_\Gamma G_N(\Delta^*;\xi,\eta)\Delta^*_\eta U(\eta)d\omega(\eta)\label{eqn:neumannbvrep}
\\&-\int_{\partial\Gamma}G_N(\Delta^*;\xi,\eta)\frac{\partial}{\partial\nu(\eta)}U(\eta)d\sigma(\eta),\nonumber
\end{align}
which yield integral representations for solutions to (DP) and (NP), respectively, under the condition that $U$ is of class ${\rm C}^{(2)}(\overline{\Gamma})$. It remains to construct the auxiliary functions $\Phi_D$ and $\Phi_N$. Some general construction principles on the sphere can be found, e.g., in \cite{gutkinnewton04, kidambinewton00}. In this chapter, we focus on spherical caps $\Gamma_\rho(\zeta)$. The procedure is similar to the construction of a Dirichlet Green function for a disc in $\mathbb{R}^2$. For $\xi\in\Gamma_\rho(\zeta)$, we need to find a reflection point $\check{\xi}\in(\Gamma_\rho(\zeta))^c$ and a scaling factor $\check{r}\in\mathbb{R}$ such that
\begin{align}
 1-\xi\cdot\eta=\check{r}\left(1-\check{\xi}\cdot\eta\right),\quad\eta\in\partial\Gamma_\rho(\zeta),\xi\in\Gamma_\rho(\zeta).\label{eqn:dgcond}
\end{align}
Indeed, under this assumption, it is clear that
\begin{align}\label{eqn:phidirichletdef}
\Phi_D(\xi,\eta)=\frac{1}{4\pi}\ln(\check{r}(1-\check{\xi}\cdot\eta))+\frac{1}{4\pi}(1-\ln(2))
\end{align}
satisfies the desired conditions \eqref{eqn:dg1} and \eqref{eqn:dg2}. The reflection point $\check{\xi}$ can be obtained by a stereographic projection of $\xi$ onto $\mathbb{R}^2$, then applying a Kelvin transform to the projection point, and eventually projecting it back to the sphere (cf. Figure \ref{fig:imagetrafo} for an illustration). $\check{\xi}$ represents the spherical Kelvin transformation of $\xi$. The scaling factor $\check{r}$ is obtained by solving \eqref{eqn:dgcond}. Alternatively, the entire Dirichlet Green function $G_D(\Delta^*;\cdot,\cdot)$ can be obtained from a stereographic projection of the Dirichlet Green function for the Laplace operator on a disc in $\mathbb{R}^2$. But this route would not supply us with a spherical counterpart to the Kelvin transform. We can conclude our considerations with the following theorem.
\begin{thm}\label{thm:dirichletgreenfunccap}
Let $\Gamma=\Gamma_{\rho}(\zeta)$ be a spherical cap with center $\zeta\in\Omega$ and radius $\rho\in(0,2)$. Furthermore, for $\xi\in\Gamma_{\rho}(\zeta)$ we set
\begin{align}
\check{\xi}&=\frac{1}{\check{r}}\xi-\frac{\check{r}-1}{\check{r}(\rho-1)}\zeta,\label{eqn:calcreflectpoints}
\\\check{r}&=-\frac{1+2\xi\cdot\zeta(\rho-1)+(\rho-1)^2}{\rho(\rho-2)}.\label{eqn:calcscalfact}
\end{align}
Then
\begin{align}
G_D(\Delta^*;\xi,\eta)=\frac{1}{4\pi}\ln(1-\xi\cdot\eta)-\frac{1}{4\pi}\ln(\check{r}(1-\check{\xi}\cdot\eta)),
\end{align}
and a solution $U\in{\rm C}^{(2)}(\overline{\Gamma})$ of the Dirichlet problem \emph{(DP)} can be represented by 
\begin{align}
U(\xi)&=\frac{1}{2\pi}\frac{\xi\cdot\zeta+\rho-1}{\sqrt{\rho(2-\rho)}}\int_{\partial\Gamma_{\rho}(\zeta)}\frac{1}{1-\xi\cdot\eta}F(\eta) \,d\sigma(\eta),\quad\xi\in\Gamma_\rho(\zeta).\label{eqn:solrepspheredirichlet}
\end{align}
\end{thm}

\begin{rem}
Applying Theorem \ref{thm:dirichletgreenfunccap} for $\zeta=\xi$ leads to the Mean Value Property II from Theorem \ref{thm:meanvalpropsphere2}.
\end{rem}

A Neumann Green function for the Beltrami operator cannot be obtained by a simple stereographic projection of the Neumann Green function for the Laplace operator on a disc in $\mathbb{R}^2$. But some computations based on the previously obtained auxiliary function $\Phi_D$ yield the following theorem.

\begin{thm}\label{lem:neumanngreenfunccap}
Let $\Gamma=\Gamma_{\rho}(\zeta)$ be a spherical cap with center $\zeta\in\Omega$ and radius $\rho\in(0,2)$. Furthermore, let $\check{\xi}$ and $\check{r}$ be given as in Theorem \ref{thm:dirichletgreenfunccap}. Then, a Neumann Green function is given by
\begin{align}
G_N(\Delta^*;\xi,\eta)=\frac{1}{4\pi}\ln(1-\xi\cdot\eta)+\frac{1}{4\pi}\ln(\check{r}(1-\check{\xi}\cdot\eta))+\frac{1-\rho}{2\pi\rho}\ln(1+\zeta\cdot\eta).
\end{align}
A solution solution $U\in{\rm C}^{(2)}(\overline{\Gamma})$ of the Neumann problem \emph{(NP)} can be represented by 
\begin{align}
U(\xi)=&\frac{1}{2\pi\rho}\int_{\Gamma_\rho(\zeta)}U(\eta)d\omega(\eta)
\\&-\int_{\partial\Gamma_\rho(\zeta)}\left(\frac{1}{2\pi}\ln(1-\xi\cdot\eta)+\frac{1-\rho}{2\pi\rho}\ln(2-\rho)\right)F(\eta)\,d\sigma(\eta),\quad\xi\in\Gamma_\rho(\zeta).\nonumber
\end{align}
\end{thm}

\section{Spherical Decompositions and First Order Differential Equations}\label{sec:sdecomp}

In this section we treat differential equations for the surface gradient $\nabla^*$ and the surface curl gradient $\L^*$. They come up, e.g., when dealing with vertical deflections and geostrophic ocean flow. Additionally, we take a look at some spherical decompositions of vector fields that are particularly useful in geosciences. 

\subsection{Surface Gradient and Surface Curl Gradient}\label{sec:surfacegrad}
Different from the Poisson equation, solutions of the differential equations with respect to the surface gradient and the surface curl gradient are uniquely determined up to an additive constant on regular regions $\Gamma\subset\Omega$, without the necessity of boundary values. Also the existence of a solution can be easily guaranteed. This is summarized in the following two lemmas.

\begin{lem}[Uniqueness]\label{lem:nablagammaunique}
Let $U$ be of class ${\rm C}^{(1)}(\Gamma)$. Then 
\begin{align}
\nabla^*U(\xi)=0, \quad\xi\in\Gamma,
\end{align}
if and only if $U$ is constant on $\Gamma$. The same holds true for $\L^*U(\xi)=0$, $\xi\in\Gamma$.
\end{lem}

\begin{thm}[Existence]\label{thm:existnabladgl}\hfill
\begin{itemize}
 \item[(a)] Let $f\in{\rm c}^{(1)}(\Gamma)$ be a tangential vector field satisfying 
\begin{align}
\L^*\cdot f(\xi)=0, \quad\xi\in\Gamma.
\end{align}
Then there exists a function $U$ of class ${\rm C}^{(2)}(\Gamma)$, which is uniquely determined up to an additive constant, such that
\begin{align}
f(\xi)=\nabla^*U(\xi),\quad\xi\in\Gamma.
\end{align}
\item[(b)] Let $f\in{\rm c}^{(1)}(\Gamma)$ be a tangential vector field satisfying 
\begin{align}
\nabla^*\cdot f(\xi)=0,\quad\xi\in\Gamma.
\end{align}
Then there exists a function $U$ of class ${\rm C}^{(2)}(\Gamma)$, which is uniquely determined up to an additive constant, such that
\begin{align}
f(\xi)=\L^*U(\xi),\quad\xi\in\Gamma.
\end{align}
\end{itemize}
\end{thm}

From Theorem \ref{thm:locfundthmd}(b), we know a possible expression of the solutions to the differential equations for the surface gradient and the surface curl gradient. However, this representation requires the knowledge of $U$ on the boundary $\partial\Gamma$, which is actually not necessary according to Theorem \ref{thm:existnabladgl}. Using a Neumann Green function together with the identities in Theorem \ref{thm:locfundthmd}(b) directly implies the following results.

\begin{thm}\label{thm:gstarrepdgl}\hfill
\begin{itemize}
\item[(a)] Let $f$ of class ${\rm c}^{(1)}(\overline{\Gamma})$ be a tangential vector field satisfying $\L^*\cdot f(\xi)=0$, $\xi\in\Gamma$. Then a solution of
\begin{align}
f(\xi)=\nabla^*U(\xi),\quad\xi\in\Gamma,
\end{align}
is given by
\begin{align}
U(\xi)=\frac{1}{\|\Gamma\|}\int_\Gamma U(\eta)d\omega(\eta)-\int_\Gamma\left(\nabla^*_\eta G_N(\Delta^*;\xi,\eta)\right)\cdot f(\eta)d\omega(\eta).
\end{align}
\item[(b)] Let $f$ of class ${\rm c}^{(1)}(\overline{\Gamma})$ be a tangential vector field satisfying $\nabla^*\cdot f(\xi)=0$, $\xi\in\Gamma$. Then a  solution of
\begin{align}
f(\xi)=\L^*U(\xi),\quad\xi\in\Gamma,
\end{align}
is given by
\begin{align}
U(\xi)=\frac{1}{\|\Gamma\|}\int_\Gamma U(\eta)d\omega(\eta)-\int_\Gamma\left(\L^*_\eta G_N(\Delta^*;\xi,\eta)\right)\cdot f(\eta)d\omega(\eta).
\end{align}
\end{itemize}
\end{thm}

\begin{rem}
If we deal with the entire sphere $\Gamma=\Omega$, the same results as in the preceding theorem hold true. For the integral representations, one simply has to substitute the Neumann Green function by the fundamental solution $G(\Delta^*;\cdot)$.
\end{rem}

\subsection{Helmholtz and Hardy-Hodge Decomposition}

We begin with the spherical Helmholtz decomposition of a vector field $f$. It essentially describes the split-up of the vector field into a radial and two tangential components, of which one is surface curl-free and the other one surface divergence-free. In geomagnetism, this has applications, e.g., in the separation of polar ionospheric current systems into field-aligned currents (which are nearly radial in polar regions) and Pedersen and Hall currents (see, e.g., \cite{amm97, backus96, gerhards11a, olsen97}). In other areas, the spherical Helmholtz decomposition has a natural connection as well: geostrophic ocean flow, e.g., is purely tangential and surface divergence-free while the vertical deflection of the geoidal normal vector is approximately purely tangential and surface curl-free. For convenience, we use the following notations for the Helmholtz operators acting on a scalar function $F$ at a point $\xi\in\Omega$:
\begin{align} \label{eqn:helmops}
 o^{(1)}F(\xi)=\xi F(\xi),\qquad o^{(2)}F(\xi)=\nabla^*F(\xi),\qquad o^{(3)}F(\xi)=\L^* F(\xi).
\end{align}
Writing $f=o^{(1)}F_1+o^{(2)}F_2+o^{(3)}F_3$ on a subdomain $\Gamma$ and using the orthogonality of the three operators, we obtain $\Delta^*F_3(\xi)=\L^*\cdot f(\xi)$, $\xi\in\Gamma$. Latter can be solved by the methods of the previous section. We need to prescribe boundary data on $F_3$ in order to obtain uniqueness of the scalar function $F_3$. All in all, we can formulate Decomposition Theorem \ref{thm:spherehelmloc}.  More details can be found, e.g., in \cite{freedengerhards12, gerhards11a}.

\begin{thm}[Spherical Helmholtz Decomposition]\label{thm:spherehelmloc} 
Let $f$ be of class ${\rm c}^{(2)}(\overline{\Gamma})$. Then there exist scalar fields $F_1$ of class ${\rm C}^{(2)}(\overline{\Gamma})$ and $F_2$, $F_3$ of class ${\rm C}^{(2)}(\Gamma)$ such that
\begin{align}
f(\xi)&=o^{(1)} F_1(\xi)+o^{(2)} F_2(\xi)+o^{(3)} F_3(\xi),\quad\xi\in\Gamma.
\end{align}
Uniqueness of $F_1,F_2,F_3$ is guaranteed by the properties
\begin{align}
\int_\Gamma F_2(\eta)d\omega(\eta)=0
\end{align}
and
\begin{align}\label{eqn:helmuniquedir}
F_3^-(\xi)=F(\xi),\quad\xi\in\partial\Gamma,
\end{align}
for a fixed function $F$ of class ${\rm C}^{(0)}(\partial\Gamma)$. The Helmholtz scalars $F_1$, $F_2$, and $F_3$ can be then represented by
\begin{align}
F_2(\xi)=&-\int_\Gamma\left(\nabla_\eta^*G_N(\Delta^*;\xi,\eta)\right)\cdot f(\eta)d\omega(\eta)
\\&+\int_{\partial\Gamma}F(\eta)\,\tau_\eta\cdot\nabla_\eta^*G_N(\Delta^*;\xi,\eta)d\sigma(\eta),\quad\xi\in\Gamma\nonumber
\\F_3(\xi)=&-\int_\Gamma \left(\L_\eta^*G_D(\Delta^*;\xi,\eta)\right)\cdot f(\eta)d\omega(\eta)
\\&+\int_{\partial\Gamma}G_D(\Delta^*;\xi,\eta)\tau_\eta\cdot f(\eta)d\sigma(\eta)\nonumber
\\&+\int_{\partial\Gamma}F(\eta)\frac{\partial}{\partial\nu_\eta}G_D(\Delta^*;\xi,\eta)d\sigma(\eta),\quad\xi\in\Gamma\nonumber
\end{align}
for $\xi\in\Gamma$. Additionally, if $\int_\Gamma F_1(\eta)d\omega(\eta)=0$, then
\begin{align}
F_1(\xi)&=\xi\cdot f(\xi)=\Delta^*_\xi\int_\Gamma G(\Delta^*;\xi\cdot\eta)\,\eta\cdot f(\eta)d\omega(\eta),\quad\xi\in\Gamma.
\end{align}
\end{thm}

\begin{rem}
Clearly, the type of boundary conditions that have to be prescribed to obtain uniqueness of the Helmholtz decomposition can be varied. They can be imposed on $F_2$ instead of $F_3$, or the Dirichlet boundary conditions can be substituted by Neumann boundary conditions. Neumann boundary conditions are occasionally more advantageous as they allow the imposition of boundary information on the normal and tangential direction of the vectorial quantities $o^{(2)}F_2$ and $o^{(3)}F_3$, respectively, which are in some cases better accessible from the given data than the scalars $F_2$ or $F_3$. Representations analogous to Theorem \ref{thm:spherehelmloc} can be derived by Green's formulas and the results from Section \ref{sec:spot2}.
\end{rem}

\begin{rem}\label{rem:helmrepglob}
For the particular case $\Gamma=\Omega$, the results from Theorem \ref{thm:spherehelmloc} hold true as well well if the boundary integrals $\int_{\partial\Gamma}\ldots d\sigma$ are dropped and the Neumann and Dirichlet Green functions are substituted by the fundamental solution $G(\Delta^*;\cdot)$. For the uniqueness, condition \eqref{eqn:helmuniquedir} has to be substituted by $\int_\Omega F_3(\eta)d\omega(\eta)=0$. We then obtain
\begin{align}
F_2(\xi)&=-\int_\Omega \left(\nabla^*_\eta G(\Delta^*;\xi\cdot\eta)\right)\cdot f(\eta)d\omega(\eta),\quad\xi\in\Omega,
\\F_3(\xi)&=-\int_\Omega\left(\L^*_\eta G(\Delta^*;\xi\cdot\eta)\right)\cdot f(\eta)d\omega(\eta),\quad\xi\in\Omega.
\end{align}
Thus, in the global case $\Gamma=\Omega$, the Helmholtz scalars $F_2$ and $F_3$ are determined uniquely up to an additive constant without further constraints. The vectorial quantities $o^{(2)}F_2$ and $o^{(3)}F_3$ are actually uniquely determined. This is not true for general subdomains $\Gamma\subset\Omega$.
\end{rem}
 
Next, we turn to a different spherical decomposition, the so-called spherical Hardy-Hodge decomposition (the name is adopted from the Euclidean decomposition presented in \cite{baratchart13}, although its spherical version is known and used significantly longer, e.g., in \cite{backus96, freedenschreiner09, gerhards11a, mayermaier06, olsen10b} and references therein). It is based on the set of operators
\begin{align}
 \oox=o^{(1)}\left(\ddd+\frac{1}{2}\right)-\oy,\qquad \ooy=o^{(1)}\left(\ddd-\frac{1}{2}\right)+\oy,\qquad \ooz=\oz, 
\end{align}
where the operator $\ddd$ is given by $\ddd=\left(-\Delta^*+\frac{1}{4}\right)^{\frac{1}{2}}$. A decomposition in terms of these operators can be interpreted as a decomposition of a spherical vectorial signal with respect to sources lying inside a given sphere (reflected by the $\oox$-contributions), sources lying in the exterior of the sphere ($\ooy$-contributions), and sources on the sphere ($\ooz$-contributions). For the gravitational field measured at satellite altitude, e.g., only the $\oox$-contribution is of relevance. Concerning the Earth's crustal magnetization, only the $\ooy$-contribution of the magnetization generates a magnetic effect at satellite altitude. The generated magnetic field itself, however, only consists of $\oox$-contributions since its source (i.e., the magnetization) is located inside the satellite's orbit. The decomposition and the integral representation of its scalar functions can be closely related to the spherical Helmholtz decomposition. For details, we refer the reader to \cite{freedengerhards12, freedenschreiner09, gerhards11a, gerhards12}. Yet, the non-local structure of the operator $\ddd$ makes it very difficult to obtain results on subdomains $\Gamma\subset\Omega$. Therefore, the following theorem only treats the decomposition for the case $\Gamma=\Omega$.  

\begin{thm}[Spherical Hardy-Hodge Decomposition]\label{thm:hardyhodge}
Let $f$ be of class ${\rm c}^{(1)}({\Omega})$. Then there exist scalar fields $\tilde{F}_1$, $\tilde{F}_2$, $\tilde{F}_3$ of class ${\rm C}^{(2)}(\Omega)$ such that
\begin{align}
f(\xi)&=\oox \tilde{F}_1(\xi)+\ooy \tilde{F}_2(\xi)+\ooz \tilde{F}_3(\xi),\quad\xi\in\Omega.
\end{align}
Uniqueness of $\tilde{F}_1,\tilde{F}_2,\tilde{F}_3$ is guaranteed by the properties
\begin{align}
\int_\Omega \tilde{F}_3(\eta)d\omega(\eta)&=0,
\\\int_\Omega \tilde{F}_1(\eta)-\tilde{F}_2(\eta)d\omega(\eta)&=0.
\end{align}
The Hardy-Hodge scalars $\tilde{F}_1$, $\tilde{F}_2$, and $\tilde{F}_3$ can then be represented by
\begin{align}
\tilde{F}_1&=\frac{1}{2}\ddd^{-1}F_1+\frac{1}{4}\ddd^{-1}F_2-\frac{1}{2}F_2, \label{eqn:hh1}
\\\tilde{F}_2&=\frac{1}{2}\ddd^{-1}F_1+\frac{1}{4}\ddd^{-1}F_2+\frac{1}{2}F_2,  \label{eqn:hh2}
\\\tilde{F}_3&=F_3, \label{eqn:hh3}
\end{align}
where $F_1$, $F_2$, $F_3$ are the Helmholtz scalars from Theorem \ref{thm:spherehelmloc} and Remark \ref{rem:helmrepglob}.
\end{thm}

\begin{rem}
The operator $\ddd^{-1}$ can be represented as the convolution operator
\begin{align}
 \ddd^{-1}F(\xi)=\frac{1}{2\pi}\int_\Omega\frac{1}{\sqrt{2(1-\xi\cdot\eta)}}F(\eta)d\omega(\eta),\quad\xi\in\Omega,
\end{align}
acting on a function $F$ of class ${\rm C}^{(0)}(\Omega)$. Thus, equations \eqref{eqn:hh1}--\eqref{eqn:hh3} together with Theorem \ref{thm:spherehelmloc} and Remark \ref{rem:helmrepglob} form integral representations of the Hardy-Hodge scalars.
\end{rem}

\section{Complete Function Systems}\label{sec:cfs}

In the Euclidean setting, spherical harmonics form a complete function system in $\L^2(\Omega_R)$, and their harmonic extensions into the ball $\mathcal{B}_R=\{x\in\mathbb{R}^3:|x|<R\}$ and its exterior $\mathcal{B}_R^c=\{x\in\mathbb{R}^3:|x|>R\}$ (so-called inner and outer harmonics, respectively) form suitable function systems to approximate functions that are harmonic with respect to the Laplace operator. The limit- and jump-relations of layer potentials enable the extension of the completeness results to more general manifolds than the sphere. With the considerations of the previous sections at hand, we are now able to formulate analogous completeness results for function systems on general curves $\partial\Gamma$. We obtain completeness for certain function systems in $\L^2(\partial\Gamma)$ whose harmonic extensions into $\Gamma\subset\Omega$ are particularly well-suited for the approximation of functions that are harmonic with respect to the Beltrami operator.

\begin{figure}
\begin{center}
\scalebox{0.6}{\includegraphics{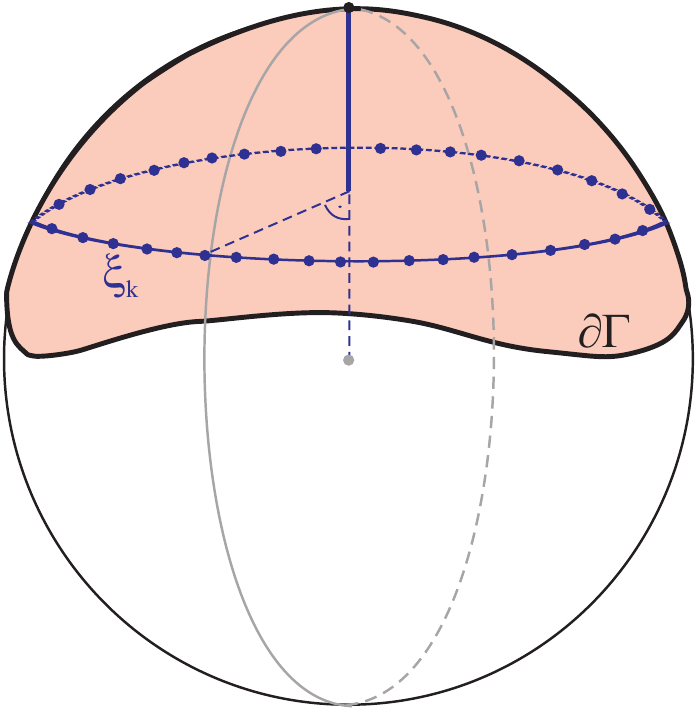}}
\end{center}
\caption{Example for a fundamental system $\{\xi_k\}_{k\in\mathbb{N}}$ (with respect to $\Gamma$). }\label{fig:bjercap}
\end{figure}

First, we need the notion of a fundamental system: Suppose that $\{\xi_k\}_{k\in\mathbb{N}}\subset\Gamma$ is a set of points satisfying
\begin{align}
\textnormal{dist}(\{\xi_k\}_{k\in\mathbb{N}},\partial\Gamma)>0.
\end{align}
If, for any harmonic function $F$ in $\Gamma$, the condition $F(\xi_k)=0$, $k\in\mathbb{N}$, implies that $F(\xi)=0$ for all $\xi\in\Gamma$, then we call $\{\xi_k\}_{k\in\mathbb{N}}$ a \emph{fundamental system (with respect to $\Gamma$)}. Assuming that $\Sigma\subset\Gamma$ is a regular region with dist$(\Sigma,\partial\Gamma)>0$, an example for such a fundamental system is given by a dense point set $\{\xi_k\}_{k\in\mathbb{N}}\subset\partial\Sigma$. A particularly simple choice for $\Sigma$ is a spherical cap within $\Gamma$ (cf. Figure \ref{fig:bjercap}).

We begin with the completeness of function systems based on the fundamental solution for the Beltrami operator.

\begin{thm}\label{thm:completefusystemln}
Let $\{\xi_k\}_{k\in\mathbb{N}}$ be a fundamental system with respect to $\Gamma$. Then the following statements hold true:
\begin{itemize}
\item[(a)] The function system $\{ G_k\}_{k\in\mathbb{N}_0}$ given by 
\begin{align}
 G_k(\xi)=\frac{1}{4\pi}\ln(1-\xi_k\cdot\xi), \,k\in\mathbb{N},\qquad G_0(\xi)=\frac{1}{4\pi},\nonumber
\end{align}
is complete, and hence closed in $\L^2(\partial\Gamma)$.
\item[(b)] The function system $\{\tilde{ G}_k\}_{k\in\mathbb{N}_0}$, given by 
\begin{align}
\tilde{ G}_k(\xi)=\frac{1}{4\pi}\frac{\partial}{\partial\nu(\xi)}\ln(1-\xi_k\cdot\xi),\, k\in\mathbb{N},\qquad\tilde{ G}_0(\xi)=\frac{1}{4\pi},\nonumber
\end{align}
is complete, and hence closed in $\L^2(\partial\Gamma)$.
\end{itemize}
\end{thm}

\begin{rem}\label{rem:phikmod}
Let us assume that $\{\xi_k\}_{k\in\mathbb{N}}$ is a fundamental system with respect to $\Gamma^c$. Then the functions $\tilde{ G}_k$ from Theorem \ref{thm:completefusystemln} are harmonic in $\Gamma$ and, thus, particularly suitably for the approximation of harmonic functions in $\Gamma$. The functions $G_k$ from Theorem \ref{thm:completefusystemln} need to be modified since they only satisfy $\Delta^* G_k(\xi)=-\frac{1}{4\pi}$, for $\xi\in\Gamma$ and $k\in\mathbb{N}$. Any auxiliary function $ G$ of class ${\rm C}^{(2)}(\overline{\Gamma})$ that satisfies $\Delta^* G(\xi)=\frac{1}{4\pi}$, $\xi\in\Gamma$, can be added to $ G_k$ without changing the completeness property. In other words, e.g.,
\begin{align}
G_k^{\rm (mod)}(\xi)= G_k(\xi)-\frac{1}{4\pi}\ln(1-\xi\cdot\bar{\xi}),\,k\in\mathbb{N},\qquad  G_0^{\rm (mod)}(\xi)= G_0(\xi),\nonumber
\end{align}
with a fixed $\bar{\xi}\in\Gamma^c$, forms a complete function system in $\L^2(\partial\Gamma)$ that additionally satisfies $\Delta^* G_k^{\rm (mod)}(\xi)=0$, $\xi\in\Gamma$.  
\end{rem}

Next, we want to transfer the results from Theorem \ref{thm:completefusystemln} to inner harmonics for spherical caps. In order to achieve this, we first need to clarify what we mean by inner harmonics for spherical caps. The sine and cosine functions obviously take the role of spherical harmonics on a circle in $\mathbb{R}^2$. Their harmonic continuations into the disc $\mathcal{D}_R=\{x\in\mathbb{R}^2:|x|<R\}$ with radius $R>0$ and into its exterior $\mathcal{D}_R^c=\{x\in\mathbb{R}^2:|x|>R\}$ (the so-called inner and outer harmonics, respectively) are given by
\begin{align}
 &H_{n,1}^{(int)}(R;x)=\frac{1}{R\sqrt{\pi}}\left(\frac{r}{R}\right)^n\cos(n\varphi),\,n\in\mathbb{N}_0,\quad x\in \mathcal{D}_R,
 \\&H_{n,2}^{(int)}(R;x)=\frac{1}{R\sqrt{\pi}}\left(\frac{r}{R}\right)^n\sin(n\varphi),\,n\in\mathbb{N},\quad x\in \mathcal{D}_R,
 \\&H_{n,1}^{(ext)}(R;x)=\frac{1}{R\sqrt{\pi}}\left(\frac{R}{r}\right)^n\cos(n\varphi),\,n\in\mathbb{N}_0,\quad x\in \mathcal{D}_R^c,
 \\&H_{n,2}^{(ext)}(R;x)=\frac{1}{R\sqrt{\pi}}\left(\frac{R}{r}\right)^n\sin(n\varphi),\,n\in\mathbb{N},\quad x\in \mathcal{D}_R^c,
\end{align}
where $x=(r\cos(\varphi),r\sin(\varphi))^T$, $r\geq0$, $\varphi\in[0,2\pi)$. Inner harmonics on a spherical cap $\Gamma_\rho(\zeta)$ with radius $\rho\in(0,2)$ and center $\zeta\in\Omega$ can then be obtained by a simple stereographic projection. More precisely,
\begin{align}
H_{n,k}^{\rho,\zeta}(\xi)=H_{n,k}^{(int)}\left(\rho^{\frac{1}{4}}(2-\rho)^{\frac{1}{4}};p_{stereo}(\zeta;\xi)\right),\quad \xi\in\Gamma_\rho(\zeta),
\end{align}
denotes an \emph{inner harmonic (of degree $n$ and order $k$) on $\Gamma_\rho(\zeta)$}. The applied stereographic projection $p_{stereo}(\zeta;\cdot):\Omega\setminus\{-\zeta\}\to\mathbb{R}^2$ is defined via
\begin{align}
p_{stereo}(\zeta;\xi)=\left(\frac{2\xi\cdot (\mathbf{t}\varepsilon^1)}{1+\xi\cdot\zeta},\frac{2\xi\cdot (\mathbf{t}\varepsilon^2)}{1+\xi\cdot\zeta}\right),
\end{align}
where $\varepsilon^1=(1,0,0)^T$, $\varepsilon^2=(0,1,0)^T$, $\varepsilon^3=(0,0,1)^T$ denotes the canonical basis in $\mathbb{R}^3$ and $\mathbf{t}\in\mathbb{R}^{3\times3}$ a rotation matrix with $\mathbf{t}\varepsilon^3=\zeta$. From the harmonicity of $H_{n,k}^{(int)}(R;\cdot)$ in $\mathcal{D}_R$ it follows that $ H_{n,k}^{\rho,\zeta}$ is harmonic in $\Gamma_\rho(\zeta)$. Note that, as always, harmonicity in the Euclidean space $\mathbb{R}^2$ is meant with respect to the Laplace operator while it is meant with respect to the Beltrami operator when we are intrinsic on the sphere $\Omega$. Opposed to the Euclidean case, outer harmonics for spherical caps do not play a distinct role. Actually, for a spherical cap $\Gamma_\rho(\zeta)$, the corresponding outer harmonics coincide with the inner harmonics for the spherical cap $\left(\Gamma_\rho(\zeta)\right)^c=\Gamma_{2-\rho}(-\zeta)$, which is why we do not consider them separately. The relation
\begin{align}
 \ln(1-\xi\cdot\eta)=&-\ln(2)+\ln(1+\xi\cdot\zeta)+\ln(1-\eta\cdot\zeta)
 \\&-\sqrt{\rho(2-\rho)}\pi\sum_{n=1}^\infty \sum_{k=1}^2\frac{2}{n}H_{n,k}^{\rho,\zeta}(\xi)H_{n,k}^{2-\rho,-\zeta}(\eta),\nonumber
\end{align}
for $\xi\in\Omega\setminus\{-\zeta\}$, $\eta\in\Omega\setminus\{\zeta\}$, and $|p_{stereo}(\zeta;\xi)|<|p_{stereo}(\zeta;\eta)|$, eventually allows to transfer the completeness results from Theorem \ref{thm:completefusystemln} to inner harmonics on spherical caps (for details, the reader is referred to \cite{freedengerhards12}).

\begin{thm}\label{thm:completefusystemih}
Let $\Gamma_\rho(\zeta)$ be a spherical cap with $\overline{\Gamma}\subset\Gamma_\rho(\zeta)$. Then the following statements hold true:
\begin{itemize}
\item[(a)] The inner harmonics $\big\{H_{0,1}^{\rho,\zeta}\big\}\cup\big\{H_{n,k}^{\rho,\zeta}\big\}_{n\in\mathbb{N},k=1,2}$ form a complete, and hence closed function system in $\L^2(\partial\Gamma)$.
\item[(b)] The normal derivatives of the inner harmonics, i.e., $\big\{H_{0,1}^{\rho,\zeta}\big\}\cup\big\{\frac{\partial}{\partial\nu}H_{n,k}^{\rho,\zeta}\big\}_{n\in\mathbb{N},k=1,2}$ form a complete and hence closed function system in $\L^2(\partial\Gamma)$.
\end{itemize}
\end{thm}

We conclude this section by stating the use of the function systems from above for the approximation of solutions to the spherical boundary value problems (DP) and (NP) from Section \ref{sec:spot2}.

\begin{thm}
Let $\{\Phi_k\}_{k\in\No}$ denote one of the function systems introduced in Theorem \ref{thm:completefusystemln}(b), Remark \ref{rem:phikmod}, or Theorem \ref{thm:completefusystemih}, and  $U\in{\rm C}^{(2)}(\Gamma)\cap{\rm C}^{(0)}(\overline{\Gamma})$ be a solution of one of the boundary value problems (DP) or (NP). Then, for every $\varepsilon>0$, there exist $M\in\mathbb{N}_0$ and coefficients $a_k\in\mathbb{R}$, $k=0,1,\ldots,M$, such that
\begin{align}    
\left\|U-\sum_{k=0}^Ma_k\Phi_k\right\|_{\L^2(\Gamma)}<\varepsilon.
\end{align}
The choice of $M$ and  the coefficients $a_k$, $k=1,\ldots,M$, can be based solely on an approximation of $U$ or $\frac{\partial}{\partial\nu}U$ on the boundary $\partial\Gamma$.
\end{thm}

\begin{rem}
All the density and approximation results that were obtained in this section in an $\L^2$-context also hold true in a ${\rm C}^{(0)}$-context with respect to the uniform topology and can be shown by the tools supplied throughout this chapter (see, e.g., \cite{freedengerhards12, freedenmichel04a}).
\end{rem}

\section{Applications in Geoscience}\label{sec:app}

In this section, we present some applications of the previous tools to the approximation of different quantities of interest in physical geodesy. More precisely, we use techniques from Section \ref{sec:surfacegrad} to reconstruct the disturbing potential from given vertical deflections over South America and the mean dynamic ocean topography (MDT) from given geostrophic ocean flow patterns over the Pacific Ocean, respectively. We will be rather brief about the geophysical derivations of the underlying spherical differential equations and refer the reader to classical literature such as \cite{heiskanenmoritz67, wellenhof05, pedlosky79, stewart}. The particular formulations of our setting can also be found, e.g., in \cite{fehlinger09, fehlinger07, freeden15, freedenschreiner09}. Opposed to the latter, our reconstructions in Sections \ref{sec:vd} and \ref{sec:goc} are based on the approach in Section \ref{sec:surfacegrad} via Neumann Green functions and does not require boundary information for the spherical caps under consideration. In Section \ref{sec:vort}, based on the results from Section \ref{sec:cfs}, we address a model problem motivated by point vortex motion on the sphere.

\subsection{Vertical Deflections}\label{sec:vd}

\begin{figure}
\begin{center}
\includegraphics[scale=0.3]{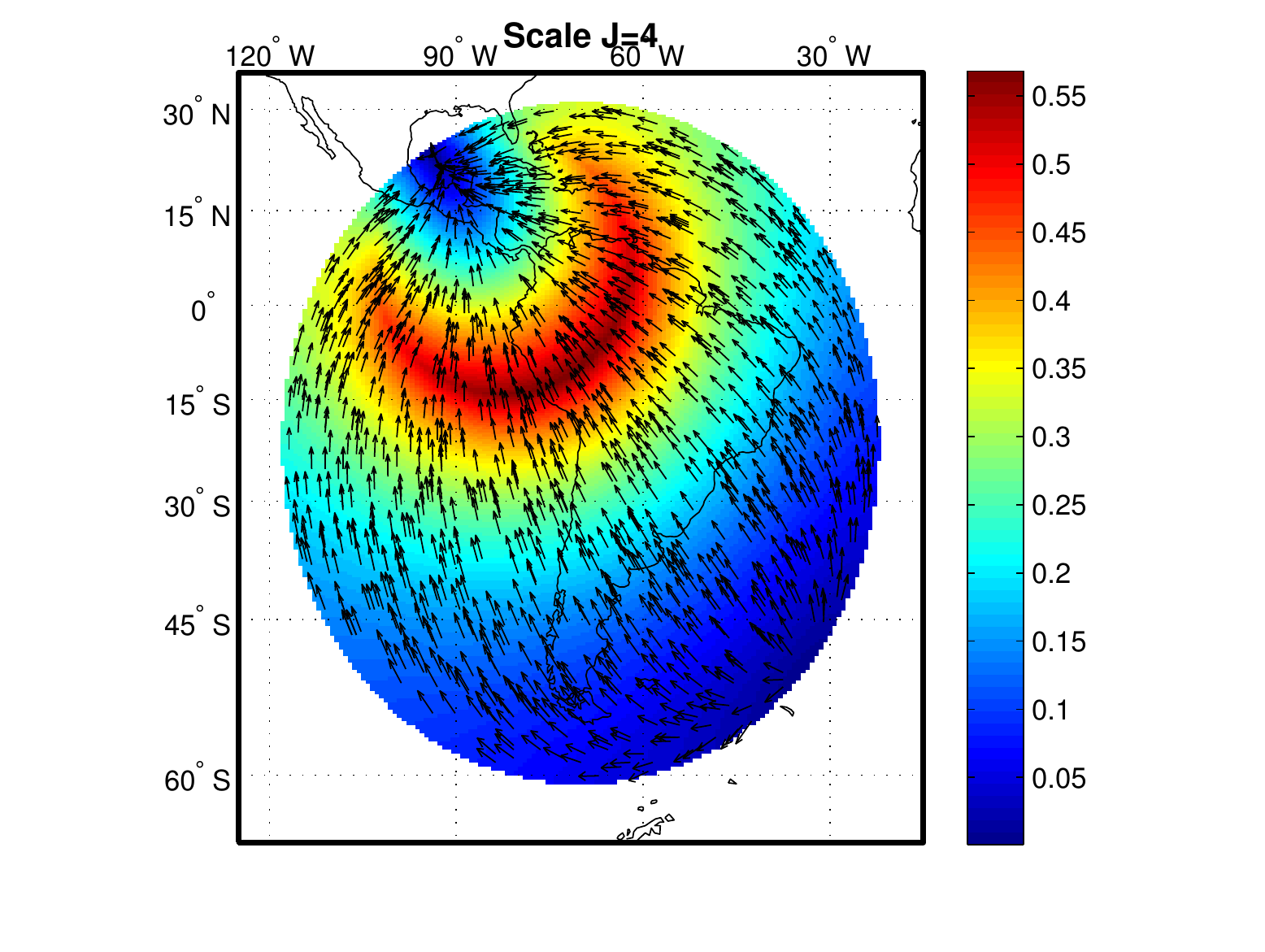} \includegraphics[scale=0.3]{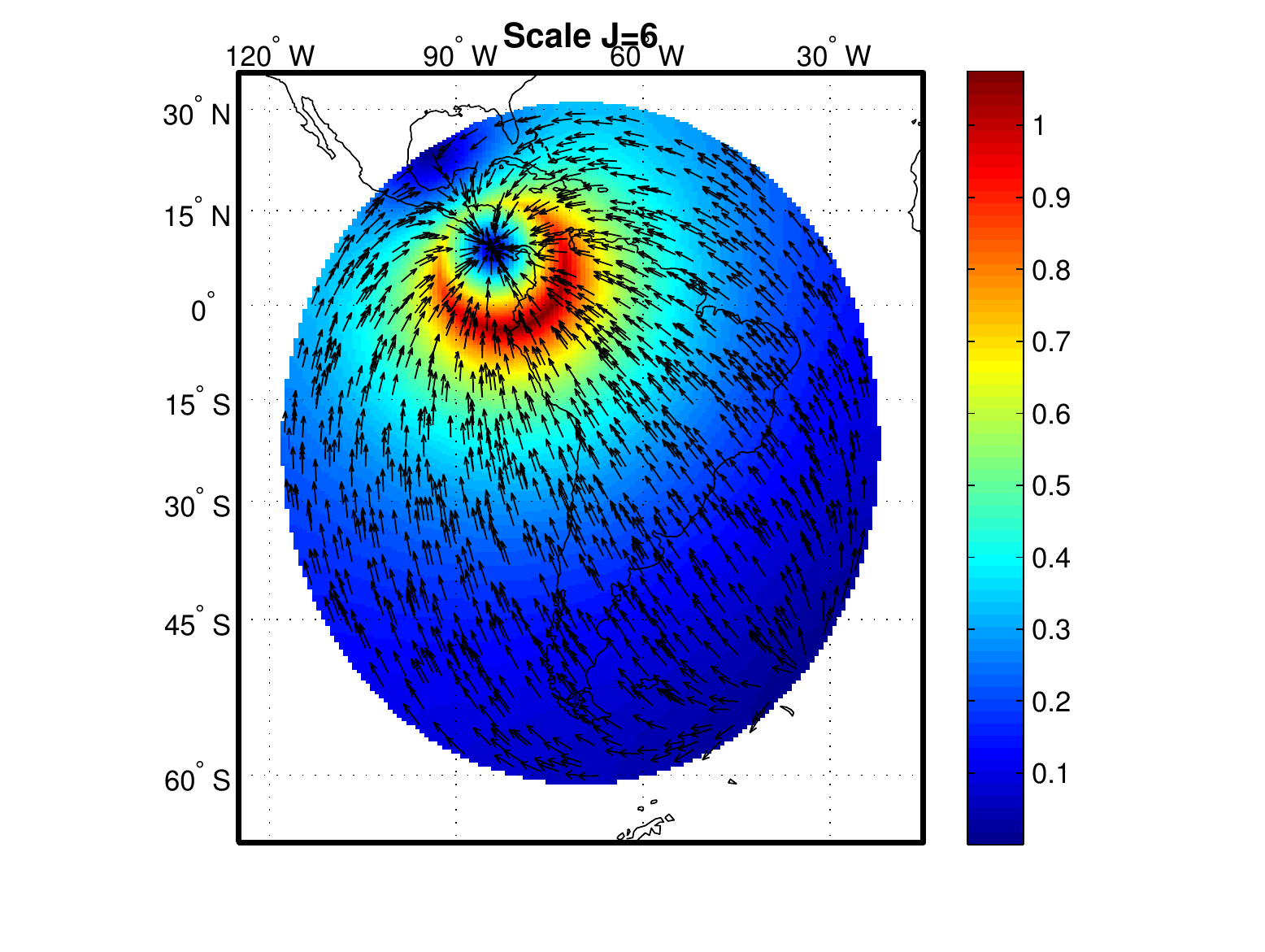} \includegraphics[scale=0.3]{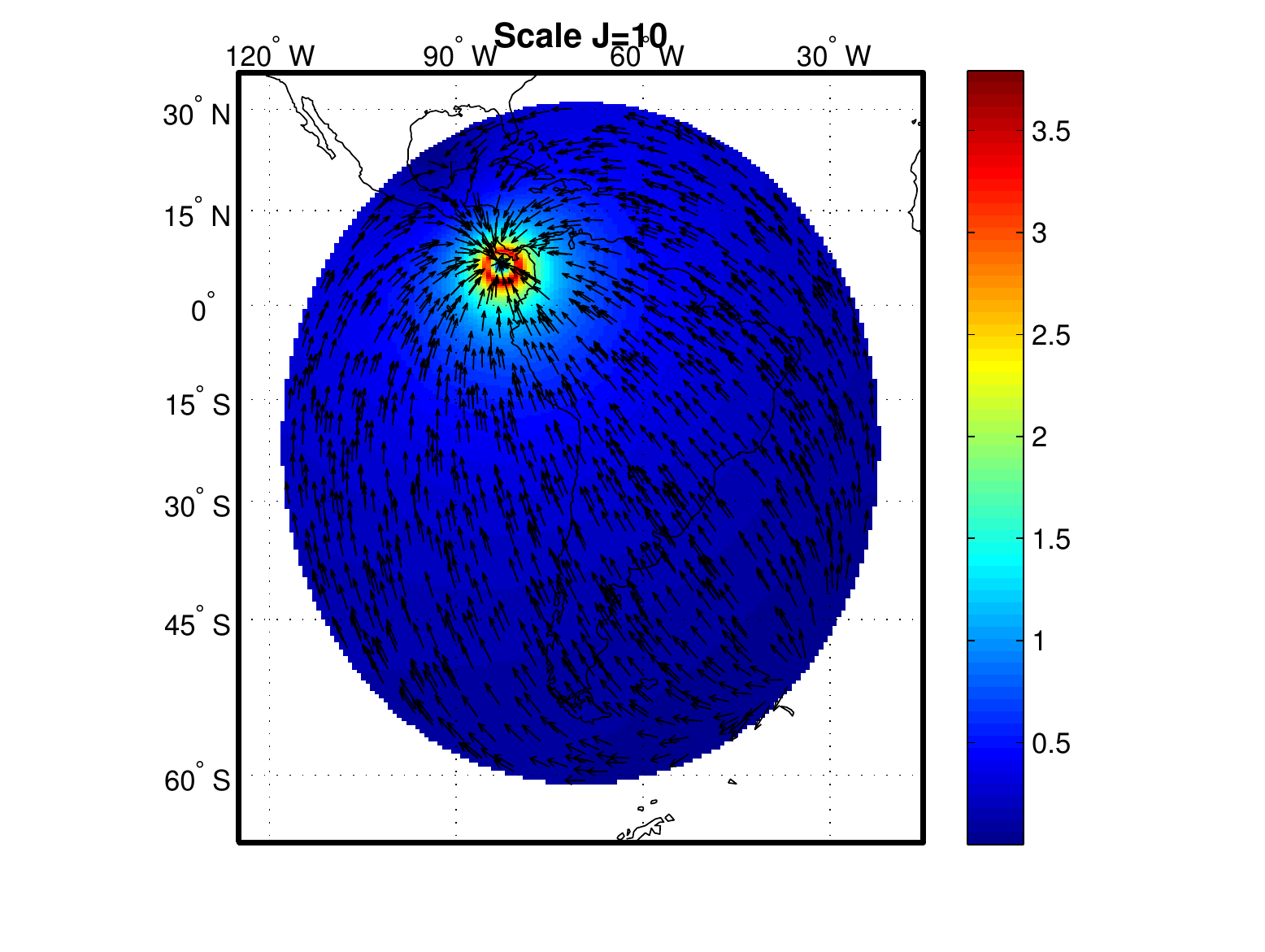}
\end{center}
\caption{The kernel $\nabla^*G_N^J(\Delta^*;\xi,\cdot)$ for scales $J=4,6,10$ and a fixed evaluation point $\xi$ located at $7^\circ$N, $74^\circ$W (colors indicate the absolute value and arrows the orientation).}\label{fig:kernels}
\end{figure}

The Earth's gravity potential $W=U+T$ is typically split into a normal gravity potential $U$ corresponding to a reference ellipsoid $\mathcal{E}$ (i.e., $U(x)=$ const. for $x\in\mathcal{E}$) and a smaller remaining disturbing potential $T$. The vertical deflection $\Theta(x)$ measures the angular distance between the normal vector $\nu_{\mathcal{G}}(x)$ at a point $x$ on the geoid $\mathcal{G}$ (i.e., $W(x)=$ const. for $x\in\mathcal{G}$) and the corresponding ellipsoidal normal vector $\nu_{\mathcal{E}}(x)$ with respect to $\mathcal{E}$. Assuming that $\nu_{\mathcal{G}}-\nu_{\mathcal{E}}$ and  $\nu_{\mathcal{E}}$ are nearly orthogonal and that the deviation of the reference ellipsoid from a sphere is negligible, one can derive the following relation between the disturbing potential and the deflections of the vertical:
\begin{align}
 \nabla^*T(R\xi)=-\frac{GM}{R}\Theta(R\xi),\quad\xi\in\Omega,\label{eqn:pdevertdeflect}
\end{align}
where $R$ is the Earth's mean radius, $G$ the gravitational constant, and $M$ the Earth's mass. For more details, the reader is referred to, e.g., \cite{freedenschreiner09, heiskanenmoritz67, wellenhof05}. We are particularly interested in solving \eqref{eqn:pdevertdeflect} for the disturbing potential $T$ in a subregion $\Gamma\subset\Omega$ (or, in other words, in a subregion $\Gamma_R$ of the spherical Earth's surface $\Omega_R$) from knowledge of the vertical deflections $\Theta$ only in that subregion. Theorem \ref{thm:gstarrepdgl} yields the representation
\begin{align}\label{eqn:solvevertdeflect}
T(R\xi)=\frac{1}{\|\Gamma\|}\int_\Gamma T(R\eta)d\omega(\eta)+\frac{GM}{R}\int_\Gamma\left(\nabla^*_\eta G_N(\Delta^*;\xi,\eta)\right)\cdot \Theta(R\eta)d\omega(\eta),\quad\xi\in\Gamma,
\end{align}
of which the first summand on the right hand side simply represents the constant mean disturbing potential $T^\Gamma_{mean}$ in $\Gamma_R$. We focus on the special case that $\Gamma=\Gamma_{\rho}(\zeta)$ is a spherical cap  with center $\zeta\in\Omega$ and radius $\rho\in(0,2)$, so that Theorem \ref{lem:neumanngreenfunccap} supplies us with an explicit representation of the Neumann Green function $G_N(\Delta^*;\cdot,\cdot)$. 

\begin{figure}
\begin{center}
\includegraphics[scale=0.42]{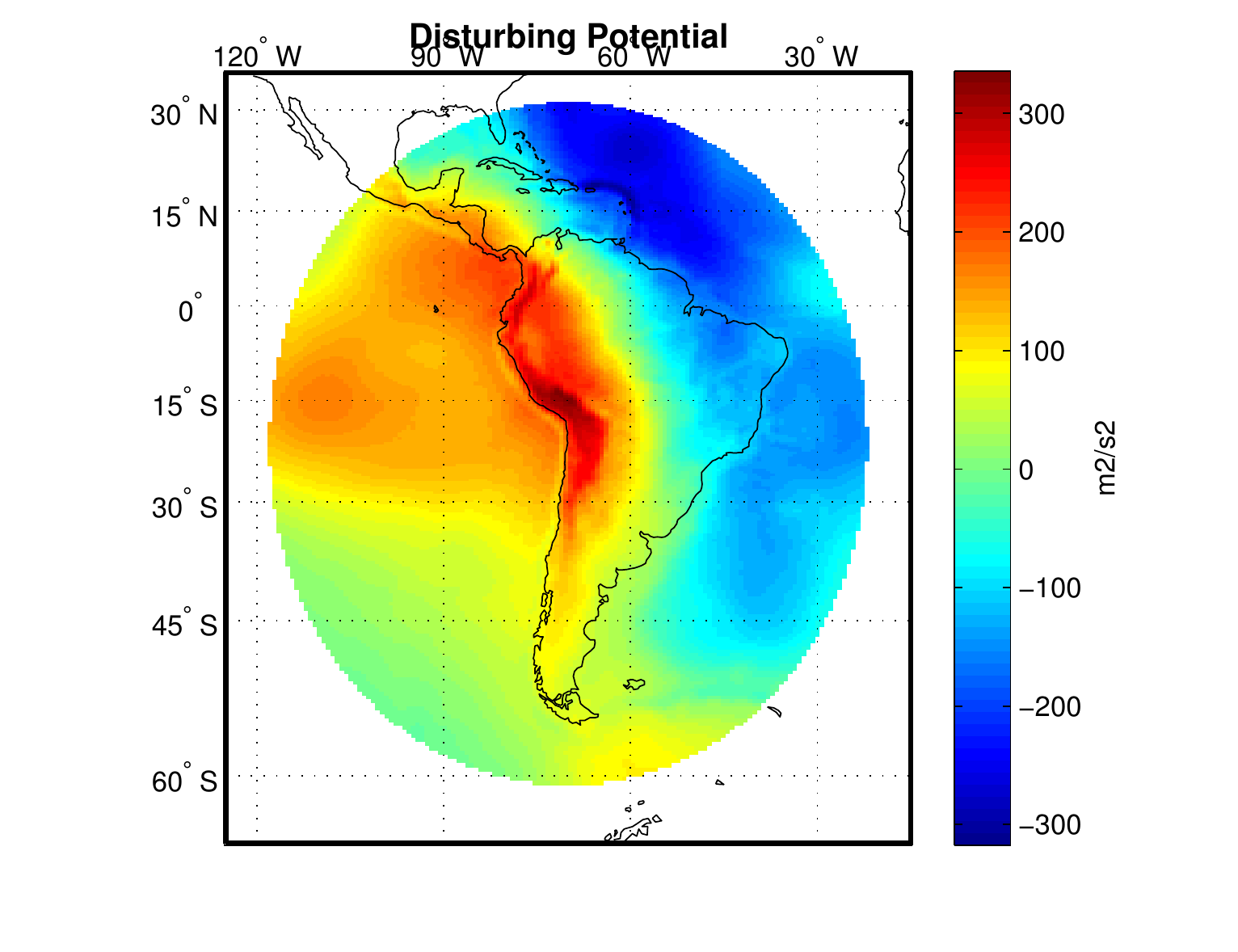}\quad \includegraphics[scale=0.42]{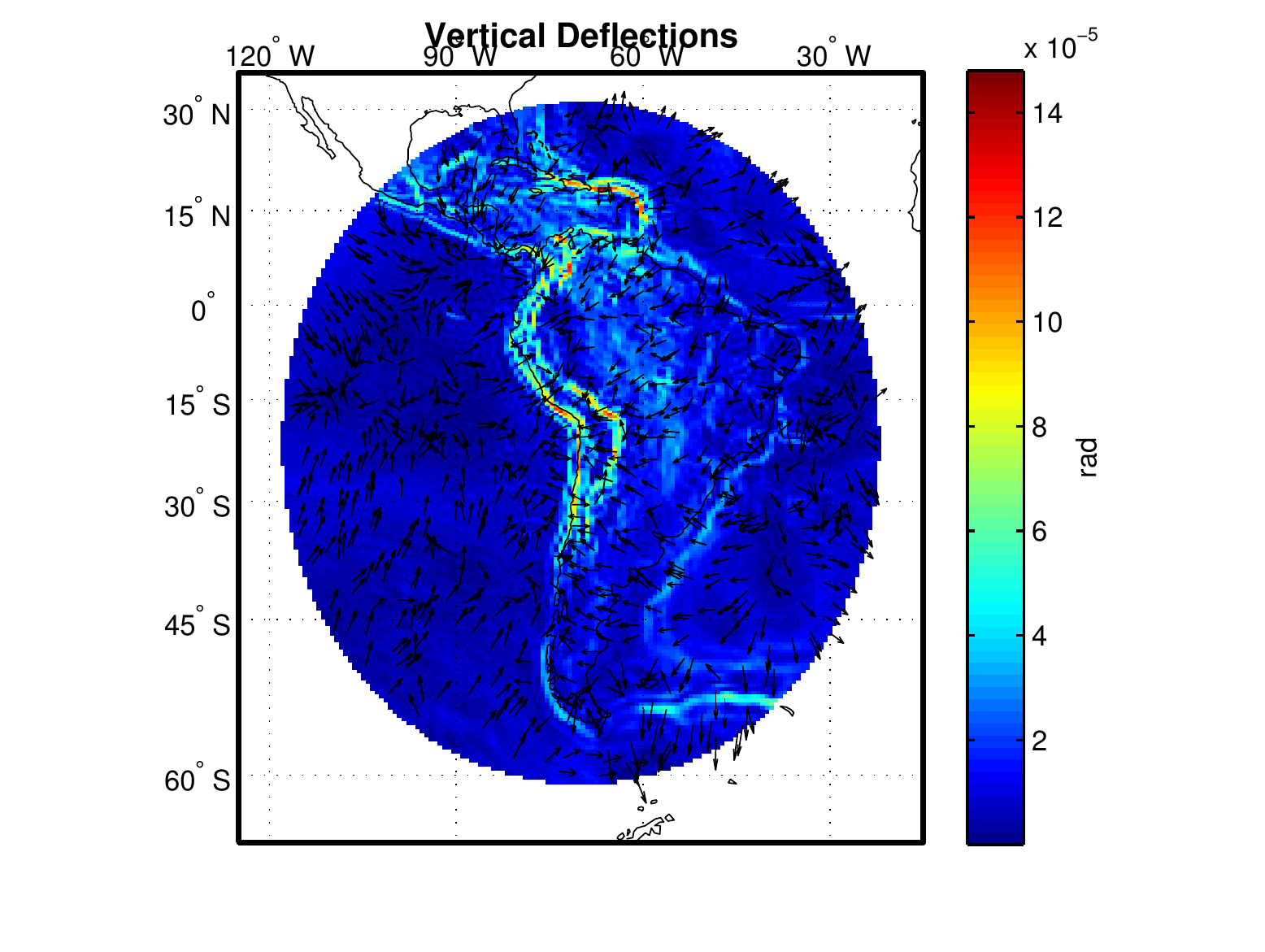}
\end{center}
\caption{The 'true' disturbing potential $T$ (left) and the corresponding vertical deflections $\Theta$ (right; colors indicate the absolute value and arrows the orientation).}\label{fig:vertdeflect}
\end{figure}

Concerning the numerical evaluation of \eqref{eqn:solvevertdeflect}, we first need to discretize the integral since $\Theta$ is typically only available in a discrete set of measurement points. For the tests in this section, we assume $\Theta$ to be given on a Gauss-Legendre grid in the spherical cap $\Gamma_\rho(\zeta)$, so that we can use the quadrature rule from \cite{hesse12}. Second, the numerical integration can become instable due to the singularity of the Neumann Green function $G_N(\Delta^*;\xi,\eta)$ at $\xi=\eta$ (originating in its contribution $\frac{1}{4\pi}\ln(1-\xi\cdot\eta)$). This can be circumvented by a regularization around this singularity via a truncated Taylor expansion. More precisely, for scaling parameters $ J=0,1,2,\ldots$, we define the regularized Neumann Green function 
\begin{align}\label{eqn:GNJ}
G_N^J(\Delta^*;\xi,\eta)=\left\{\begin{array}{lr}\frac{1}{4\pi}\ln(1-\xi\cdot\eta)+\frac{1}{4\pi}\ln(\check{r}(1-\check{\xi}\cdot\eta))
\\+\frac{1-\rho}{2\pi\rho}\ln(1+\zeta\cdot\eta),& 1-\xi\cdot\eta\geq 2^{-J},
                                 \\[1.35ex]\frac{2^J}{4\pi}(1-\xi\cdot\eta)-\frac{J}{4\pi}\ln(2)-\frac{1}{4\pi}
                                 \\+\frac{1}{4\pi}\ln(\check{r}(1-\check{\xi}\cdot\eta))+\frac{1-\rho}{2\pi\rho}\ln(1+\zeta\cdot\eta),&1-\xi\cdot\eta< 2^{-J}.
                                 \end{array}\right.
\end{align}
The regularization $G_N^J(\Delta^*;\cdot,\cdot)$ of the Neumann Green function $G_N(\Delta^*;\cdot,\cdot)$ closely relates to the regularization of the fundamental solution $G(\Delta^*;\cdot)$ briefly mentioned after Theorem \ref{thm:nablaintexchange}. A stable approximation of $T$ at scale $J$ is then given by
\begin{align}\label{eqn:solvevertdeflectJ}
T_J(R\xi)=T_{mean}^{\Gamma}+\frac{GM}{R}\int_{\Gamma_\rho(\zeta)}\left(\nabla^*_\eta G_N^J(\Delta^*;\xi,\eta)\right)\cdot \Theta(R\eta)d\omega(\eta),\quad\xi\in\Gamma_\rho(\zeta),
\end{align}
and satisfies $\lim_{J\to\infty}\sup_{\xi\in\tilde{\Gamma}}|T_J(R\xi)-T(R\xi)|=0$ for every subset $\tilde{\Gamma}\subset\Gamma_\rho(\zeta)$ with dist$(\tilde{\Gamma},\partial\Gamma_\rho(\zeta))>0$. Thus, higher scales $J$ yield a more precise approximation of $T$ and the difference $T_{J+1}-T_J$ between two consecutive scales reveals features of more and more local origin. The kernel $\nabla^* G_N^J(\Delta^*;\xi,\cdot)$ is illustrated in Figure \ref{fig:kernels}.

\begin{figure}
\includegraphics[scale=0.42]{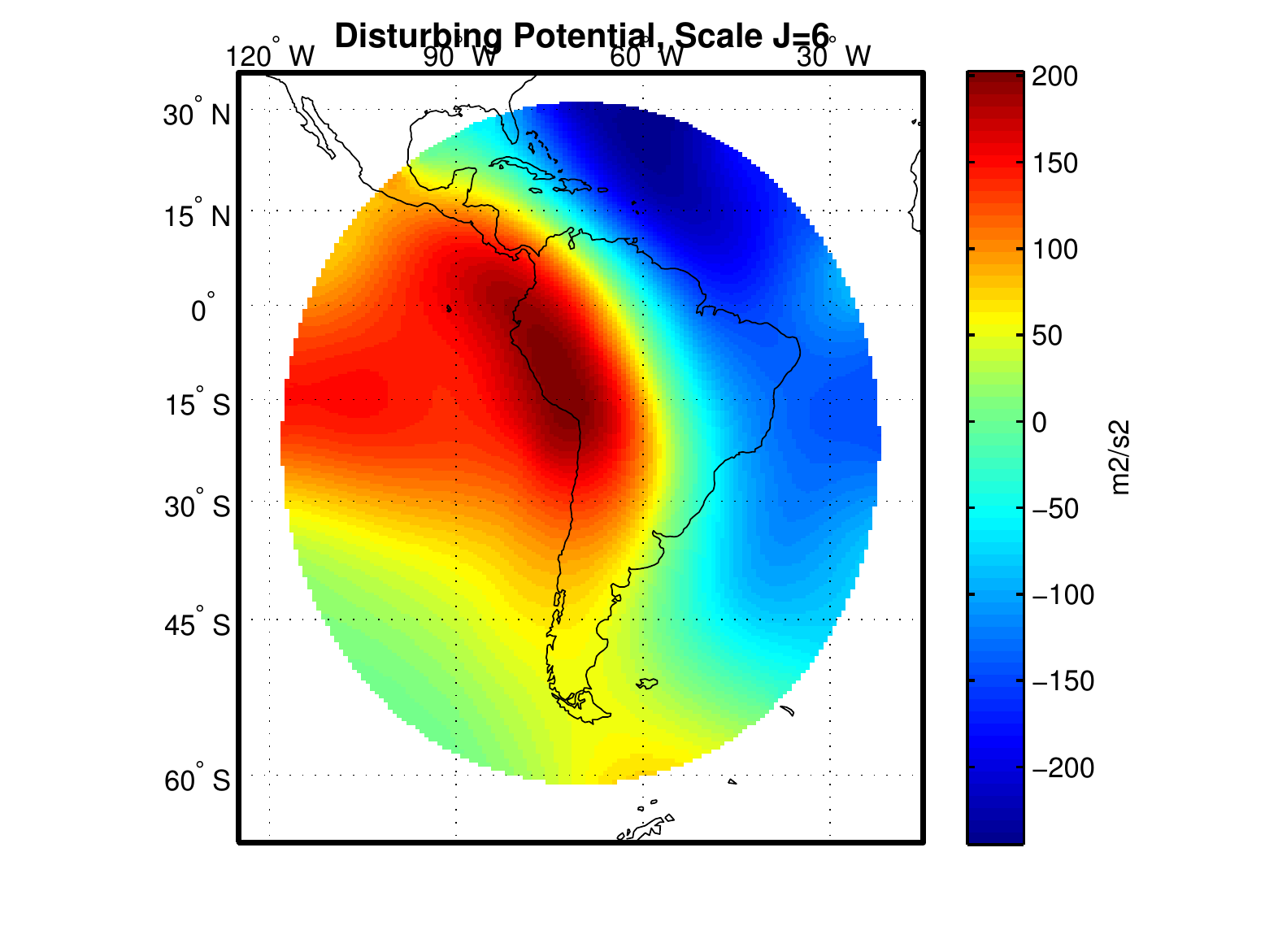}\quad \includegraphics[scale=0.42]{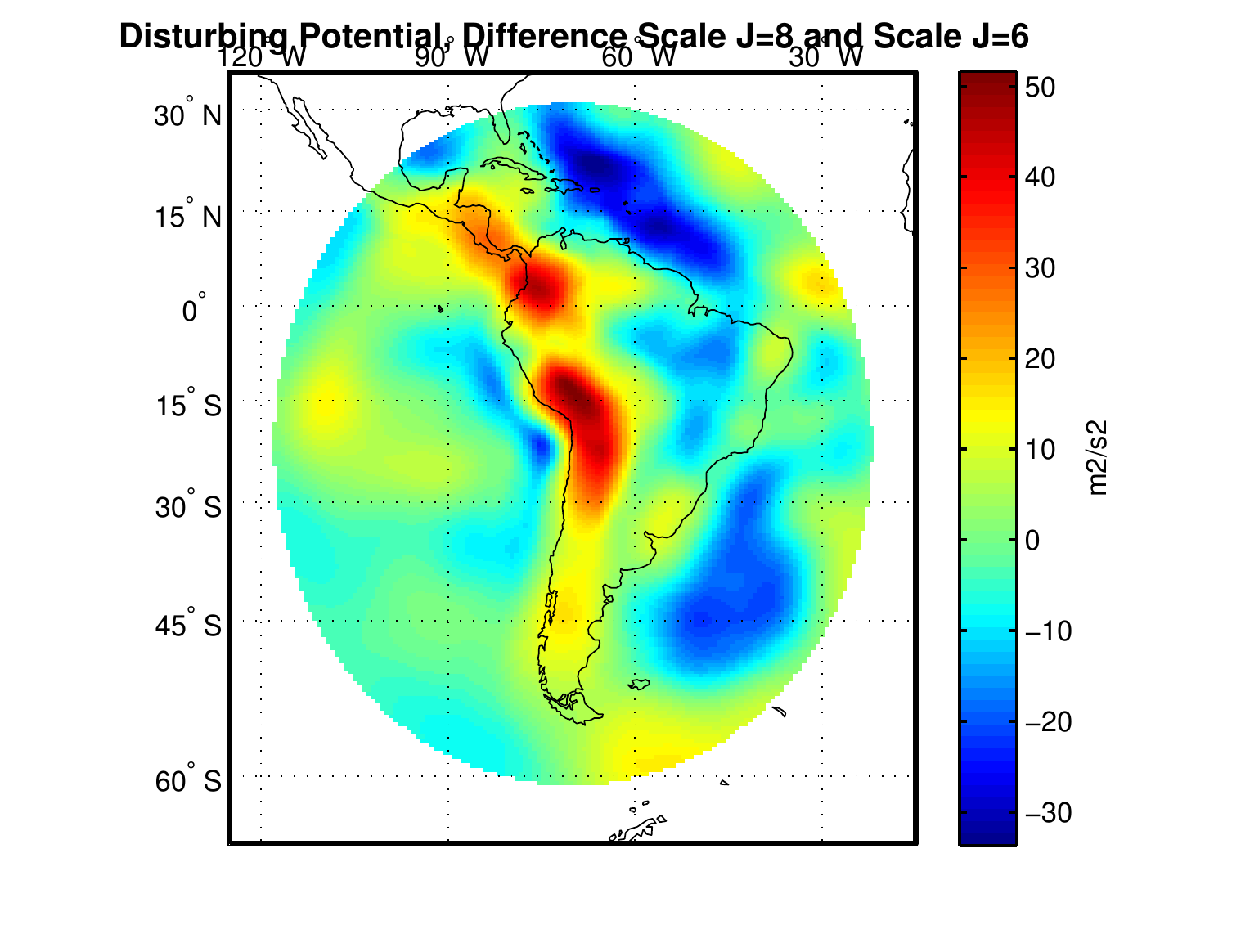}
\\\includegraphics[scale=0.42]{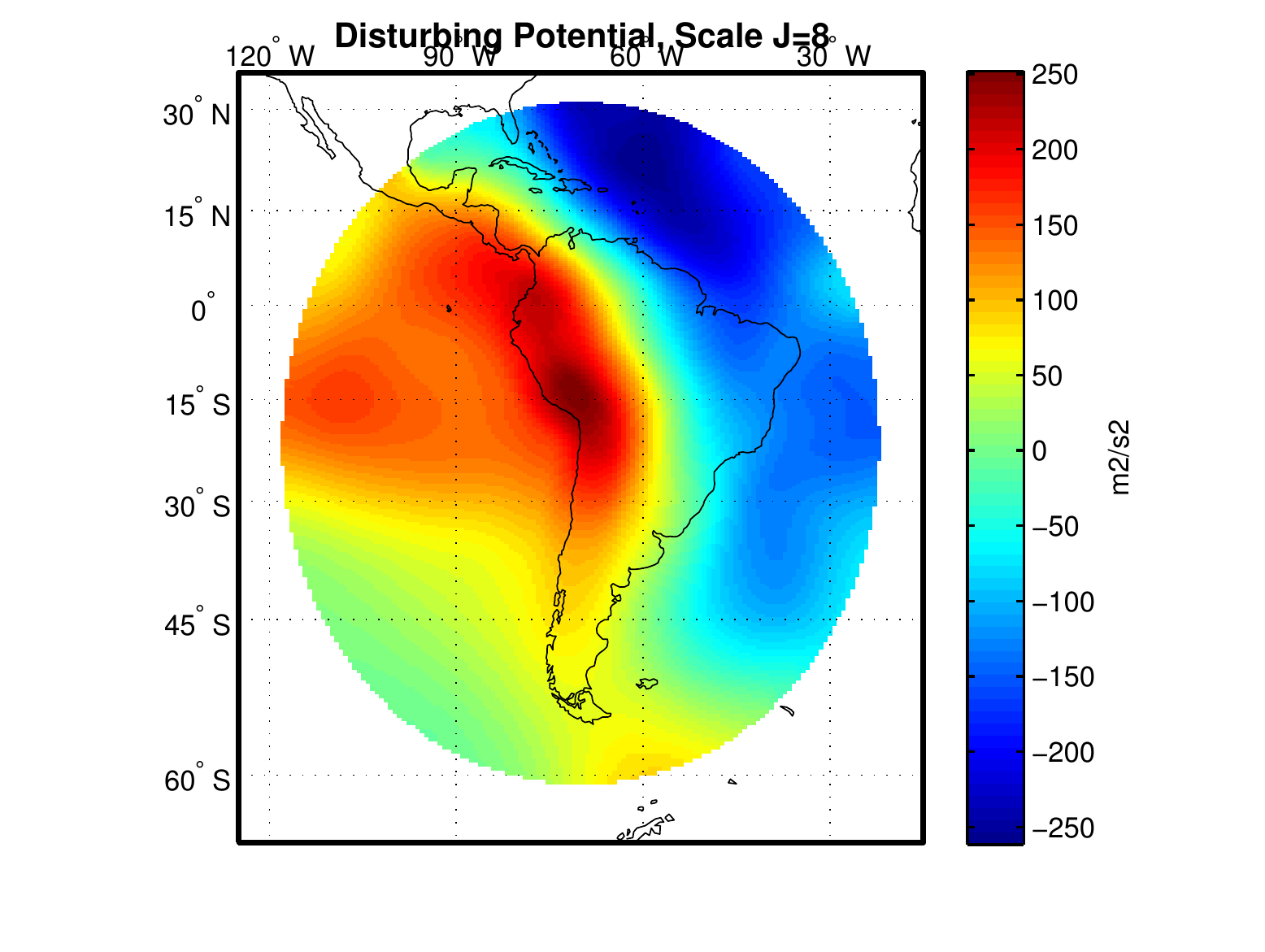}\quad \includegraphics[scale=0.42]{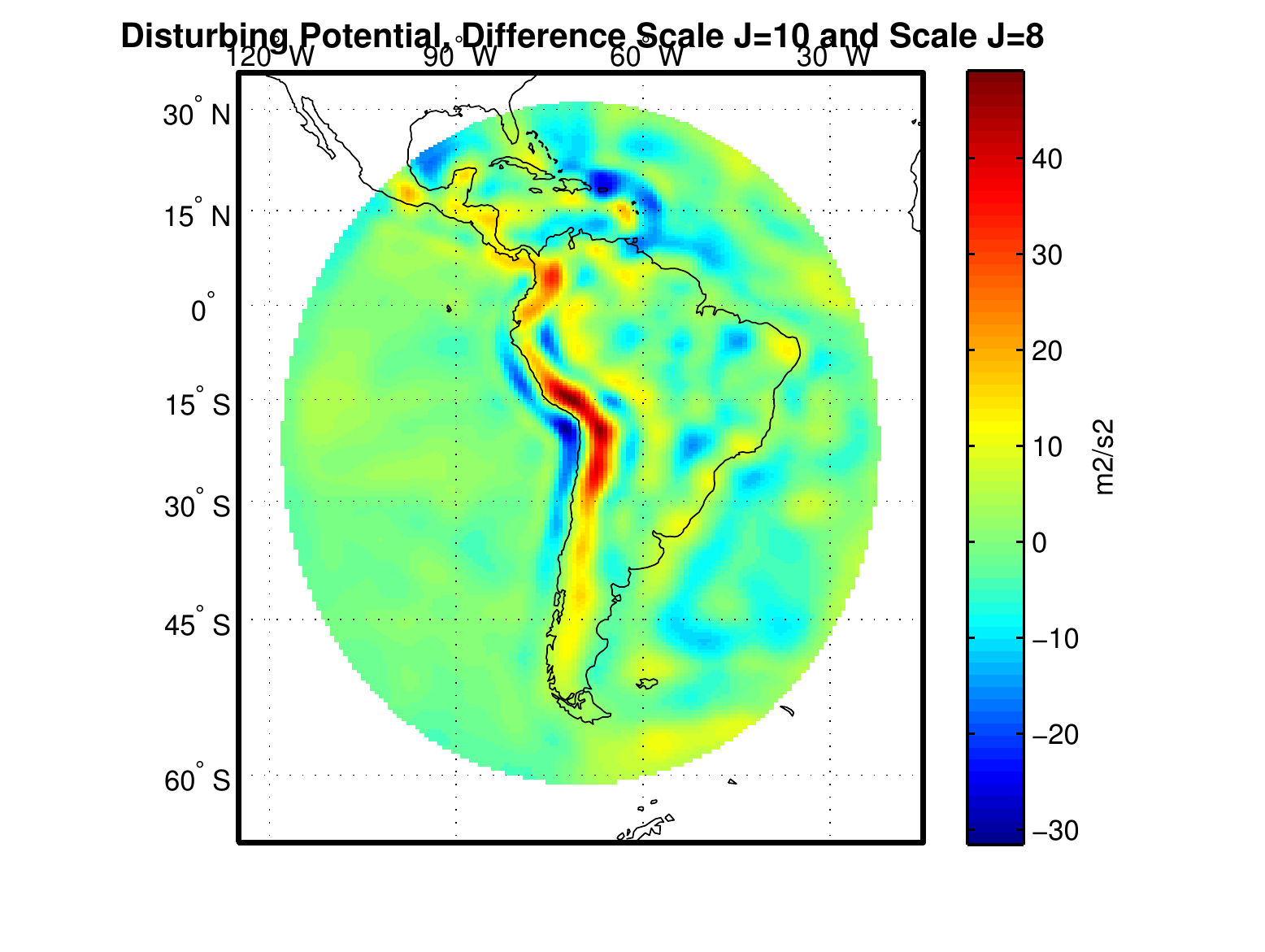}
\\\includegraphics[scale=0.42]{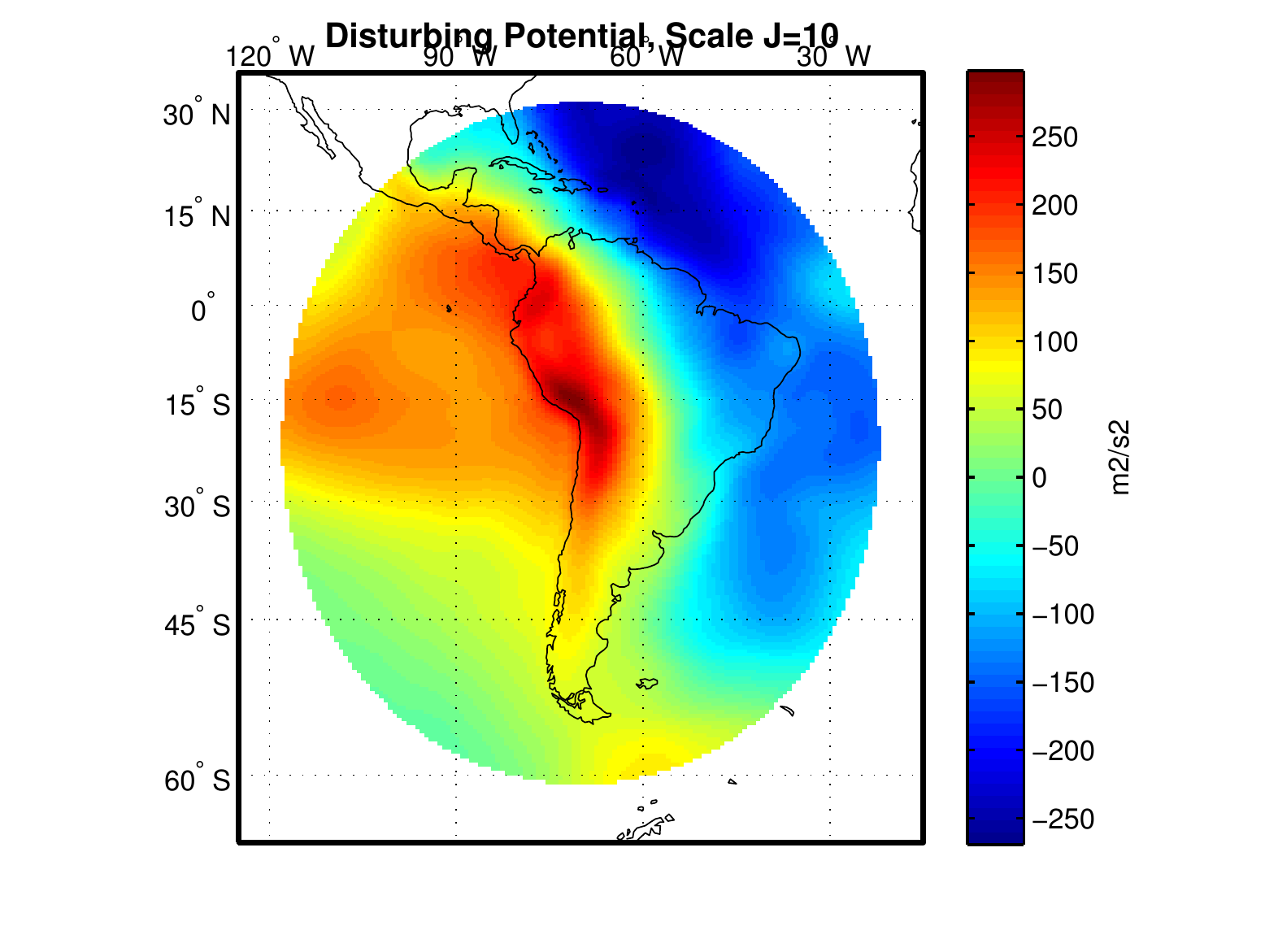}\quad \includegraphics[scale=0.42]{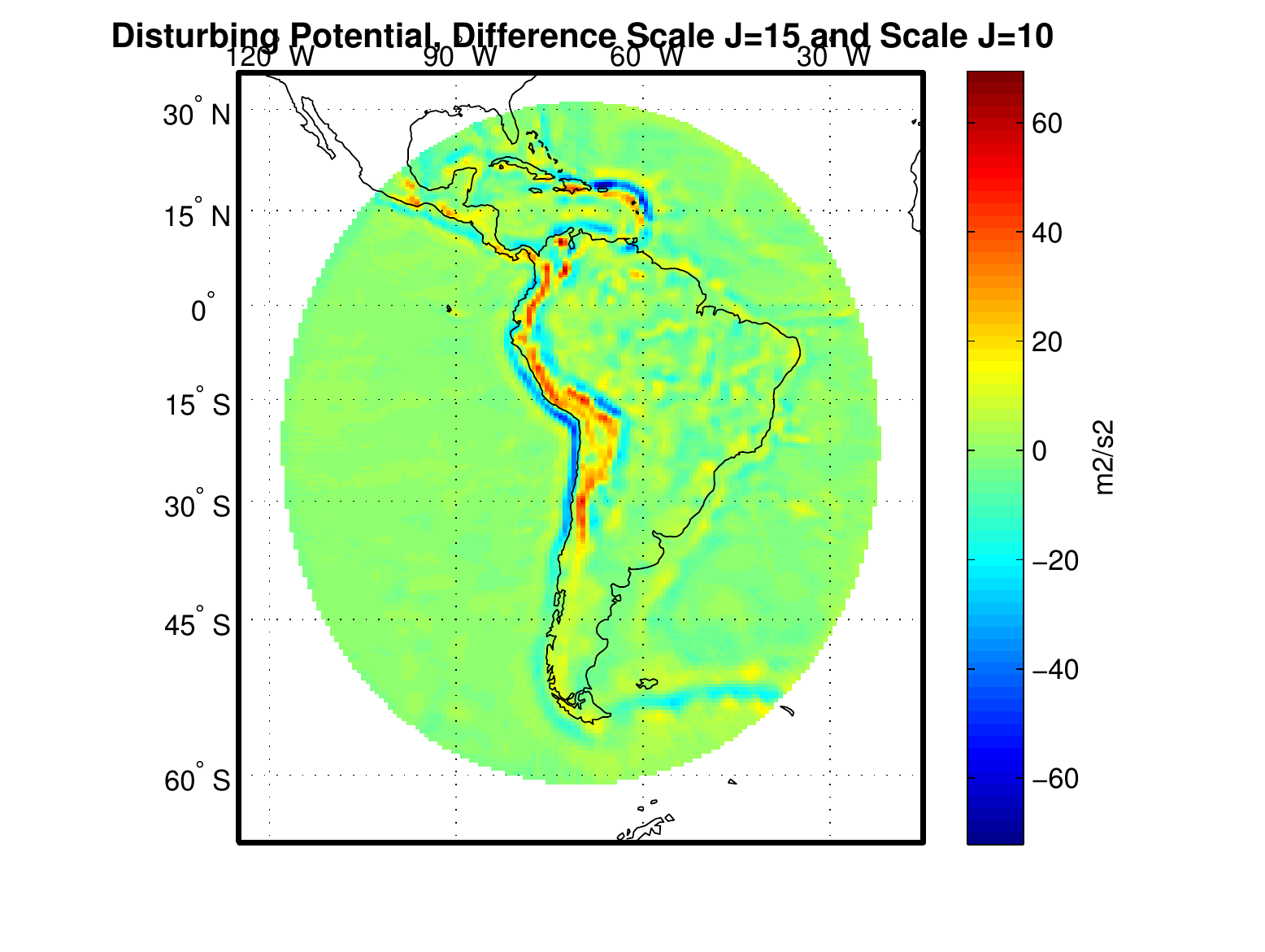}
\\\includegraphics[scale=0.42]{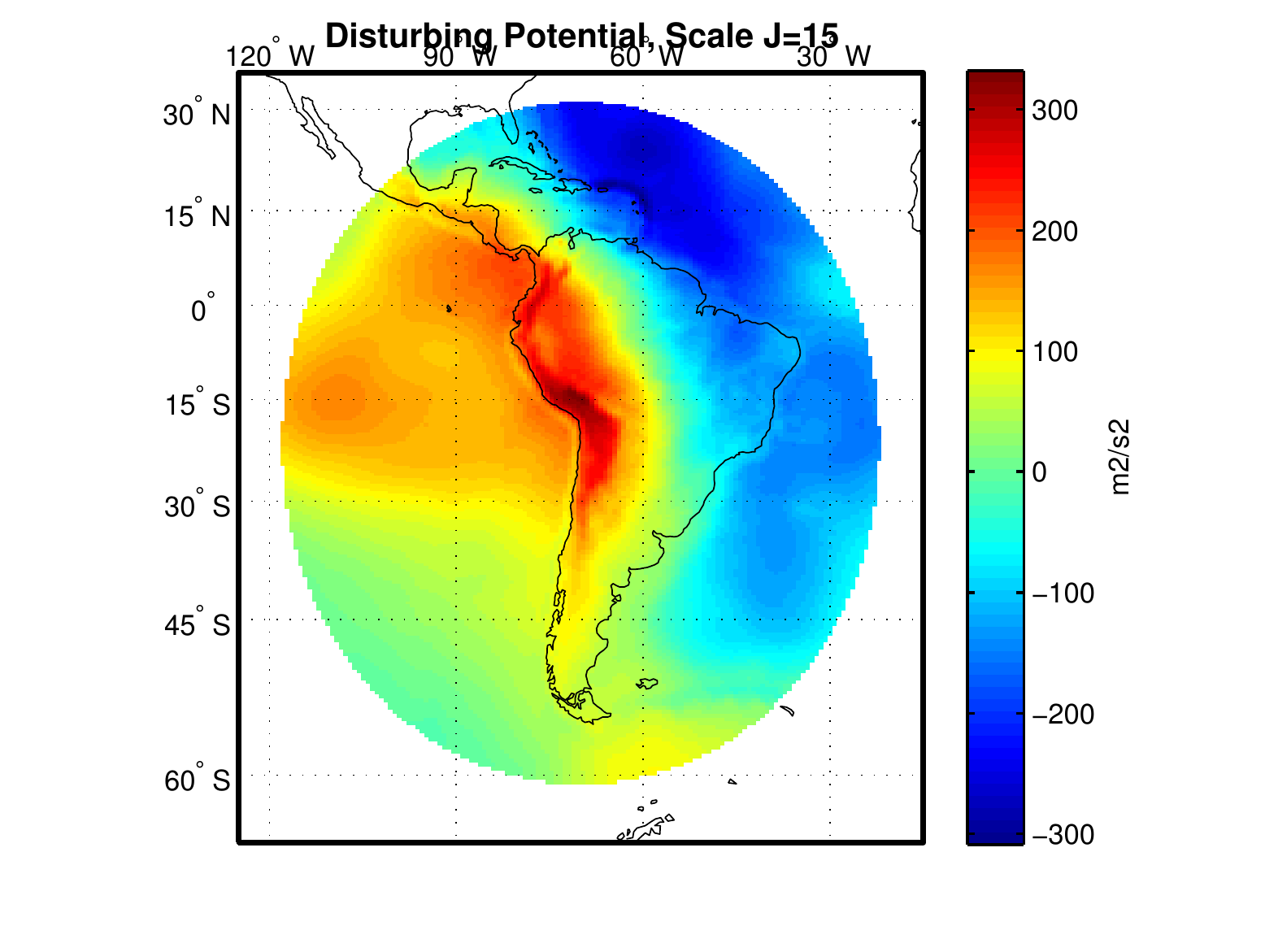}
\caption{Reconstructions of the disturbing potential $T_J$ at scales $J=6,8,10,15$ (left) and the differences $T_8-T_6$, $T_{10}-T_8$, $T_{15}-T_{10}$ between the reconstructions at these scales (right).}\label{fig:reconstdistpot}
\end{figure}

\begin{figure}
\includegraphics[scale=0.42]{true_DistPot-eps-converted-to.pdf}\quad\includegraphics[scale=0.42]{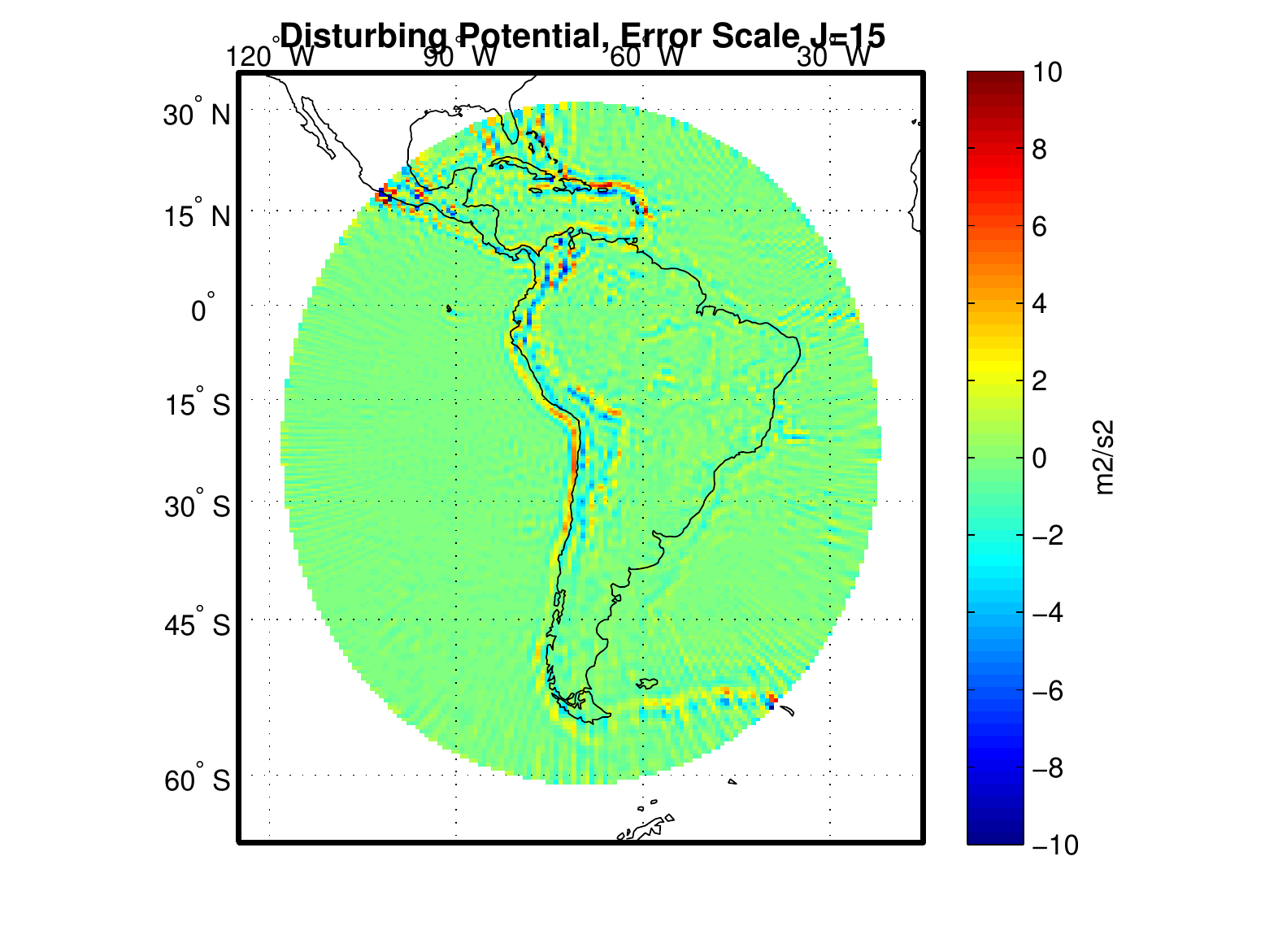}
\caption{The 'true' disturbing potential $T$ (left) and the reconstruction error $T-T_{15}$ (right).}\label{fig:truedistpot}
\end{figure}

In order to illustrate the reconstruction of the disturbing potential by the approximations $T_J$, we first compute a 'true' disturbing potential $T$ from EGM2008 (cf. \cite{egm2008}\footnote{data accessed via http://earth-info.nga.mil/GandG/wgs84/gravitymod/egm2008/egm08\_wgs84.html}) as a reference, using spherical harmonic degrees $n=3,\ldots,250$. From this $T$, we obtain our input vertical deflections $\Theta$ via \eqref{eqn:pdevertdeflect} on a Gauss-Legendre grid of 63,252 points in a spherical cap over South America (cf. Figure \ref{fig:vertdeflect}). The approximations $T_J$ for different scales $J$ are shown in Figure \ref{fig:reconstdistpot}. One can clearly see the refinement of the local features in the differences of the reconstructions $T_J$. Furthermore, the error $T-T_{15}$ in Figure \ref{fig:truedistpot} indicates a good approximation of $T$ and does not reveal any artefacts due to the local reconstruction without use of any boundary information.

\subsection{Geostrophic Ocean Flow}\label{sec:goc}

In subregions $\Gamma_R\subset\Omega_R$ of the ocean with a sufficiently large horizontal extent, away from the top and bottom Ekman layers and coastal regions, the geostrophic balance holds true: the horizontal pressure gradients in the ocean balance the Coriolis force resulting from horizontal currents. The Coriolis force term in a point $x\in\Gamma_R$ is given as the tangential contribution of $-2R\rho \,w\wedge v(x)$, where $v(x)$ is the horizontal ocean flow velocity and $w=|w|\varepsilon^3$ the Earth's rotation vector. $\rho$ denotes the density and is assumed to be constant. The pressure $P(x)$ in $x\in\Gamma_R$ can be regarded as being proportional to the mean dynamic topography (MDT) $H(x)$, which denotes the height of the sea surface relative to the Geoid $\mathcal{G}$ and can be determined from altimetry measurements.  More precisely, $P(x)=\rho G H(x)$, where $G$ denotes the gravitational constant. Using the geostrophic balance, we therefore obtain
\begin{align}
 -2R\rho(w\cdot\xi)\,\xi\wedge v(R\xi)=\rho G \nabla^* H(R\xi),\qquad \xi\in\Gamma,
\end{align}
or, equivalently,
\begin{align}\label{eqn:geoflow}
 \frac{2R}{G}|w|(\xi\cdot\varepsilon^3) v(R\xi)=\L^* H(R\xi),\qquad \xi\in\Gamma.
\end{align}
For more details on the geophysical background, the reader is referred, e.g., to \cite{pedlosky79,stewart}. In order to compute the MDT $H$ from knowledge of the ocean flow velocity $v$ in $\Gamma$, we need to solve Equation \eqref{eqn:geoflow}. Theorem \ref{thm:gstarrepdgl} yields the representation
\begin{align*}
H(R\xi)=\frac{1}{\|\Gamma\|}\int_\Gamma H(R\eta)d\omega(\eta)-\frac{2R}{G}|w|\int_\Gamma (\eta\cdot\varepsilon^3)\left(\L^*_\eta G_N(\Delta^*;\xi,\eta)\right)\cdot v(R\eta)  d\omega(\eta),\quad\xi\in\Gamma,
\end{align*}
of which the first summand on the right hand side simply represents the constant mean MDT $H^\Gamma_{mean}$ in $\Gamma_R$. Again, we focus on the special case that $\Gamma=\Gamma_{\rho}(\zeta)$ is a spherical cap  with center $\zeta\in\Omega$ and radius $\rho\in(0,2)$, so that we can apply the considerations from the previous section, i.e., we obtain an approximation at scale $J$ by 
\begin{align}\label{eqn:solvemdtJ}
H_J(R\xi)=H_{mean}^{\Gamma}-\frac{2R}{G}|w|\int_{\Gamma_\rho(\zeta)}(\eta\cdot\varepsilon^3)\left(\L^*_\eta G_N^J(\Delta^*;\xi,\eta)\right)\cdot v(R\eta) d\omega(\eta),\quad\xi\in\Gamma_\rho(\zeta),
\end{align}
where $G_N^J(\Delta^*;\cdot,\cdot)$ is given as in \eqref{eqn:GNJ}.

\begin{figure}
\begin{center}
\includegraphics[scale=0.42]{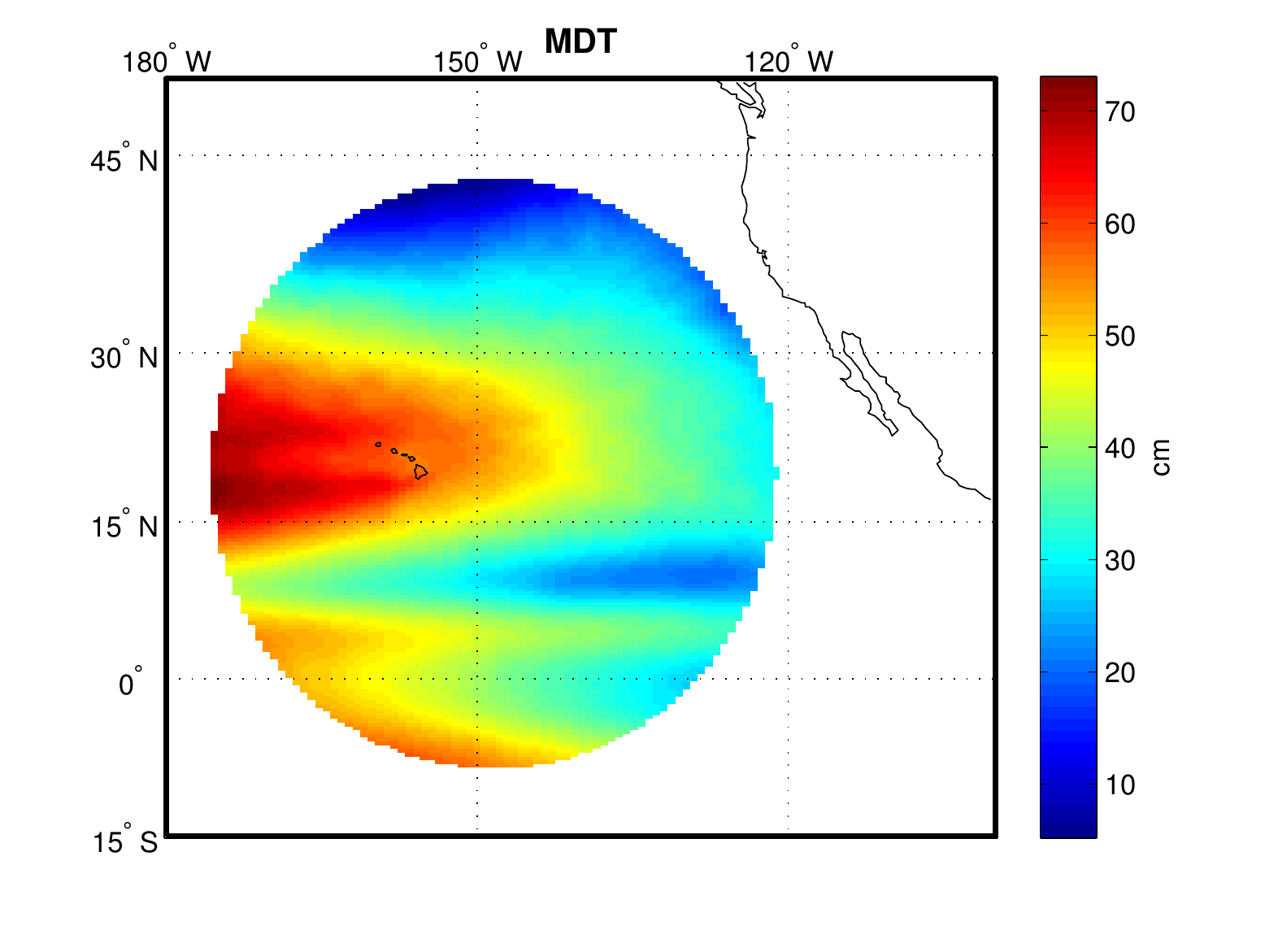}\quad \includegraphics[scale=0.42]{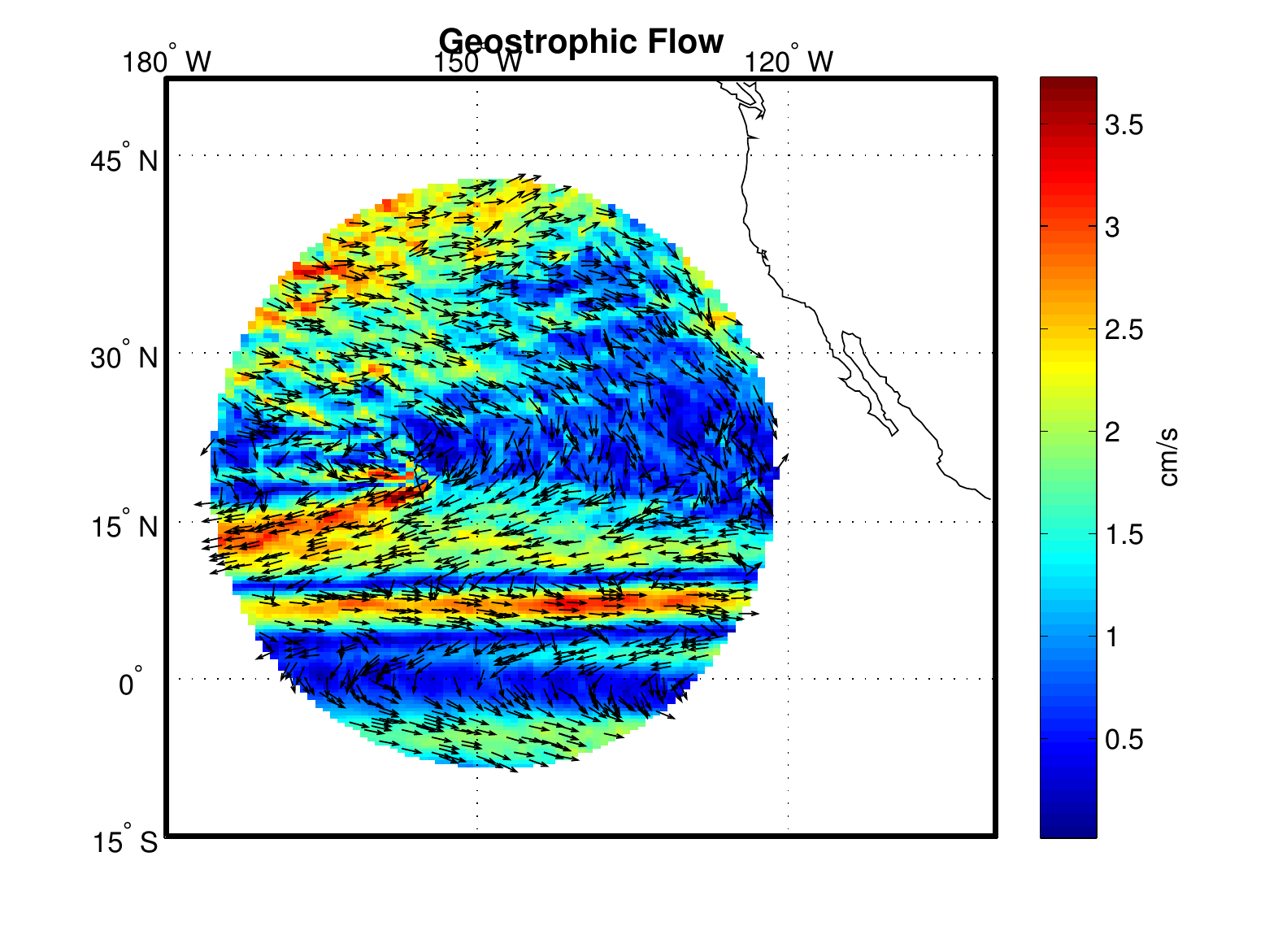}
\end{center}
\caption{The 'true' MDT $H$ (left) and the corresponding scaled geostrophic ocean flow velocity $\xi\mapsto(\xi\cdot\varepsilon^{3})v(\xi)$ (right; colors indicate the absolute values and arrows the orientation).}\label{fig:geoflow}
\end{figure}

\begin{figure}
\includegraphics[scale=0.42]{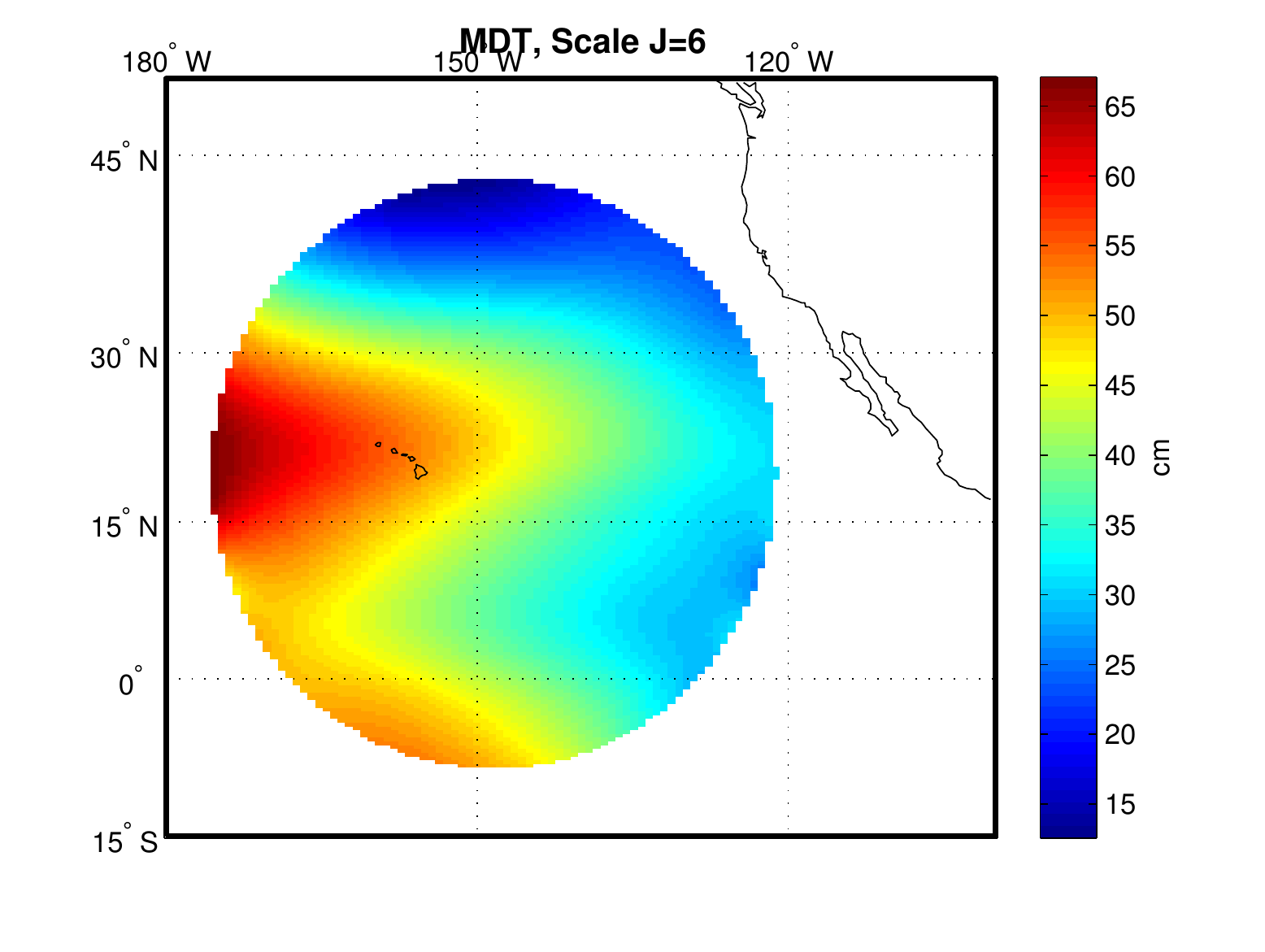}\quad \includegraphics[scale=0.42]{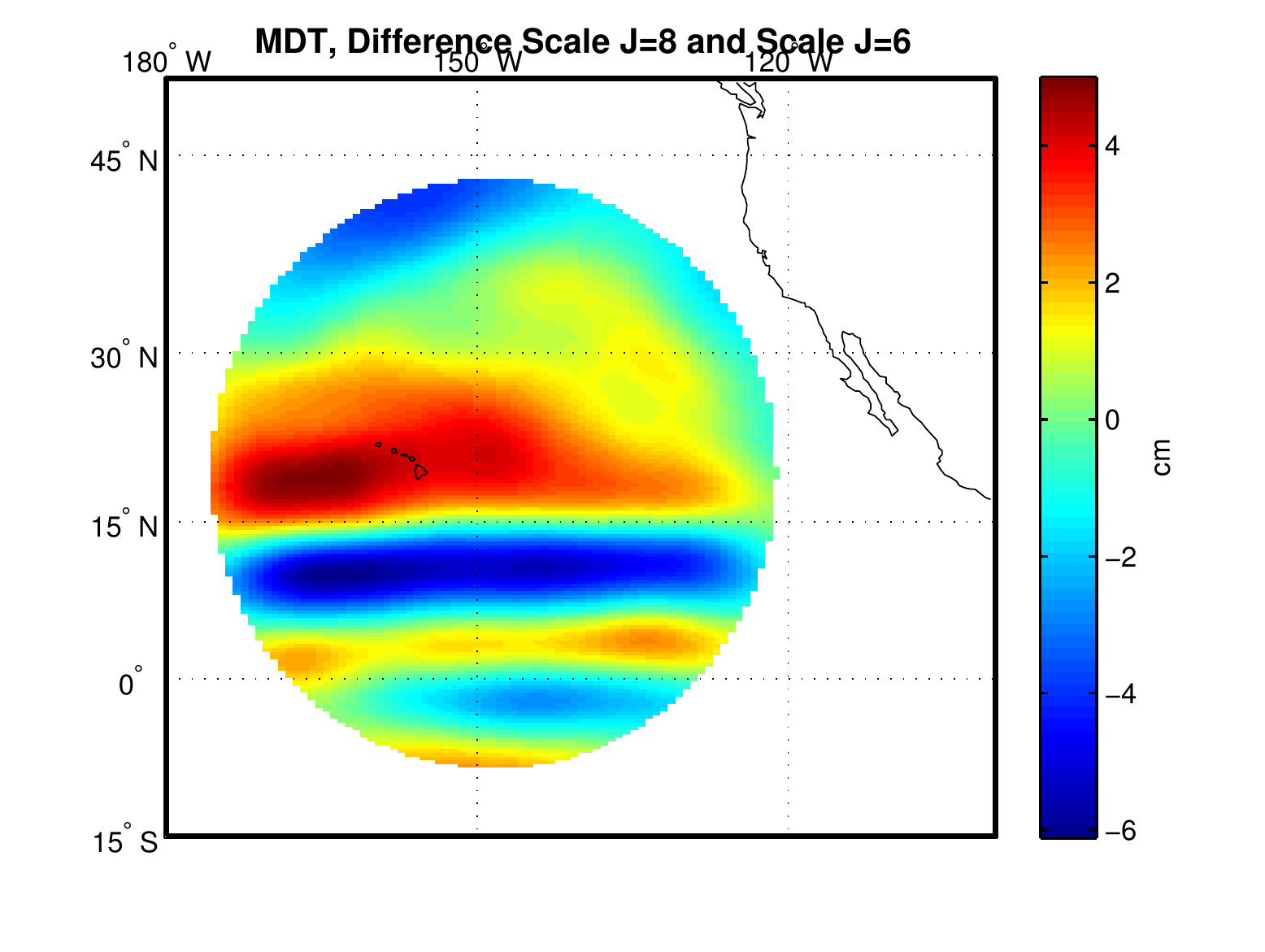}
\\\includegraphics[scale=0.42]{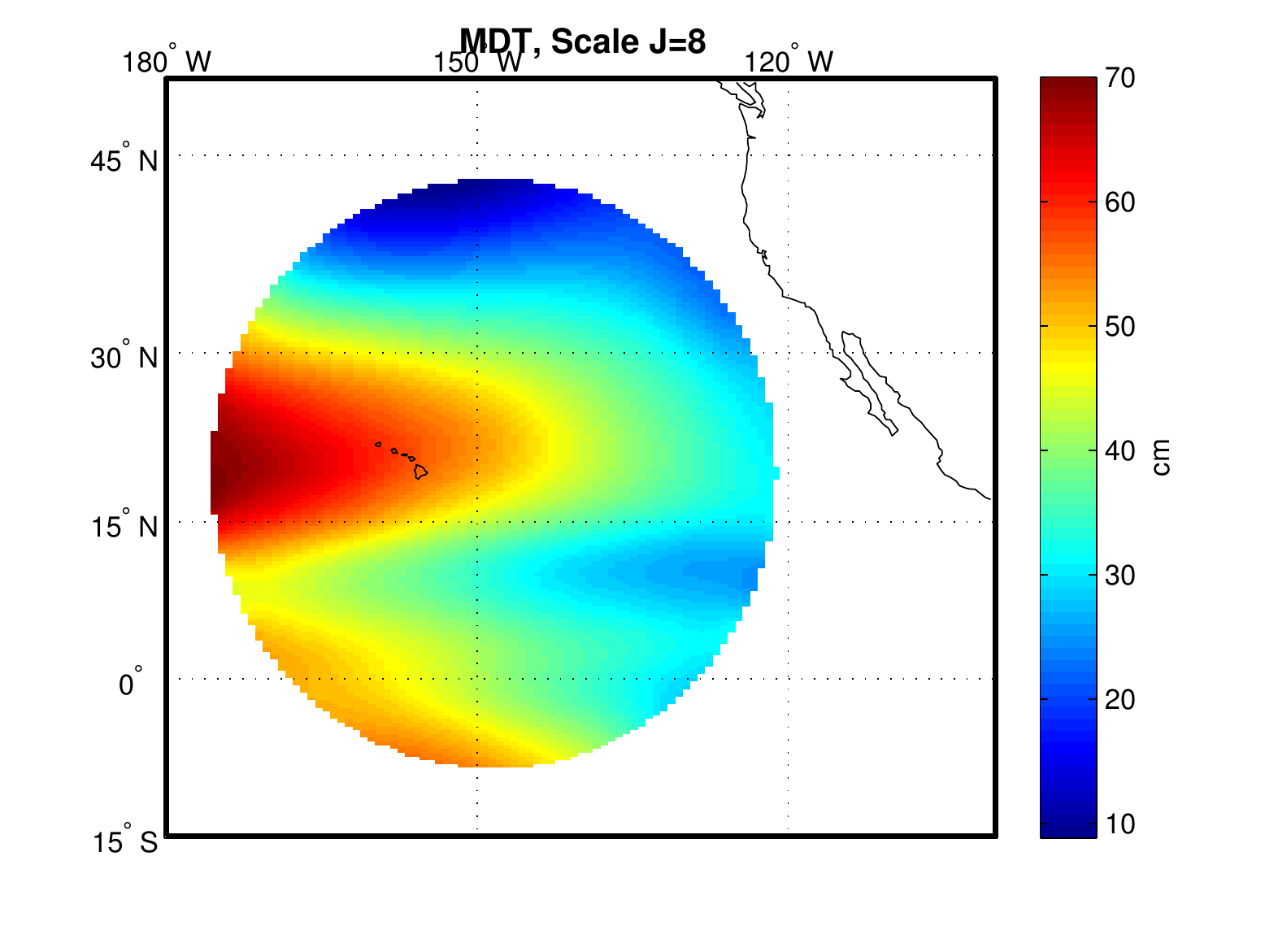}\quad \includegraphics[scale=0.42]{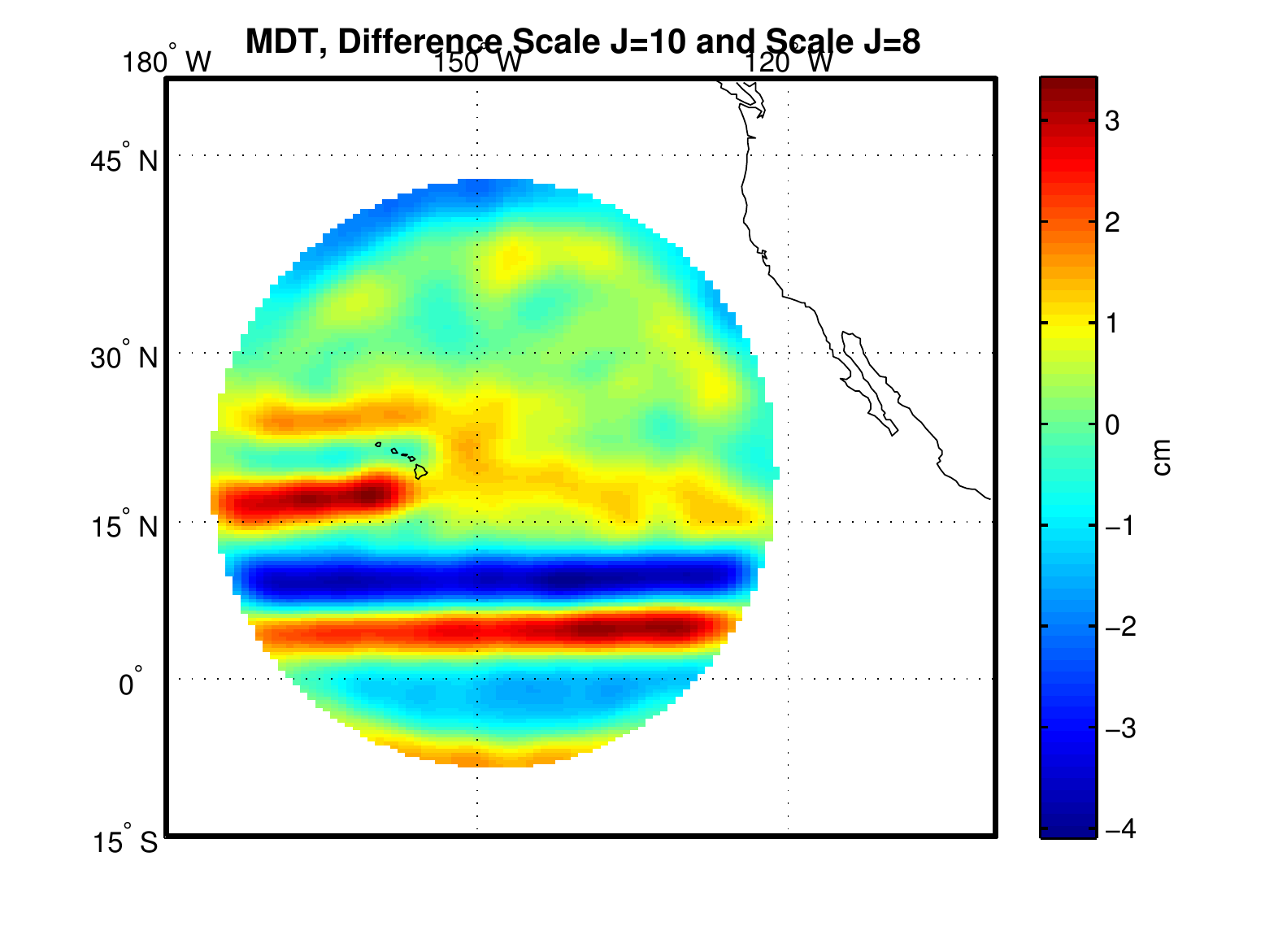}
\\\includegraphics[scale=0.42]{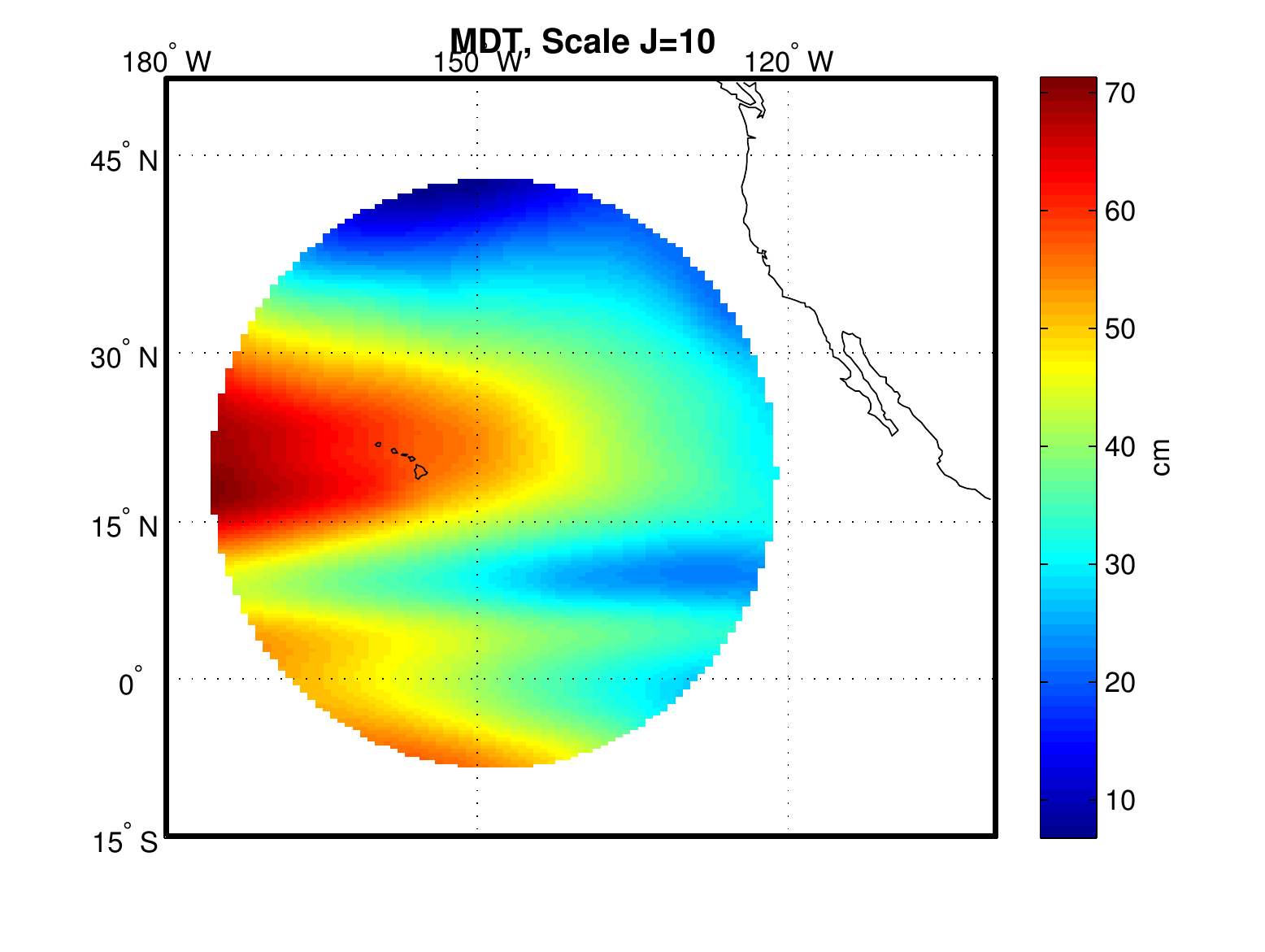}\quad \includegraphics[scale=0.42]{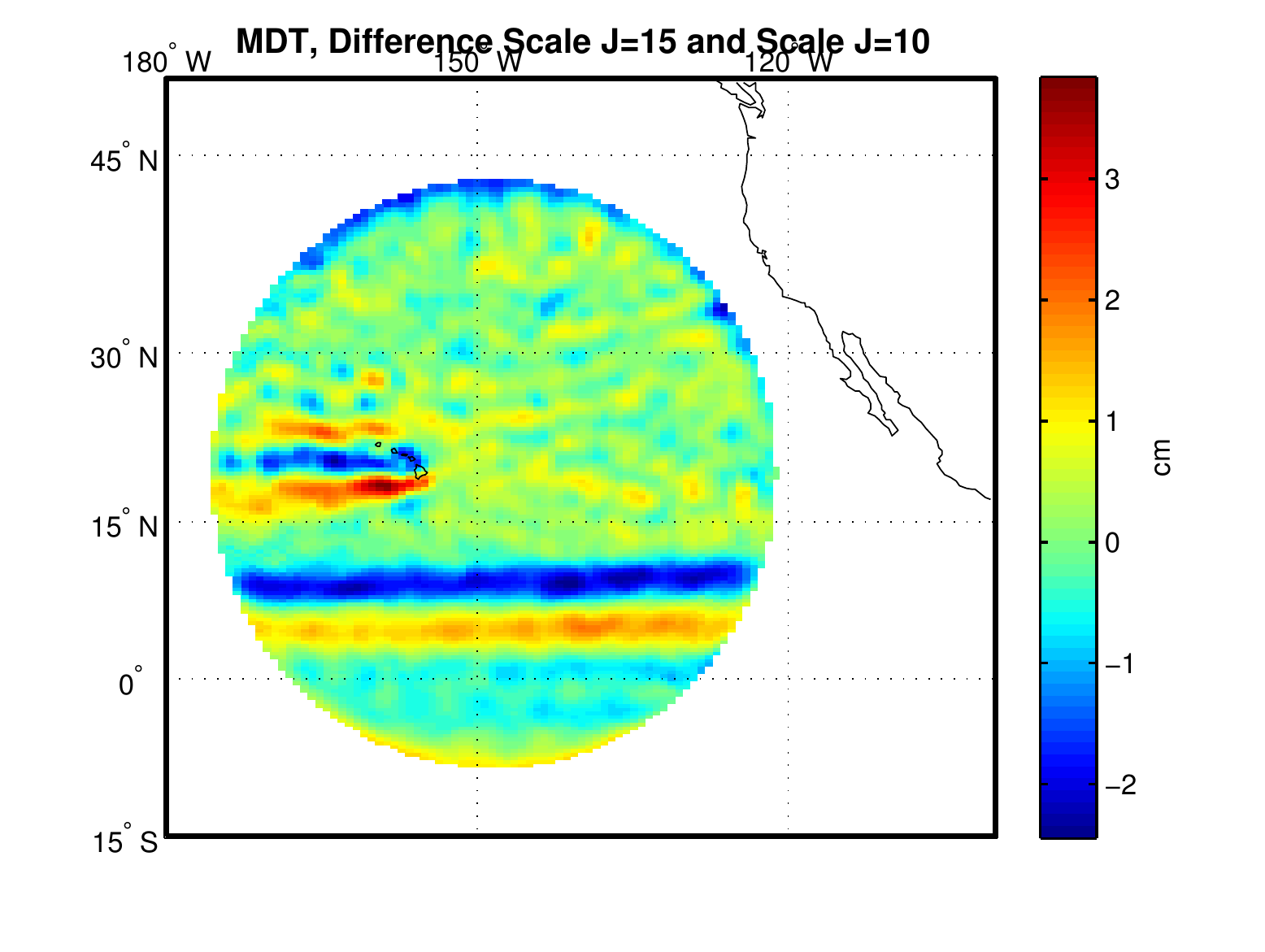}
\\\includegraphics[scale=0.42]{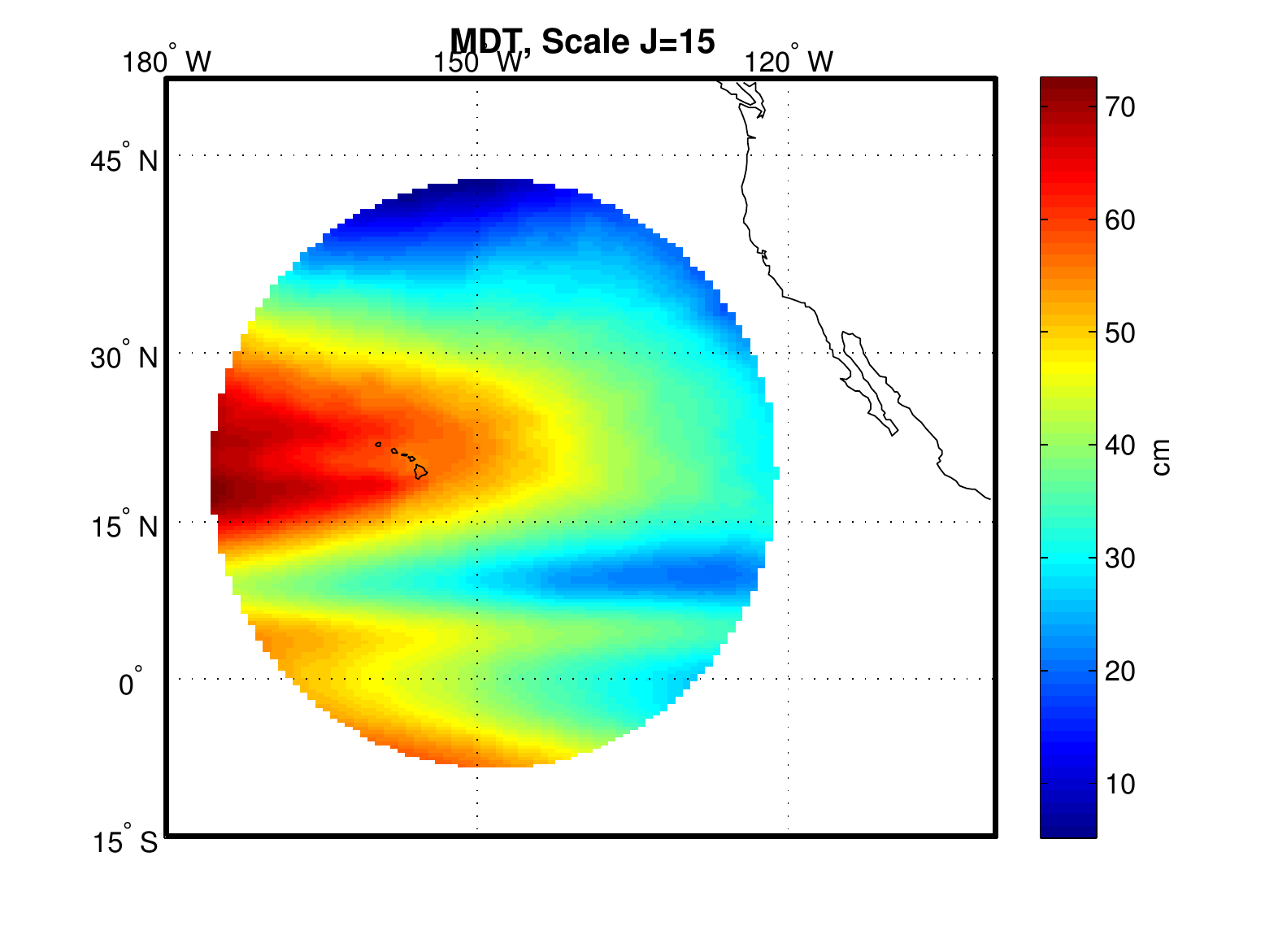}
\caption{MDT reconstructions $H_J$ at scales $J=6,8,10,15$ (left) and the differences $H_8-H_6$, $H_{10}-H_8$, $H_{15}-H_{10}$ between the reconstructions at these scales (right).}\label{fig:reconstmdt}
\end{figure}

\begin{figure}
\includegraphics[scale=0.42]{true_MDT-eps-converted-to.pdf}\quad\includegraphics[scale=0.42]{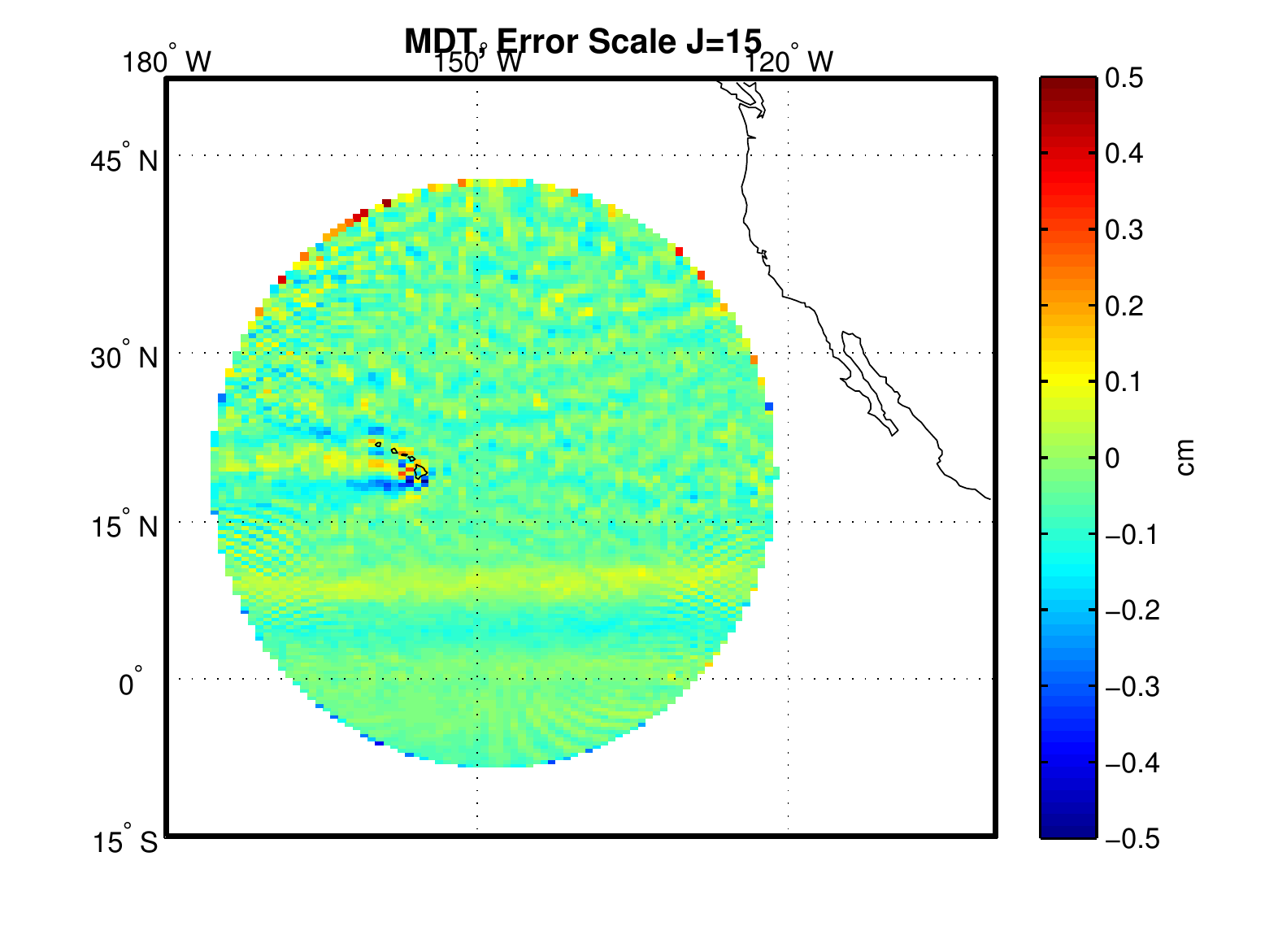}
\caption{The true MDT $H$ (left) and the reconstruction error $H-H_{15}$ (right).}\label{fig:truemdt}
\end{figure}

In order to illustrate the reconstruction of the MDT by the approximations $H_J$, we first compute a 'true' MDT $H$ from \cite{maximenko09}\footnote{data accessed via http://apdrc.soest.hawaii.edu/projects/DOT} as a reference. From this $H$, we can obtain our input ocean flow velocity $v$ via \eqref{eqn:geoflow} on a Gauss-Legendre grid of 63,252 points in a spherical cap over the Western Pacific Ocean (cf. Figure \ref{fig:geoflow}). The approximations $H_J$ for different scales $J$ are shown in Figure \ref{fig:reconstmdt}. The error $H-H_{15}$ in Figure \ref{fig:truemdt} indicates a good approximation of $H$ with larger errors only around the Hawaiian islands (where the geostrophic balance does not hold in the first place).

\subsection{Point Vortex Motion}\label{sec:vort}

Vorticity describes the rotational motion of a fluid. In the ocean, for horizontal flows $v$ which extend over regions $\Gamma_R\subset\Omega_R$ at spatial scales of several tens or hundreds of kilometers, the following relation for the vorticity $\omega$ holds true:
\begin{align}
 \omega(R\xi)=\L^*\cdot v(R\xi),\quad \xi\in\Gamma.
\end{align}
The incompressible horizontal flow $v$  itself can be represented by a stream function $\Psi$ via $v=\L^*\Psi$, so that we obtain
\begin{align}
 \omega(R\xi)= \Delta^*\Psi(R\xi),\quad \xi\in\Gamma.
\end{align}
The geostrophic flow from Section \ref{sec:goc} is an example for such a current. For more geophysical background on vorticity, the reader is again referred to \cite{pedlosky79, stewart}. 

A single point vortex at location $R{\eta}$ on the sphere is associated with a vorticity $\omega(R\xi)=\bar{\omega}\left(\delta(1-\xi\cdot{\eta})-\frac{1}{4\pi R}\right)$ of strength $\bar{\omega}\in\mathbb{R}$ (by $\delta$ we denote the Dirac distribution) and a corresponding stream function  $\Psi(R\xi)=\frac{\bar{\omega}}{R}G(\Delta^*;\xi\cdot{\eta})$, $\xi\in\Omega\setminus\{{\eta}\}$. If we consider a point vortex at location $R{\eta}$ in a subdomain $\Gamma_R\subset\Omega_R$ that produces no flow across the boundary $\partial\Gamma_R$ (e.g., a coastline), the vorticity would be $\omega(R\xi)=\bar{\omega}\left(\delta(1-\xi\cdot{\eta})\right)$ and the corresponding stream function $\Psi(R\xi)=\frac{\bar{\omega}}{R}G_D(\Delta^*;{\eta},\xi)$, $\xi\in\Gamma\setminus\{{\eta}\}$. In \cite{gemmrich08,kropinski14}, this motivated solving the model problem 
\begin{align}
\Delta^*\tilde{\Psi}(R\xi)&=0,\quad \xi\in\Gamma,\label{eqn:vort1}
\\\tilde{\Psi}^-(R\xi)&=\sum_{i=1}^N\frac{\bar{\omega}_i}{R}G(\Delta^*;\xi\cdot{\eta}_i)-\frac{\bar{\omega}_i}{4\pi R}\ln(1-\xi\cdot\bar{\xi}),\quad \xi\in\partial\Gamma,\label{eqn:vort2}
\end{align}
for a fixed $\bar{\xi}\in\Gamma^c$ and point vortices of strengths $\bar{\omega}_i$ located at $R\eta_i\in\Gamma_R$, $i=1,\ldots,N$. The actual stream function is then given by $\Psi(R\xi)=\sum_{i=1}^N\frac{\bar{\omega}_i}{R}G(\Delta^*;\xi\cdot{\eta}_i)-\frac{\bar{\omega}_i}{4\pi R}\ln(1-\xi\cdot\bar{\xi})-\tilde{\Psi}(R\xi)=\sum_{i=1}^N\frac{\bar{\omega}_i}{R}G_D(\Delta^*;{\eta}_i,\xi)$, $\xi\in\Gamma\setminus\{\eta_1,\ldots,\eta_N\}$. More details on point vortex motion on the entire sphere (and more general closed manifolds) can be found, e.g., in \cite{dritschel88,dritschel15, kidambi98,kimuraokamoto87}, and details on point vortex motion on subdomains of the sphere with impenetrable boundaries, e.g., in \cite{gutkinnewton04,kidambinewton00}.

\begin{figure}
\begin{center}
\includegraphics[scale=0.36]{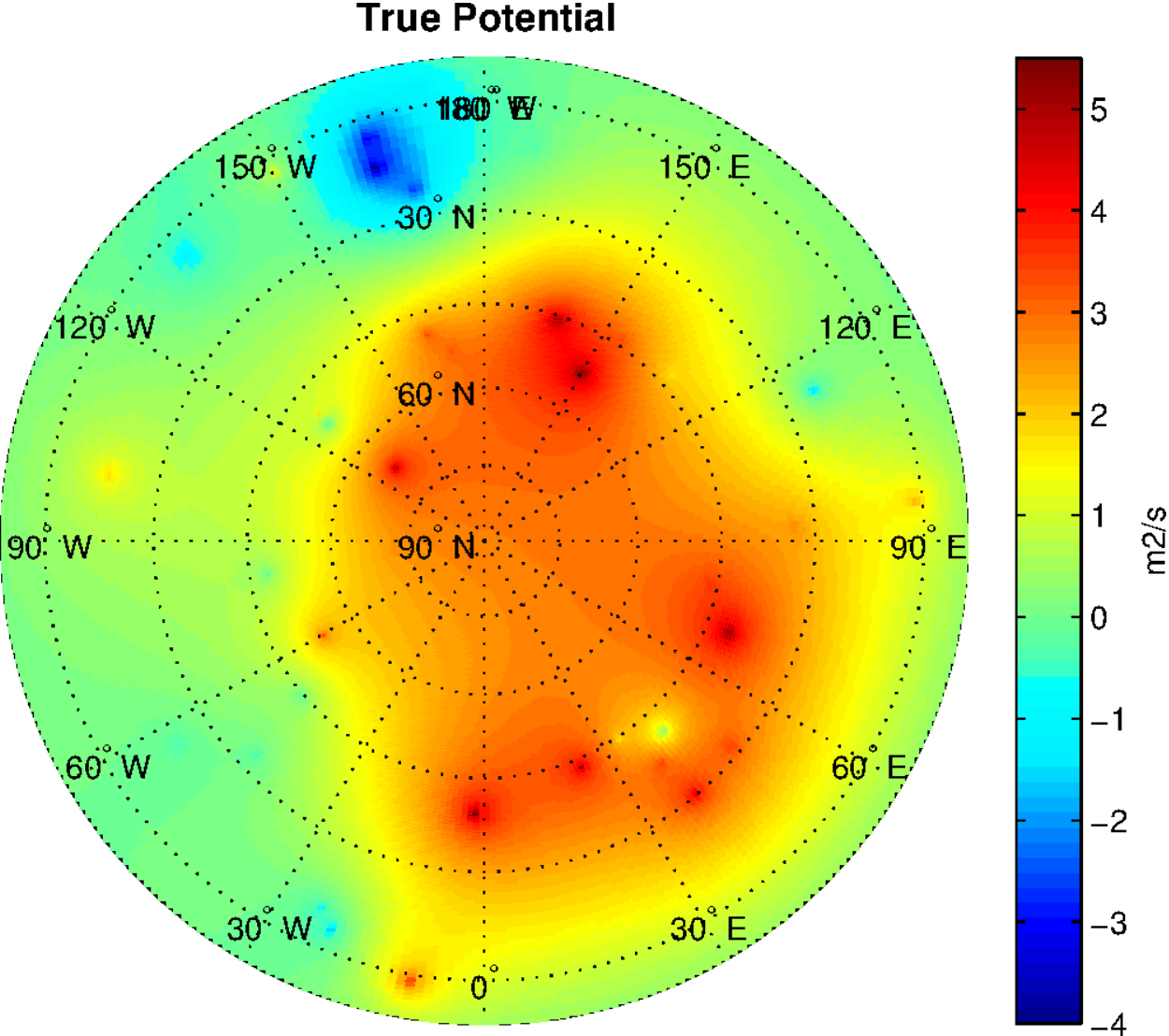}\quad\quad\quad\quad\quad \includegraphics[scale=0.42]{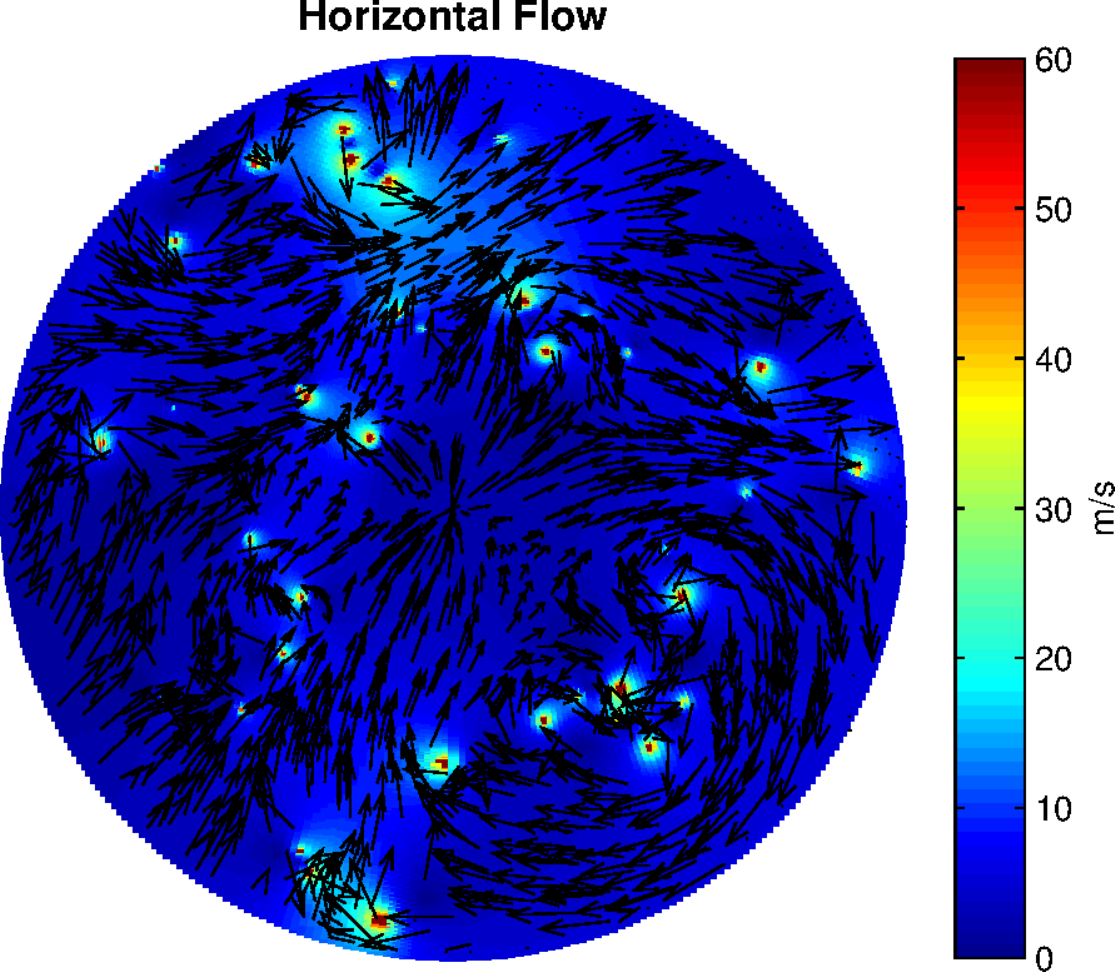}
\end{center}
\caption{The 'true' potential $\Psi$ (left) and the corresponding horizontal flow velocity $v$ (right; colors indicate the absolute values and arrows the orientation).}\label{fig:vortflow}
\end{figure}

\begin{figure}
\begin{center}
\includegraphics[scale=0.36]{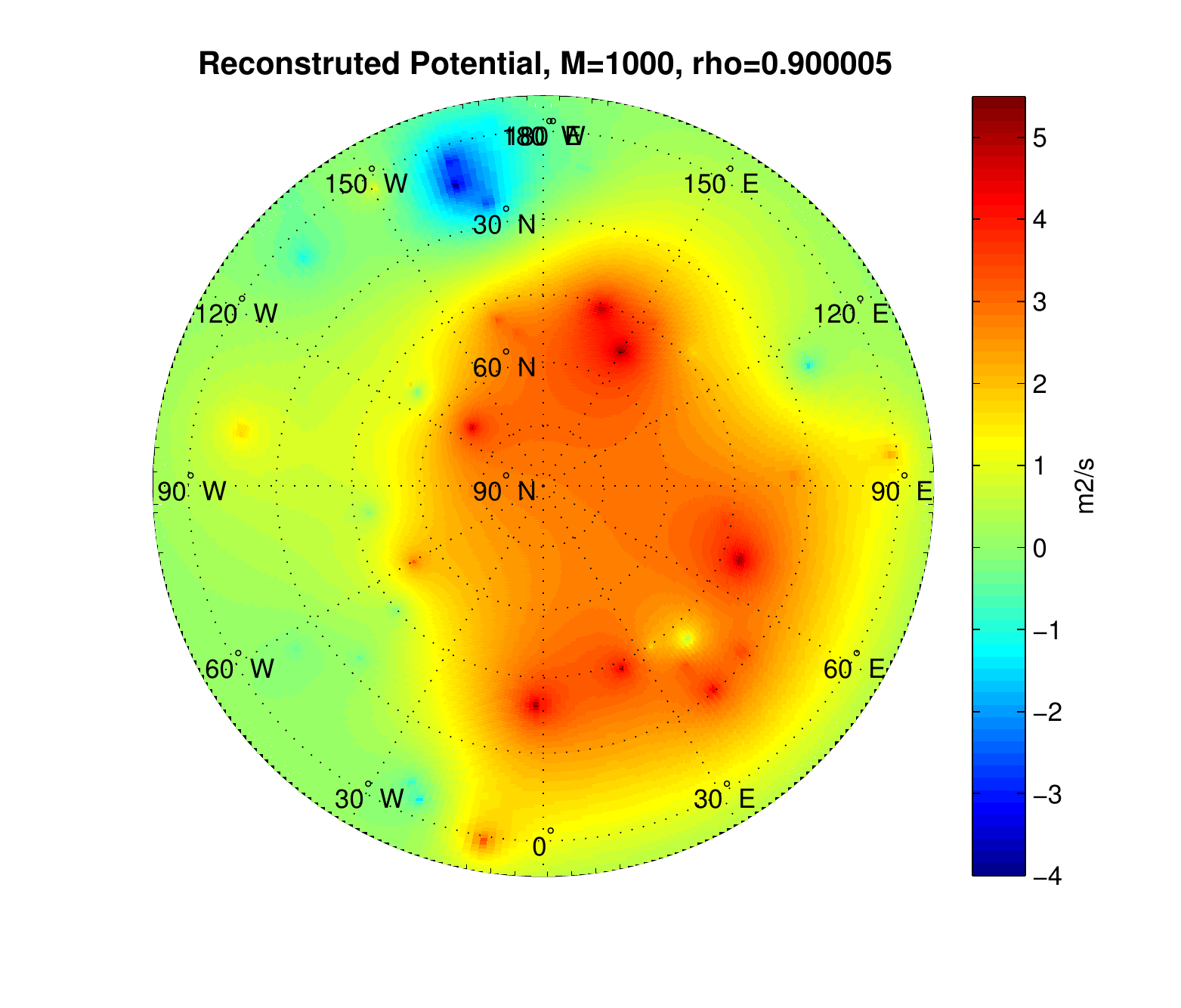}\quad\quad\includegraphics[scale=0.36]{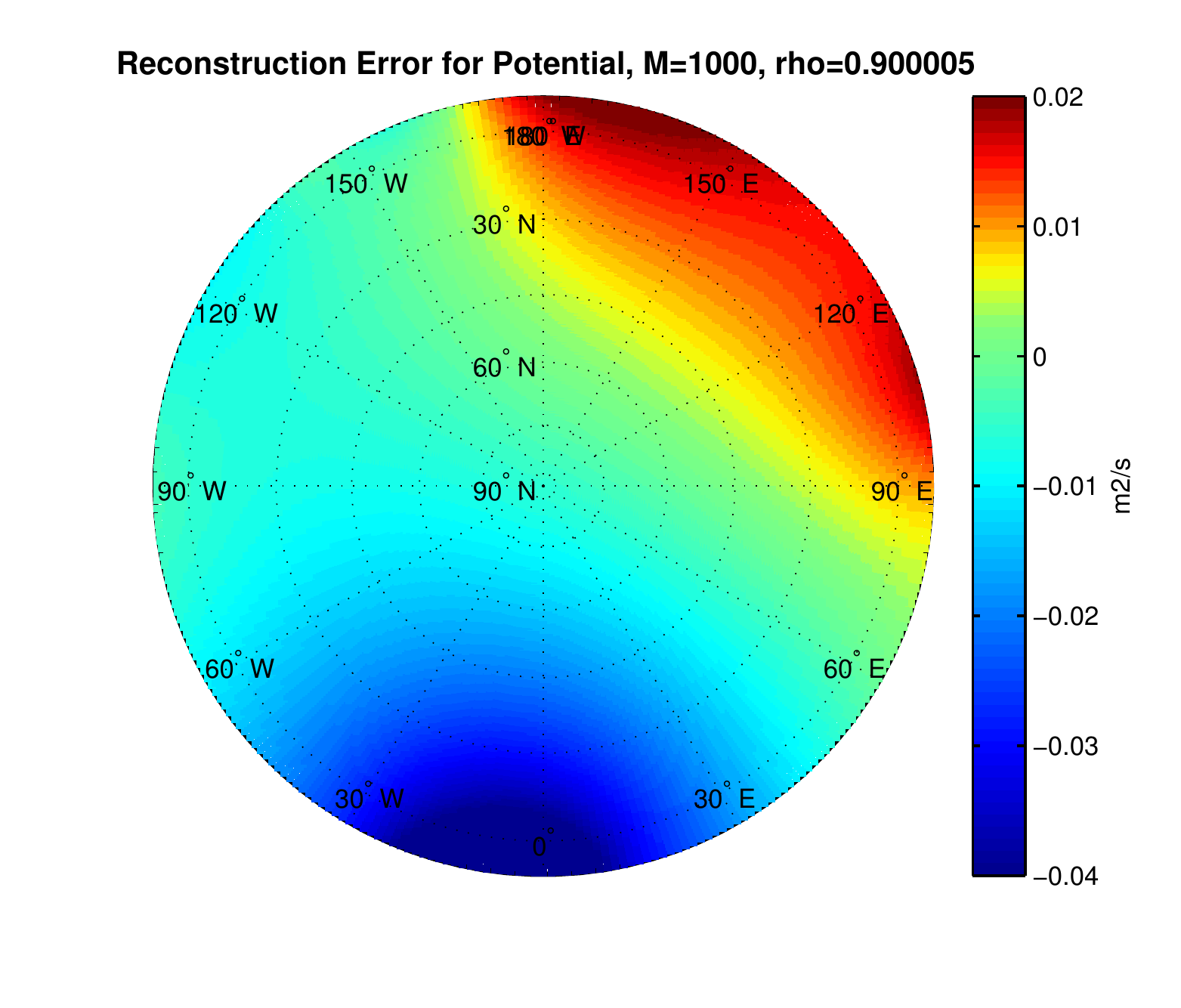}
\\\includegraphics[scale=0.36]{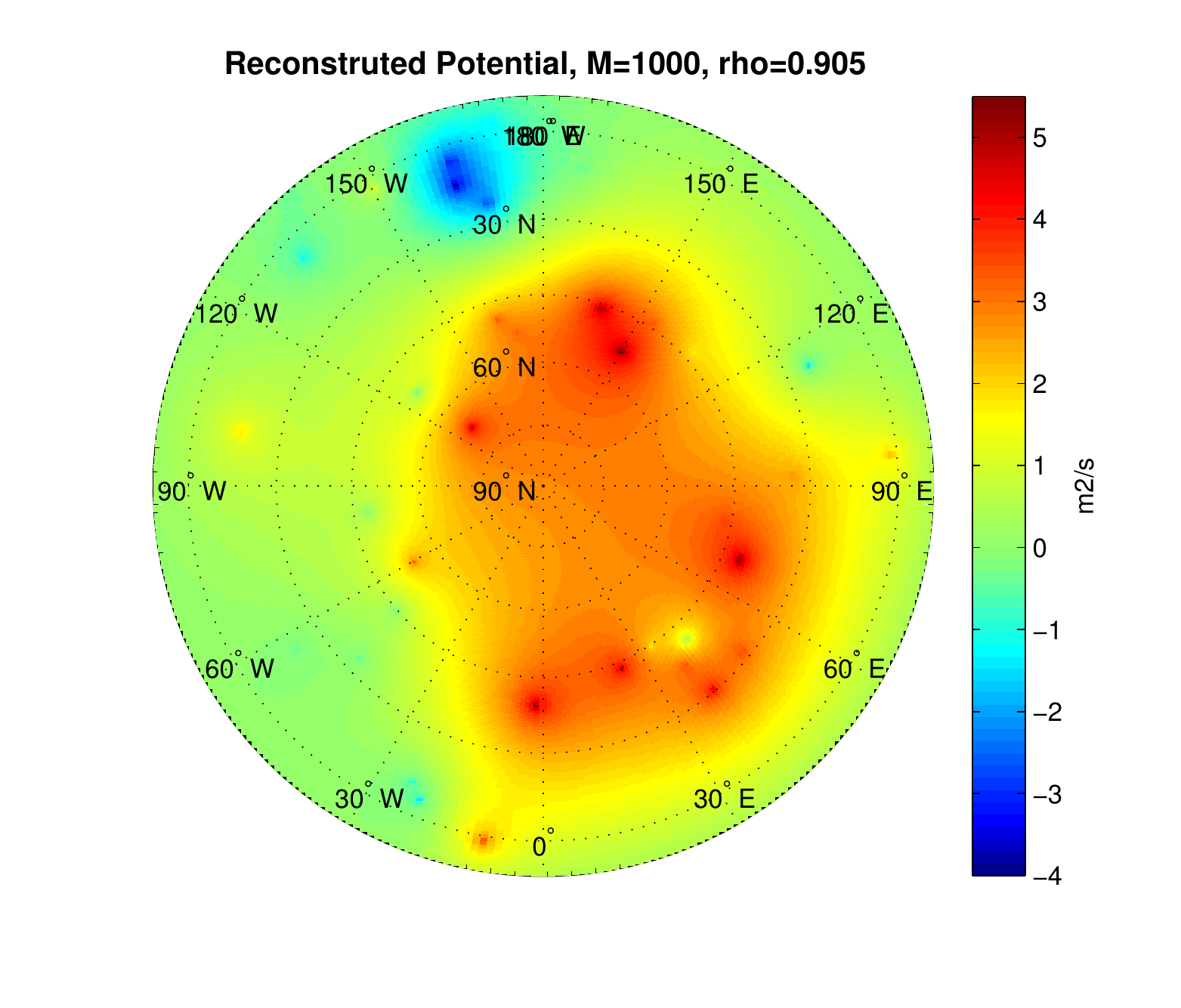}\quad\quad\includegraphics[scale=0.36]{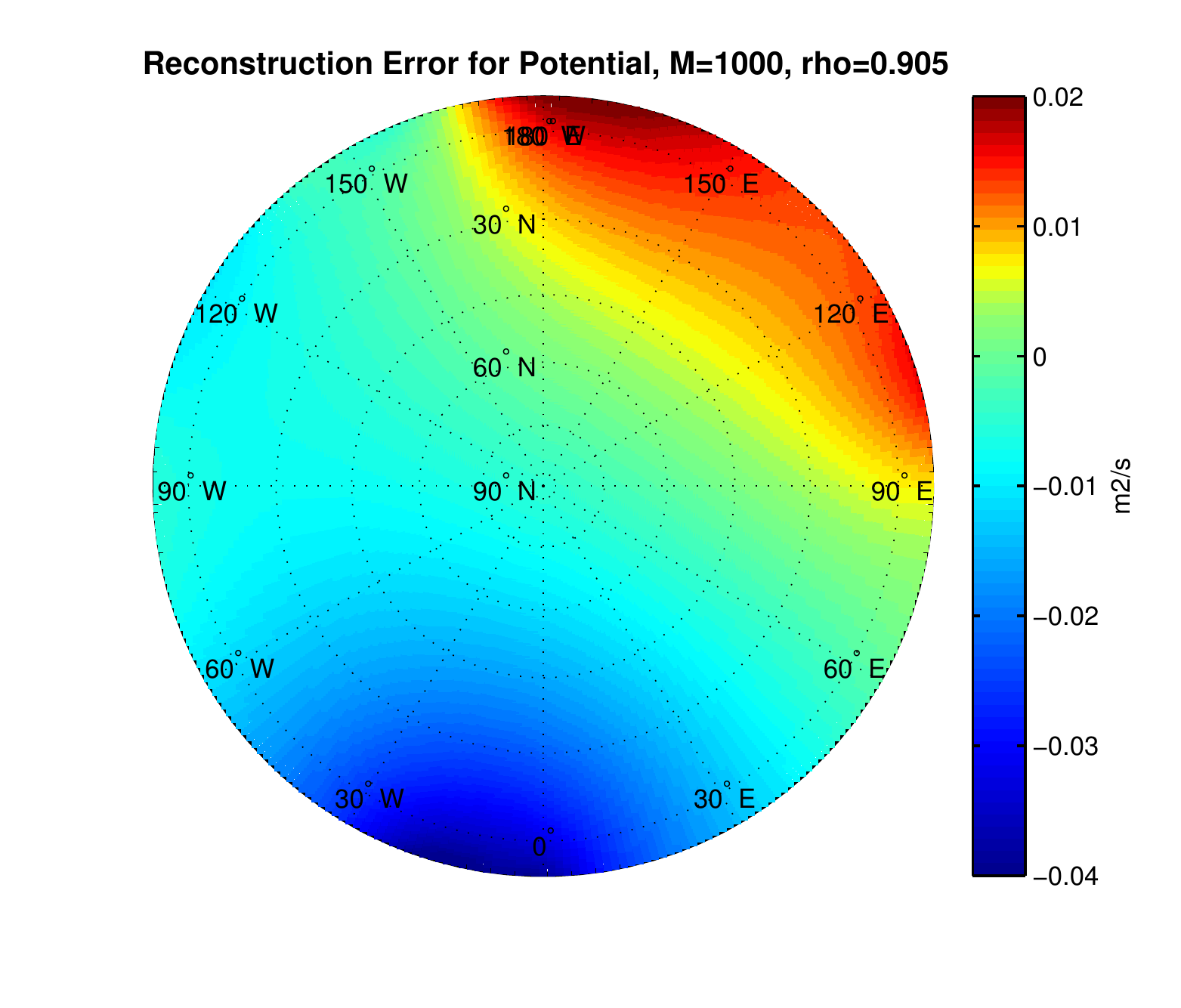}
\\\includegraphics[scale=0.36]{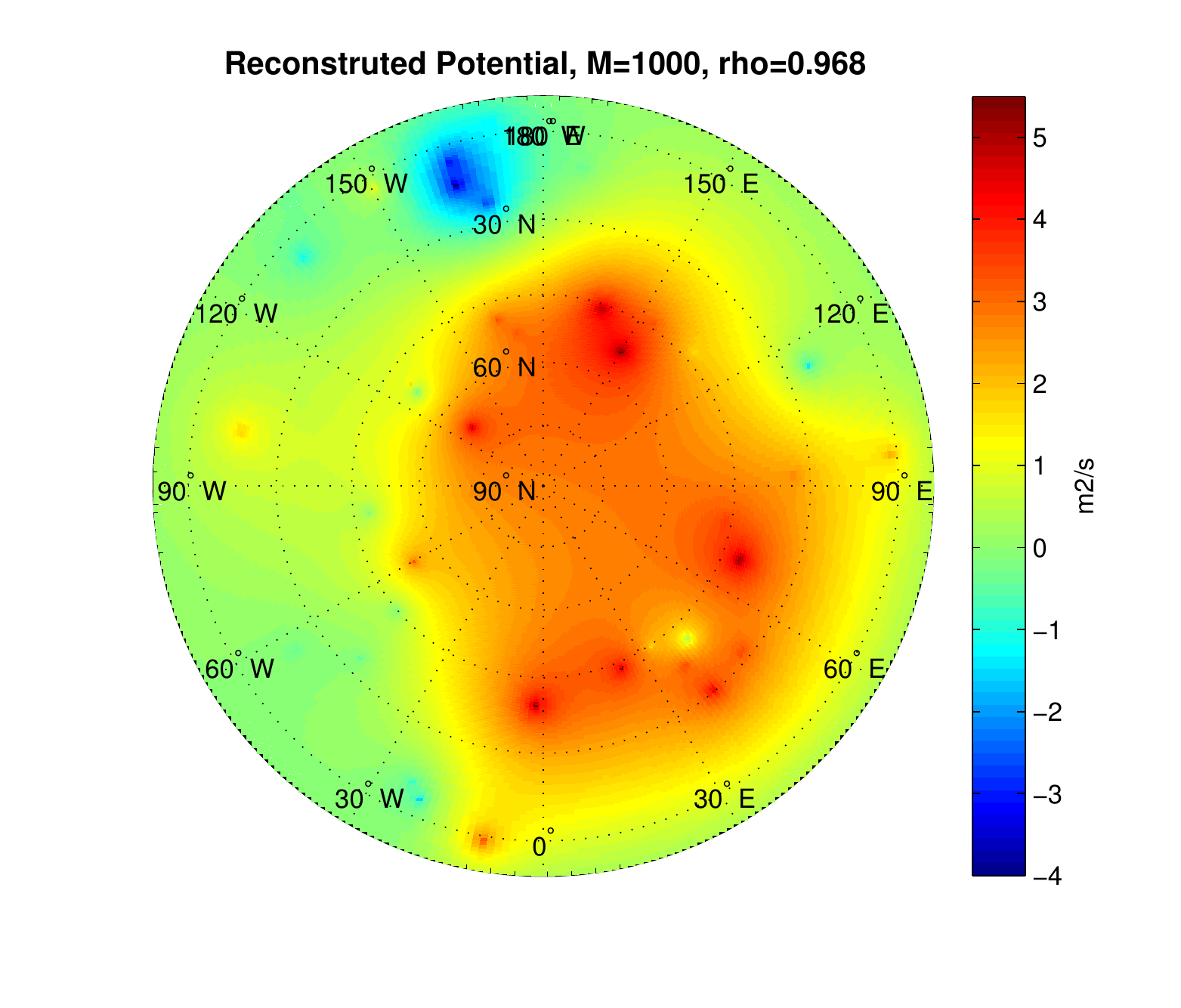}\quad\quad\includegraphics[scale=0.36]{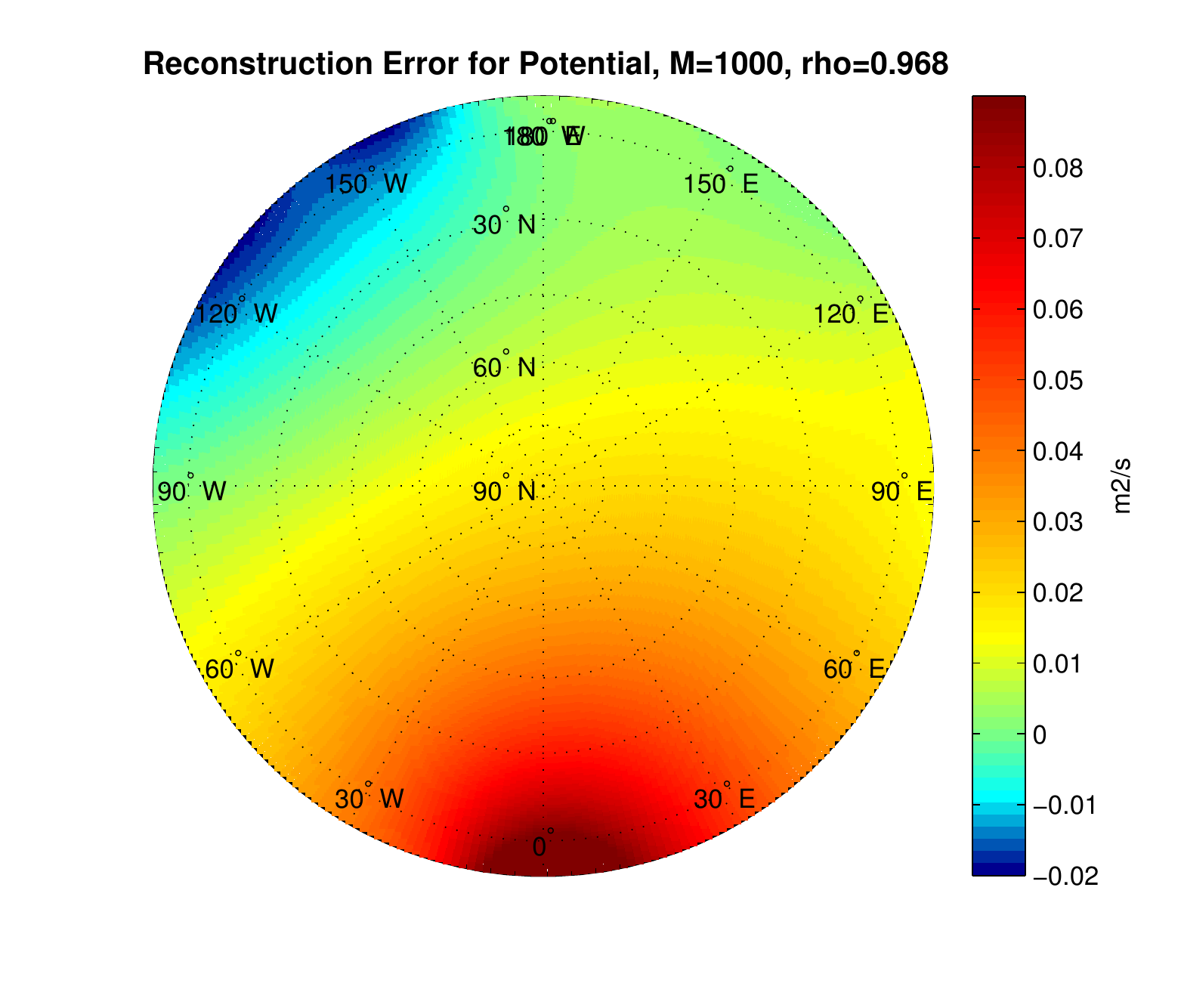}
\includegraphics[scale=0.36]{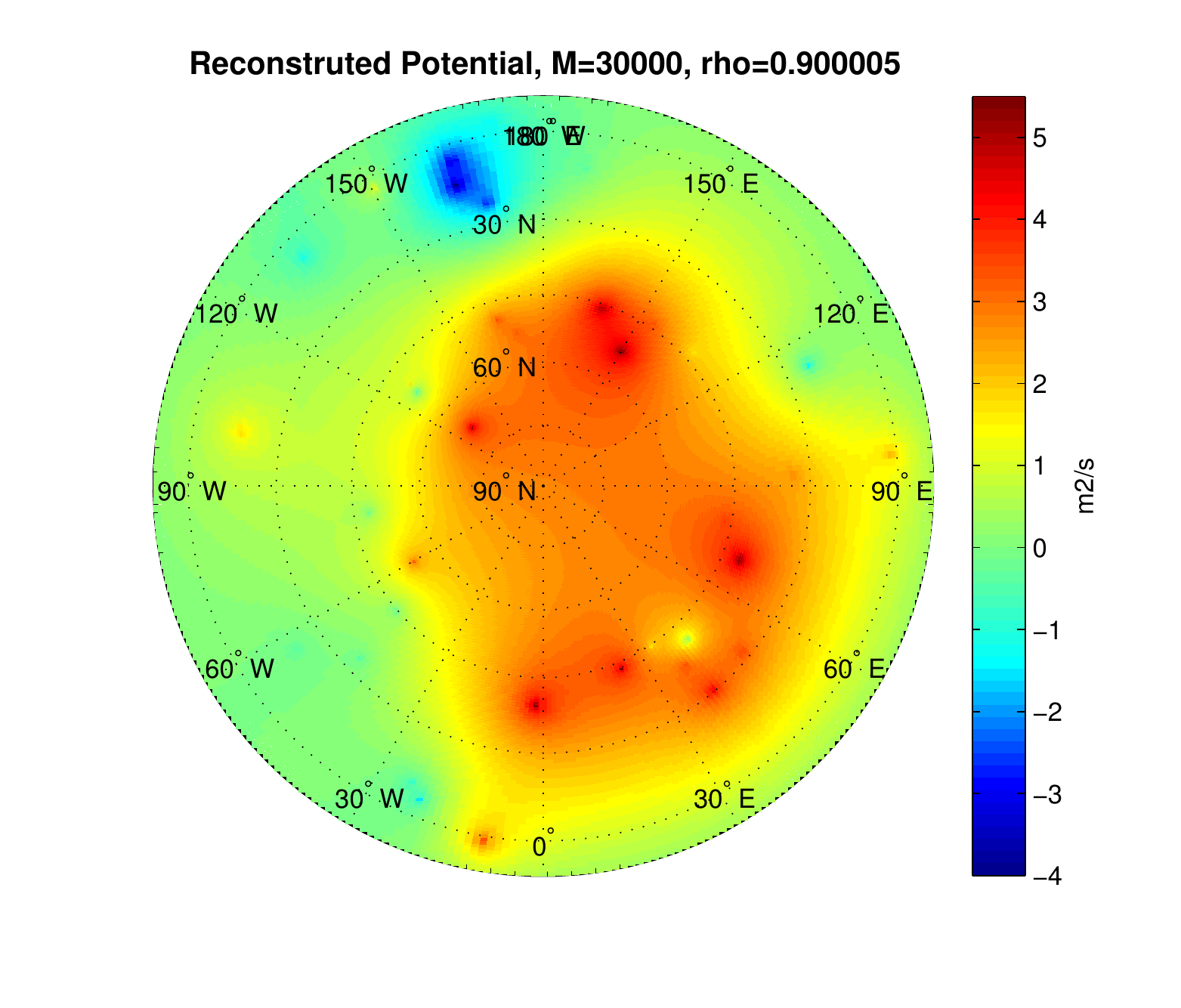}\quad\quad\includegraphics[scale=0.36]{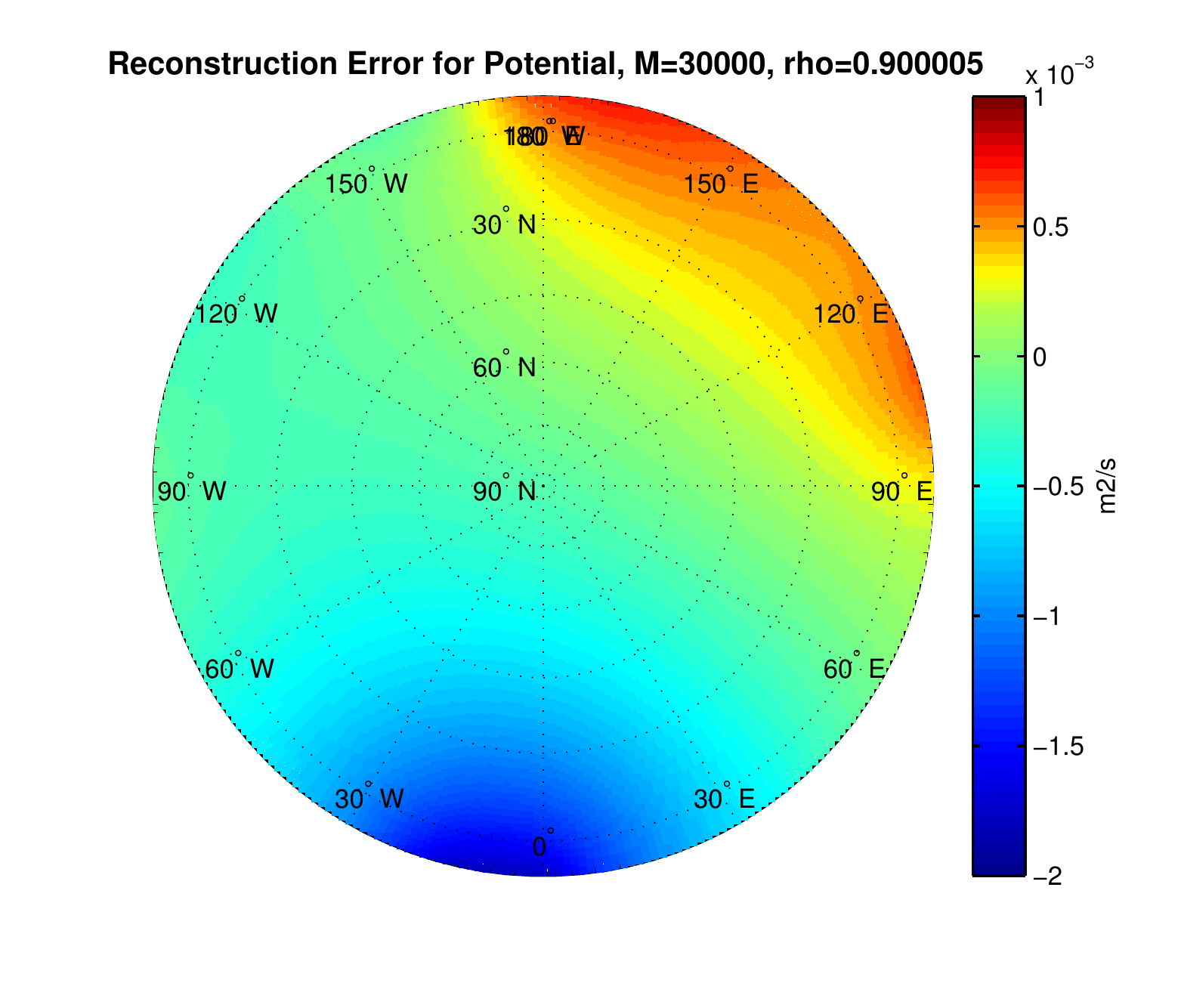}
\end{center}
\caption{The reconstructed potential $\Psi_{M,N,\bar{\rho}}$ (left) and the corresponding reconstruction errors $\Psi-\Psi_{M,N,\bar{\rho}}$ (right) for $M=1000,30000$ and $\bar{\rho}=0.900005, 0.905, 0.968$.}\label{fig:vort1}
\end{figure}

In this section, we focus on the model problem \eqref{eqn:vort1}, \eqref{eqn:vort2}. Opposed to \cite{gemmrich08,kropinski14}, where boundary integral methods have been used, we want to solve it by the method of fundamental solutions based on the results of Section \ref{sec:cfs}. More precisely, we choose $\Gamma$ to be a spherical cap in the Northern hemisphere: $\Gamma=\Gamma_\rho(\zeta)$ with center $\zeta=(0,0,1)^T$ and radius $\rho=0.9$. For simplicity, we set $R=1$. The centers $\eta_i\in\Gamma_\rho(\zeta)$, $i=1,\ldots, N$, and the corresponding strengths $\bar{\omega}_i$ of the point vortices are chosen randomly. The point $\bar{\xi}\in(\Gamma_\rho(\zeta))^c$ from \eqref{eqn:vort2} is set to $\bar{\xi}=(0,0,-1)^T$. Furthermore, we assume the boundary data \eqref{eqn:vort2} to be given in equidistantly distributed points $\xi_i\in\partial\Gamma_\rho(\zeta)$, $i=1,\ldots,M$. Eventually, we interpolate the data by the functions $G_k^{(mod)}$, $k=0,\ldots,M-1$, from Theorem \ref{thm:completefusystemln} and Remark \ref{rem:phikmod}, i.e.,
\begin{align}
 G_k^{(mod)}(\xi)=\frac{1}{4\pi}\ln(1-\xi\cdot\bar{\xi}_k)-\frac{1}{4\pi}\ln(1-\xi\cdot\bar{\xi}),\quad k=1,\ldots,M-1,
\end{align}
where the center points $\bar{\xi}_k$, $k=1,\ldots,M-1$, are chosen to be equidistantly distributed on $\partial\Gamma_{\bar{\rho}}(\zeta)$, for a radius $\bar{\rho}>\rho$. The resulting approximation $\Psi_{M,N,\bar{\rho}}$ of $\Psi$ in $\Gamma_\rho(\zeta)$ is given by 
\begin{align}
 \Psi_{M,N,\bar{\rho}}(\xi)&=\sum_{i=1}^N\bar{\omega}_iG(\Delta^*;\xi\cdot{\eta}_i)-\frac{1}{4\pi}\ln(1-\xi\cdot\bar{\xi})-\tilde{\Psi}_{M,N,\bar{\rho}}(\xi),\quad\xi\in\Gamma_\rho(\zeta),
 \\\tilde{\Psi}_{M,N,\bar{\rho}}(\xi)&=\sum_{k=0}^{M-1}a_kG_k^{(mod)}(\xi),\quad\xi\in\Gamma_\rho(\zeta),
\end{align}
where the coefficients $a_k$, $k=0,\ldots,M-1$, are obtained from the approximate solution of \eqref{eqn:vort1}, \eqref{eqn:vort2} via interpolation of the boundary data. The resulting $\Psi_{M,N,\bar{\rho}}$ and the corresponding reconstruction errors are plotted in Figure \ref{fig:vort1} for different settings of $M,\bar{\rho}$ (we fix the number of point vortices to $N=40$). The actual potential $\Psi$ and the underlying horizontal flow $v$ are shown in Figure \ref{fig:vortflow}. We restrict our test example to a spherical cap $\Gamma=\Gamma_\rho(\zeta)$ because we then know an explicit representation of $\Psi$ via the Dirichlet Green function $G_D(\Delta^*;\cdot,\cdot)$ from Section \ref{sec:cfs} and can compute the reconstruction errors. However, the approach can be easily adapted to more complex geometries of $\Gamma$.

The results in Figure \ref{fig:vort1} show a good performance for the test example of this easy to implement technique. The influence of the parameter $\bar{\rho}$ turns out to be fairly harmless for $M=1000$ source points. A significant deterioration of the reconstruction error does not occur before $\bar{\rho}=0.968$ (cf. Figure \ref{fig:vort1}). However, in general, the method of fundamental solutions can be rather sensitive to the choice of the involved parameters, in particular of the source points $\bar{\xi}_k$ and the collocation points $\xi_k$. Furthermore, it can be advantageous to use a regularized least squares method instead of a simple interpolation.  An overview on the method of fundamental solutions in general and its recent developments can be found, e.g., in \cite{chen08,fairweather98}. Latter, however, treat only the Euclidean setting. The current section is meant as a basic illustration of the method of fundamental solution for boundary value problems intrinsic on the sphere based on the techniques described in this chapter.


\begin{thebibliography}{10}

\bibitem{amm97}
O.~Amm.
\newblock Elementary currents for ionospheric fields.
\newblock {\em J. Geomag. Geoelectr.}, 49:947--955, 1997.

\bibitem{backus96}
G.~Backus, R.~Parker, and C.~Constable.
\newblock {\em Foundations of Geomagnetism}.
\newblock Cambridge University Press, 1996.

\bibitem{baratchart13}
L.~Baratchart, D.P. Hardin, E.A. Lima, E.B. Saff, and B.P. Weiss.
\newblock Characterizing kernels of operators related to thin plate
  magnetizations via generalizations of {H}odge decompositions.
\newblock {\em Inverse Problems}, 29:015004, 2013.

\bibitem{chen08}
C.S. Chen, A.~Karageorghis, and Y.S. Smyrlis.
\newblock {\em The Method of Fundamental Solutions - A Meshless Method}.
\newblock Dynamic Publishers, Inc., 2008.

\bibitem{comblen09}
R.~Comblen, S.~Legrand, E.~Deleersnijdera, and V.~Legata.
\newblock A finite element method for solving the shallow water equations on
  the sphere.
\newblock {\em Ocean Modelling}, 28:12--23, 2009.

\bibitem{dritschel88}
D.G. Dritschel.
\newblock Contour dynamics/surgery on the sphere.
\newblock {\em J. Comp. Phys.}, 78:477--483, 1988.

\bibitem{dritschel15}
D.G. Dritschel and S.~Boatto.
\newblock The motion of point vortices on closed surfaces.
\newblock {\em Proc. R. Soc. A}, 471:20140890, 2015.

\bibitem{duduchava06}
R.L. Duduchava, D.~Mitrea, and M.~Mitrea.
\newblock Differential operators and boundary value problems on hypersurfaces.
\newblock {\em Math. Nachr.}, 279:996--1023, 2006.

\bibitem{fairweather98}
G.~Fairweather and A.~Karageorghis.
\newblock The method of fundamental solutions for elliptic boundary value
  problems.
\newblock {\em Adv. Comp. Math.}, 9:69--95, 1998.

\bibitem{fehlinger09}
T.~Fehlinger, W.~Freeden, C.~Mayer, and M.~Schreiner.
\newblock On the local multiscale determination of the earths disturbing
  potential from discrete deflections of the vertical.
\newblock {\em Comp. Geosc.}, 12:473--490, 2009.

\bibitem{fehlinger07}
T.~Fehlinger, W.~S. Freeden, W.~Freeden, C.~Mayer, D.~Michel, and M.~Schreiner.
\newblock Local modelling of sea surface topography from (geostrophic) ocean
  flow.
\newblock {\em ZAMM}, 87:775--791, 2007.

\bibitem{fenglerfreeden05}
M.J. Fengler and W.~Freeden.
\newblock A nonlinear galerkin scheme involving vector and tensor spherical
  harmonics for solving the incompressible navier-stokes equation on the
  sphere.
\newblock {\em SIAM J. Sci. Comp.}, 27:967--994, 2005.

\bibitem{freeden15}
W.~Freeden.
\newblock Geomathematics: Its role, its aim, and its potential.
\newblock In W.~Freeden, M.Z. Nashed, and T.~Sonar, editors, {\em Handbook of
  Geomathematics}. Springer, 2nd edition, 2015.

\bibitem{freeden09}
W.~Freeden, T.~Fehlinger, M.~Klug, D.~Mathar, and K.~Wolf.
\newblock Classical globally reflected gravity field determination in modern
  locally oriented multiscale framework.
\newblock {\em J. Geod.}, 83:1171--1191, 2009.

\bibitem{freedengerhards12}
W.~Freeden and C.~Gerhards.
\newblock {\em Geomathematically Oriented Potential Theory}.
\newblock Pure and Applied Mathematics. Chapman \& Hall/CRC, 2012.

\bibitem{freeden05}
W.~Freeden, D.~Michel, and V.~Michel.
\newblock Local multiscale approximations of geostrophic flow: Theoretical
  background and aspects of scientific computing.
\newblock {\em Marine Geodesy}, 28:313--329, 2005.

\bibitem{freedenmichel04a}
W.~Freeden and V.~Michel.
\newblock {\em Multiscale Potential Theory (With Applications to Geoscience)}.
\newblock Birkh{\"a}user, 2004.

\bibitem{freedenschreiner09}
W.~Freeden and M.~Schreiner.
\newblock {\em Spherical Functions of Mathematical Geosciences}.
\newblock Springer, 2009.

\bibitem{ganesh11}
M.~Ganesh, Q.T. LeGia, and I.H. Sloan.
\newblock A pseudospectral quadrature method for {N}avier-{S}tokes equations on
  rotating spheres.
\newblock {\em Math. Comp.}, 80:1397--1430, 2011.

\bibitem{gemmrich08}
S.~Gemmrich, N.~Nigam, and O.~Steinbach.
\newblock Boundary integral equations for the {L}aplace-{B}eltrami operator.
\newblock In H.~Munthe-Kaas and B.~Owren, editors, {\em Mathematics and
  Computation, a Contemporary View. Proceedings of the Abel Symposium 2006}.
  Springer, 2008.

\bibitem{gerhards11a}
C.~Gerhards.
\newblock Spherical decompositions in a global and local framework: Theory and
  an application to geomagnetic modeling.
\newblock {\em Int. J. Geomath.}, 1:205--256, 2011.

\bibitem{gerhards12}
C.~Gerhards.
\newblock Locally supported wavelets for the separation of spherical vector
  fields with respect to their sources.
\newblock {\em Int. J. Wavel. Multires. Inf. Process.}, 10:1250034, 2012.

\bibitem{gerhards14c}
C.~Gerhards.
\newblock Multiscale modeling of the geomagnetic field and ionospheric
  currents.
\newblock In W.~Freeden, M.Z. Nashed, and T.~Sonar, editors, {\em Handbook of
  Geomathematics}. Springer, 2nd edition, 2015.

\bibitem{gutkinnewton04}
E.~Gutkin and K.P. Newton.
\newblock The method of images and green's function for spherical domains.
\newblock {\em J. Phys. A: Math. Gen.}, 37:11989--12003, 2004.

\bibitem{heiskanenmoritz67}
W.A. Heiskanen and H.~Moritz.
\newblock {\em Physical Geodesy}.
\newblock W.H. Freeman and Company, 1967.

\bibitem{helms69}
L.L. Helms.
\newblock {\em Introduction to Potential Theory}.
\newblock Wiley-Interscience, 1969.

\bibitem{hesse12}
K.~Hesse and R.S. Womersley.
\newblock Numerical integration with polynomial exactness over a spherical cap.
\newblock {\em Adv. Comp. Math.}, 36:451--483, 2012.

\bibitem{wellenhof05}
B.~Hofmann-Wellenhof and H.~Moritz.
\newblock {\em Physical Geodesy}.
\newblock Springer, 2nd edition, 2005.

\bibitem{ilin91}
A.A. Il'in.
\newblock The {N}avier-{S}tokes and {E}uler equations on two-dimensional closed
  manifolds.
\newblock {\em Math. {USSR} Sb.}, 69:559--579, 1991.

\bibitem{kellogg67}
O.D. Kellogg.
\newblock {\em Foundations of Potential Theory}.
\newblock Springer, reprint edition, 1967.

\bibitem{kidambi98}
R.~Kidambi and K.P. Newton.
\newblock Motion of three point vortices on a sphere.
\newblock {\em Physica D}, 116:143--175, 1998.

\bibitem{kidambinewton00}
R.~Kidambi and K.P. Newton.
\newblock Point vortex motion on a sphere with solid boundaries.
\newblock {\em Phys. Fluids}, 12:581--588, 2000.

\bibitem{kimuraokamoto87}
Y.~Kimura and H.~Okamoto.
\newblock Vortex motion on a sphere.
\newblock {\em J. Phys. Soc. Jpn.}, 56:4203--4206, 1987.

\bibitem{kropinski14}
M.C.A. Kropinski and N.~Nigam.
\newblock Fast integral equation methods for the laplace-beltrami equation on
  the sphere.
\newblock {\em Adv. Comp. Math.}, 40:577--596, 2014.

\bibitem{maximenko09}
N.~Maximenko, P.~Niiler, M.-H. Rio, L.C. Melnichenko, D.~Chambers,
  V.~Zlotnicki, and B.~Galperin.
\newblock Mean dynamic topography of the ocean derived from satellite and
  drifting buoy data using three different techniques.
\newblock {\em J. Atm. Ocean. Tech.}, 26:1910--1919, 2009.

\bibitem{mayermaier06}
C.~Mayer and T.~Maier.
\newblock Separating inner and outer {E}arth's magnetic field from {CHAMP}
  satellite measurements by means of vector scaling functions and wavelets.
\newblock {\em Geophys. J. Int.}, 167:1188–1203, 2006.

\bibitem{mitreataylor99}
M.~Mitrea and M.~Taylor.
\newblock Boundary layer methods for {L}ipschitz domains in {R}iemannian
  manifolds.
\newblock {\em J. Func. Anal.}, 163:181--251, 1999.

\bibitem{mitreataylor00}
M.~Mitrea and M.~Taylor.
\newblock Potential theory on {L}ipschitz domains in {R}iemannian manifolds:
  {S}obolev-{B}esov space results and the {P}oisson problem.
\newblock {\em J. Func. Anal.}, 176:1--79, 2000.

\bibitem{olsen97}
N.~Olsen.
\newblock Ionospheric {F}-region currents at middle and low latitudes estimated
  from {MAGSAT} data.
\newblock {\em J. Geophys. Res.}, 102:4564--4576, 1997.

\bibitem{olsen10b}
N.~Olsen, K-H. Glassmeier, and X.~Jia.
\newblock Separation of the magnetic field into external and internal parts.
\newblock {\em Space Sci. Rev.}, 152:135--157, 2010.

\bibitem{egm2008}
N.K. Pavlis, S.A. Holmes, S.C. Kenyon, and Factor J.K.
\newblock The development and evaluation of the {E}arth {G}ravitational {M}odel
  2008 ({EGM}2008).
\newblock {\em J. Geophys. Res.}, 117:B04406, 2012.

\bibitem{pedlosky79}
J.~Pedlosky.
\newblock {\em Geophysical Fluid Dynamics}.
\newblock Springer, 1979.

\bibitem{polvanidritschel93}
L.M. Polvani and D.G. Dritschel.
\newblock Wave and vortex dynamics on the surface of a sphere.
\newblock {\em J. Fluid Mech.}, 255:35--64, 1993.

\bibitem{stewart}
R.H. Stewart.
\newblock Introduction to physical oceanography.
\newblock Online.

\bibitem{wermer74}
J.~Wermer.
\newblock {\em Potential Theory}.
\newblock Springer, 1974.

\end{thebibliography}
\end{document}